\documentclass[11pt,a4paper]{amsart}
\usepackage{amsmath, amssymb, a4wide, epsfig, enumerate, psfrag}
\usepackage[utf8]{inputenc}
\usepackage[english]{babel}
\usepackage[active]{srcltx}
\usepackage{color}
\usepackage{mathtools}
\usepackage{newtxmath}

\usepackage{caption}
\usepackage{tikz}
\usepackage{tkz-euclide}
\usetikzlibrary{calc,angles,through,intersections,quotes}
\usepackage{graphicx} 

\date{\today}

\author{Bertrand Deroin \and Nicolas Tholozan}
\thanks{B.D.'s research was supported by the CY Initiative under the Emergence programme (ANR-16-IDEX-0008, project IsoFol). Both authors' research is supported by ANR-25-CE40-1360 (project IsoMoDyn).}
\address{CNRS \\ Laboratoire AGM \\
CY Cergy Paris Universit\'e \\
2 av. Adolphe Chauvin \\
95302 Cergy-Pontoise Cedex}
\email{bertrand.deroin@cyu.fr}

\address{CNRS \\ Laboratoire de Math\'ematiques d'Orsay  \\ Universit\'e Paris-Saclay \\ 307, rue Michel Magat \\ F-91405 Orsay Cedex}
\email{nicolas.tholozan@cnrs.fr}

\title{Holonomies of hyperbolic surfaces with one branch point}

\newcommand{\R}{\mathbb{R}}
\newcommand{\Z}{\mathbb{Z}}

\renewcommand{\d}{\mathrm{d}}

\newcommand{\norm}[1]{\left\Vert#1\right\Vert}

\newcommand{\D}{\mathbb{D}}

\newcommand{\SL}{\mathrm{SL}}

\newcommand{\C}{\mathbb{C}}
\newcommand{\N}{\mathbb{N}}

\renewcommand{\H}{\mathbb{H}}

\newcommand{\Isom}{\mathrm{Isom}}
\newcommand{\Aut}{\mathrm{Aut}}
\newcommand{\Hom}{\mathrm{Hom}}
\newcommand{\Rep}{\mathfrak{X}}

\newcommand{\Id}{\mathrm{Id}}

\newcommand{\Tr}{\mathrm{Tr}}
\newcommand{\Ad}{\mathrm{Ad}}

\newcommand{\INT}{\mathrm{Int}}

\newcommand{\function}[5]{
\begin{array}{rrcl}
#1 : & #2 & \to & #3 \\
\ & #4 & \mapsto & #5
\end{array}
}

\newcommand{\dev}{\mathrm{dev}}

\renewcommand{\epsilon}{\varepsilon}

\newcommand{\PSL}{\mathrm{PSL}}

\newcommand{\PU}{\mathrm{PU}}

\renewcommand{\phi}{\varphi}

\DeclareMathOperator{\vol}{vol}

\DeclareMathOperator{\im}{Im}

\newcommand{\Vol}{\mathrm{Vol}}

\renewcommand{\tilde}[1]{\widetilde{#1}}

\newcommand{\euler}{\mathbf{eu}}
\newcommand{\Tol}{\mathbf{tol}}
\newcommand{\Teich}{\mathcal{T}}
\newcommand{\MCG}{\mathrm{MCG}}

\newcommand{\EE}{\mathbf{E}}

\theoremstyle{plain}
\newtheorem{prop}{Proposition}[section]
\newtheorem{thm}[prop] {Theorem}
\newtheorem{lem}[prop] {Lemma}
\newtheorem{cor}[prop]{Corollary}

\newtheorem{fact}[prop]{Fact}

\theoremstyle{definition}
\newtheorem{defi}[prop]{Definition}

\theoremstyle{remark}

\newtheorem{rmk}[prop]{Remark}

\newtheoremstyle{mainThm}
{0.3cm}       
{0.3cm}       
{\itshape}    
{}            
{\scshape}    
{.}           
{0.5em}       
{}

\theoremstyle{mainThm}
\newtheorem*{BigThm}{Main Theorem}
\newtheorem{conj}{Conjecture}

\begin{document}

\maketitle

\begin{abstract}
    Let $S$ be a closed oriented surface of genus $g\geq 2$. We prove that any representation of $\pi_1(S)$ into $\PSL(2,\R)$ with Euler class $3-2g$ is the holonomy of a hyperbolic metric on $S$ with one simple branch point.
\end{abstract}

\section*{Introduction}

Let $S$ be a closed oriented surface of genus $g \geq 2$. A \emph{branched hyperbolic metric} on $S$ is a singular metric $\sigma$ which is locally isometric to the pull-back of the Poincar\'e metric $\sigma_\D$ on the unit disc $\D$ under the map $z\mapsto z^k$, for some positive integer $k$. A point with $k\geq 2$ is a \emph{branch point} with multiplicity $k-1$ and cone angle $2k\pi$. Let $\H$ denote the hyperbolic plane and $\Isom_+(\H)$ its group of orientation-preserving isometries. A branched hyperbolic metric determines a complex structure $J$ on $S$ and a \emph{holonomy representation} $\rho: \pi_1(S) \to \Isom_+(\H)$, such that the lift of $\sigma$ to the universal cover $\tilde S$ is the pull-back of $\sigma_\H$ under a non-constant $\rho$-equivariant holomorphic map $\dev: \tilde S \to \H$ called a \emph{developing map}. 

A longstanding question in low-dimensional geometry, popularised by Goldman (see for instance \cite{GoldmanConjectures}), asks which surface group representations arise as holonomies of branched hyperbolic metrics. To make this question more precise, recall that one can associate to every representation  $\rho: \pi_1(S)\to \Isom_+(\H)$ its \emph{Euler class} $\euler(\rho)$, an integral invariant taking all values between $\chi(S)= 2-2g$ and $-\chi(S) = 2g-2$ by the \emph{Milnor--Wood inequality}. Here, we will find it more convenient to use the \emph{volume} of a representation $\rho$, given by
\[\Vol(\rho) = -2\pi \euler(\rho) = \int_S f^*\vol_\H~,\]
where $\vol_\H$ is the hyperbolic volume form on $\H$ and $f:\tilde S \to \H$ is any smooth $\rho$-equivariant map. When $\rho$ is the holonomy of a branched hyperbolic metric, its volume equals the area of that metric and is therefore positive. The Gauss--Bonnet formula yields
\[\Vol(\rho) = 2\pi(\vert \chi(S)\vert - k)~,\]
where $k$ is the number of branch points counted with multiplicity.

Conversely, Thurston and Goldman \cite{GoldmanThese} proved that every representation of maximal volume $2\pi \vert  \chi(S)\vert$ is the holonomy of a (smooth) hyperbolic metric on $S$, which is unique up to isotopy.
In contrast, Tan \cite{Tan} constructed examples of representations with positive volume that are not holonomies of branched hyperbolic metrics. His examples are quite specific: they are representations $\rho$ with discrete image and non-trivial kernel. This led Goldman to conjecture that every representation with positive volume and dense image is the holonomy of a branched hyperbolic metric \cite[Conjecture 3.9]{GoldmanConjectures}. In fact, by analogy with known results on holonomies of branched euclidean metrics \cite{Haupt20,Kapovich20,Ghazouani18}, one may formulate the following more precise conjecture:

\begin{conj}[strong form of Goldman's conjecture] \label{conj: Goldman Conjecture}
Let $\rho:\pi_1(S) \to \Isom_+(\H)$ be a representation with positive volume. Then $\rho$ is the holonomy of a branched hyperbolic metric unless it is unfaithful, has discrete image, and satisfies
\[\Vol(\rho) = \Vol(\rho(\pi_1(S))\backslash \H)~.\]
\end{conj}

This conjecture predicts, in particular, that \emph{all} representations of volume $2\pi(2g-3)$ are geometrisable. The main purpose of this article is to prove this statement. 

\begin{BigThm}
Every representation $\rho:\pi_1(S)\to \Isom_+(\H)$ of volume $2\pi(2g-3)$ is the holonomy of a hyperbolic metric on $S$ with one branch point of angle $4\pi$.
\end{BigThm}

This result was mentioned in \cite{Tan} as a work in preparation by Goldman and Neumann. However, Goldman later stated it as a conjecture ``arising in discussions with Neumann'' \cite[Conjecture 3.9]{GoldmanConjectures}.
With our convention $\Vol(\rho)=-2\pi\euler(\rho)$, the theorem concerns Euler class $3-2g$. Note that Tan's counterexamples have volume at most $2\pi(2g-4)$.

\subsection{Three approaches}

One can think of three rather different approaches to Goldman's conjecture.

The first approach is geometric. One might try to ``build'' a hyperbolic surface with prescribed holonomy by gluing pieces in a suitable way. This is, in essence, how Goldman proves that representations of maximal volume are Fuchsian: he shows that their restrictions to the pairs of pants of some decomposition of $S$ are holonomies of hyperbolic metrics with geodesic boundaries, then glues these pairs of pants appropriately. Unfortunately, one cannot expect representations of lower volume to behave well when restricted to all the pieces of a given pair-of-pants decomposition.

The second approach is analytic. For a representation $\rho$ of positive volume and a complex structure $J$ on $S$, there exists a unique $\rho$-equivariant \emph{harmonic map} $f_{J,\rho}: (\tilde{S},J)\to \H$. Its energy $\EE_\rho(J)$ defines a smooth functional on the Teichm\"uller space of $S$. One can prove that $f_{J,\rho}$ is holomorphic if and only if $J$ is a critical point of the functional $\EE_\rho$; such a critical point is necessarily a global minimum, see Theorem \ref{t: Critical points energy}. Thus constructing a branched hyperbolic metric with holonomy $\rho$ amounts to finding a minimum of $\EE_\rho$. Unfortunately, as we will see in Proposition \ref{prop: Energy not proper}, the functional $\EE_\rho$ is never proper when $\rho$ is not maximal.

The third approach is dynamical. The locus of conjugacy classes of holonomies of branched hyperbolic metrics is open (Lemma \ref{l:ehresman-thurston}) and invariant under the mapping class group action. One can therefore try to geometrize a representation by moving it sufficiently close to a known geometric example via a (possibly large) mapping class. Eventually, a sufficiently precise description of orbit closures under mapping class group actions could lead to a complete resolution of Conjecture \ref{conj: Goldman Conjecture}, as it does for holonomies of translation surfaces \cite{Kapovich20} and branched affine structures \cite{Ghazouani18}. Unfortunately, the dynamical properties of the mapping class group action on character varieties is still elusive. Goldman conjectured that the mapping class group action on non-Fuchsian components is ergodic \cite[Conjecture 3.1]{GoldmanConjectures}. This would imply that almost every representation of positive volume is geometrizable. Unfortunately, Goldman's ergodicity conjecture has only been resolved in genus $2$ by March\'e--Wolff \cite{MarcheWolff,MarcheWolff2}.

Our proof of the main theorem somehow combines these three approaches. We begin with a minimizing sequence for the energy functional $\EE_\rho$. Although this sequence may diverge in $\Teich(S)$, its projection modulo the mapping class group accumulates to a stable, possibly nodal curve in the Deligne--Mumford compactification of the moduli space.
We prove that, on each connected component of the smooth locus of the limiting curve, suitably normalized harmonic
maps converge to weakly conformal maps
(Lemma~\ref{lem: Conformal limits harmonic maps minimizing sequence}).
The condition $\Vol(\rho) = 2\pi (2g-3)$ severely restricts the possible degenerations:
either every limiting map is non-constant and holomorphic, or precisely
one component of Euler characteristic $-1$ collapses
(Lemma~\ref{lem: decomposition co-euler class 1}).
The main difficulty is then to reconstruct a metric with the original
holonomy from these limiting data.

For this purpose, we introduce in Section \ref{ss: Hyperbolic blocks} what we call \emph{hyperbolic blocks}: compact surfaces carrying branched hyperbolic metrics whose boundary components are regular curves of constant geodesic curvature, subject to additional technical conditions.
In the limit, the non-collapsed components carry singular hyperbolic metrics that can be truncated to produce hyperbolic blocks (Proposition~\ref{prop: Cone metrics -> hyperbolic blocks}). A version of the Ehresmann--Thurston principle
(Lemma~\ref{l:perturbation of hyperbolic block}) allows us to freely deform the holonomy of a hyperbolic block and produce hyperbolic blocks with the required holonomy on the non-collapsed pieces, with  control of their boundary lengths. 

The collapsing components require explicit constructions on pairs of pants (Lemma~\ref{l: geometrization of representations with small translation length}) and one-holed tori (Lemmas~\ref{l:Uniformization punctured torus elliptic}
and~\ref{l: geometrization pairs of pants in the parabolic case}). These constructions realise the relevant representations as holonomies of hyperbolic blocks with sufficiently short boundary components. Finally, hyperbolic blocks are designed so that they can be assembled together under some conditions on boundary lengths that we need to control throughout the construction (Lemma \ref{l: gluing lemma}). Because of constant-curvature boundary condition, there is a twist parameter in the gluing, which can be adjusted to recover the prescribed global holonomy. This produces a branched hyperbolic metric on $S$ with holonomy $\rho$. The Gauss--Bonnet formula and the prescribed volume imply that it has exactly one simple branch point.

\subsection{Notes and references}

The study of holonomies of possibly branched locally homogeneous structures on surfaces has a long history. One of the earliest instances of the problem is the determination of holonomies of \emph{translation surfaces}: which classes in $\mathrm{H}^1(S,\C)$ are represented by a holomorphic $1$-form \emph{for some complex structure} on $S$. An old theorem of Haupt \cite{Haupt20} gives a complete answer analogous to Conjecture \ref{conj: Goldman Conjecture}. Kapovich \cite{Kapovich20} revisited this theorem using mapping class group dynamics: in this setting, the mapping class group action factors to the action of $\mathrm{Sp}(2g,\Z)$ on $\C^{2g}$ and Ratner's theory describes orbit closures, leading to the classification of holonomies of translation surfaces (see also \cite[Section 6]{CDF23}). Building on this approach, Ghazouani determined the holonomies of branched complex affine structures on surfaces.

Another fundamental result is the theorem of  Gallo--Kapovich--Marden \cite{GKM}. They proved that every non-elementary representation $\rho:\pi_1(S)\to \PSL(2,\C)$ is the holonomy of a \emph{complex projective structure} on $S$ with either $0$ or $1$ branch point, according to whether $\rho$ lifts to $\SL(2,\C)$. In genus $2$, their argument had a gap, which was subsequently filled by Le Fils \cite{LeFils23}. In yet another direction, Faraco--Rungi \cite{FaracoRungi} recently investigated holonomies of \emph{branched real projective structures}.

Interestingly, the Gallo--Kapovich--Marden theorem applies to representations in $\PSL(2,\R)$. In particular, it implies that every representation $\rho:\pi_1(S)\to \Isom_+(\H)$ of Euler class $3-2g$ is the holonomy of a simply branched complex projective structure. This does not prove our main theorem, however, because the developing map of this complex projective structure will typically map some pieces of the surface to the lower half-plane $\{z\in \C\mid \Im(z) < 0\}$. Using some additivity property of the volume of representations (see Proposition \ref{prop: Volume Additive}), one could prove that these negative pieces are either cylinders or a subsurface of Euler characteristic $-1$ containing the branch point. One might hope to ``erase'' these negative regions by a form of ``reverse grafting'' generalizing \cite{GoldmanProjectiveFuchsianHolonomy}. Unfortunately, the Gallo--Kapovich--Marden construction does not provide the geometric controls needed to perform such a surgery.\\

Branched hyperbolic metrics are special cases of hyperbolic metrics with cone singularities. The locol holonomy around a cone singularity whose angle is not a multiple of $2\pi$ is non-trivial. The question of which representations of the fundamental group of a punctured surface arise as holonomies of hyperbolic cone metrics with cone singularities has also been studied from several perspectives.

The uniformisation theorems of McOwen and Troyanov
\cite{McOwen88,Troyanov} identify the deformation space of hyperbolic metrics with prescribed singularities to the Teichm\"uller space of the punctured surface, under some Gauss--Bonnet condition. Luo
\cite{Luo93} and Mondello \cite{Mondello} studied the local behaviour of the holonomy map from this deformation space to the character variety. For the one-holed torus, Mathews \cite{Mathews11} characterised the representations geometrised by hyperbolic metrics with piecewise geodesic boundary and at most one corner.

We expect that the strategy of this paper can also be used to characterise the holonomies of hyperbolic metrics with cone singularities of angle in $(0,2\pi)$. In forthcoming work, Mondello and the second author resolve this question when the sum of any two angles is at most $2\pi$. In order to further investigate this problem, we introduced in \cite{DeroinTholozan19} a notion of \emph{relative Euler class} and established a version of the Milnor--Wood inequality. Interestingly, we uncovered the existence of representations of the fundamental group of a punctured sphere that violate the expected bound, which we called \emph{supra-maximal}. Their existence is in some sense an obstruction to a naive approach to the geometrization problem. Nevertheless, we proved in
\cite[Theorem~5 and Section~4]{DeroinTholozan19}
that every nonabelian supra-maximal representation is, after conjugation by an orientation-reversing isometry, the holonomy of a hyperbolic metric with conical singularities at the punctures
and additional branch points. Remarkably, this geometrisation exists in every conformal structure:
the associated finite-energy equivariant harmonic map to the lower half-plane is always holomorphic.
Mondello recovered this geometrisation using parabolic Higgs
bundles \cite[Corollary~4.17]{Mondello16}.
More recently, Fenyes and Maret \cite{FenyesMaret}
gave an explicit geometric construction of hyperbolic cone
spheres geometrising these representations, akin to some of the constructions presented in the present paper.

\subsection*{Acknowledgements}
We began working on this project more than ten years ago, and our results were already announced at a conference in honour of Fran\c cois Labourie in 2022. We thank our mathematical community for their patience.

\subsection*{Tool Use}
All the original mathematical ideas in this paper are our own and predate the recent progress of artificial intelligence in mathematical research. OpenAI's Astra was helpful in drawing the pictures of the paper and correcting some grammar and typos. 

\setcounter{tocdepth}{1}
\tableofcontents

\section{Preliminaries on surface topology}

This section recalls some standard ingredients of surface topology, mainly to fix notation and conventions.

\subsection{Orientation}

\begin{figure}[h]
\centerline{\includegraphics[scale=0.6]{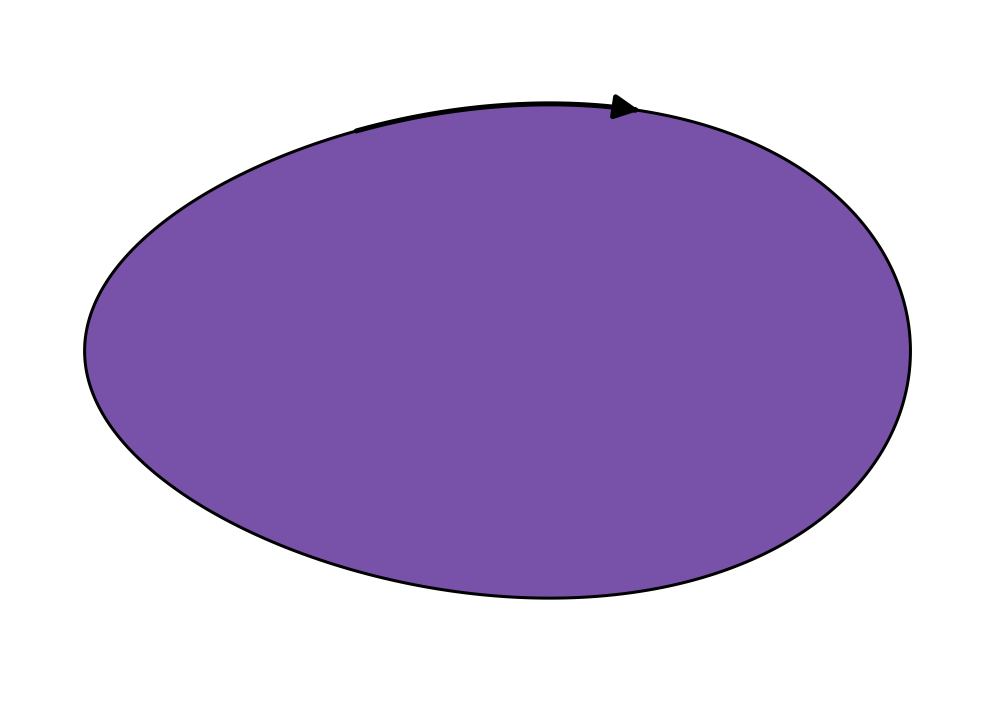}}
\caption{Orientation convention.}
\label{f: orientation}
\end{figure}

Throughout the paper, the surfaces we consider will be oriented, possibly with boundary. When $S$ is a surface with boundary, recall that an orientation of $S$ induces an orientation of its boundary. Our convention for the orientation of $\partial S$ is unusual: the surface lies on the right of \(\partial S\), as shown in Figure \ref{f: orientation}.

\subsection{Fundamental group and universal cover}

Whenever we are given a connected surface $S$, we implicitly fix a universal covering $\tilde S \to S$. We can thus define the fundamental group $\pi_1(S)$ of $S$ as the automorphism group of its universal cover (without having to refer to a basepoint).

Given a map $\phi: S' \to S$ between two surfaces, we also implicitly fix a lift $\tilde \phi : \tilde S' \to \tilde S$, equivariant with respect to a morphism $\phi_*: \pi_1(S')\to \pi_1(S)$. 

Let $S'$ be a connected subsurface of $S$ and denote by $\iota: S'\to S$ the inclusion. The subsurface $S'$ is called \emph{incompressible} if the lift $\tilde \iota : \tilde S' \to \tilde S$ of the inclusion is injective, or equivalently if the induced morphism $\iota_*: \pi_1(S') \to \pi_1(S)$ is injective. In that case, we identify $\tilde S'$ with its image under $\tilde \iota$ and consider it as a subsurface of $\tilde S$. Similarly, we identify $\pi_1(S')$ with its image under $\iota_*$, namely
\[\{ \gamma \in \pi_1(S) \mid \gamma (\tilde S') = \tilde S'\}~.\]

In particular, if $\rho$ is a representation $\pi_1(S) \to \Isom_+(\H)$, we write $\rho_{\vert S'} \coloneqq \rho \circ \iota_*$ for its restriction to $\pi_1(S')$.

Similarly, if we are given an essential oriented simple closed curve $b \subset S$, we implicitly assume that we are also given an embedding $\tilde b \subset \tilde S$. The curve $b$ (together with its lift) thus represents an element of the fundamental group of $S$ that we denote by $[b]$. If $\rho$ is a representation of $\pi_1(S)$ into $\Isom_+(\H)$, we write $\rho(b)$ for $\rho([b])$. 

If $S$ has $I$ connected components $S_i, 1\leq i \leq I$, we define 
\[\pi_1(S) \coloneqq \bigsqcup_{1\leq i \leq I} \pi_1(S_i)~,\]
which is now only a groupoid. We will use this terminology in the following way:
\begin{itemize}
    \item A \emph{representation of $\pi_1(S)$}
 into $\Isom_+(\H)$ is just a collection $(\rho_i)_{1\leq i \leq I}$ where $\rho_i$ is a representation of $\pi_1(S_i)$,
    \item In particular, if $S'=\bigsqcup_{i=1}^I S_i'$ is a disconnected
subsurface of $S$ with incompressible components, and
$\rho:\pi_1(S)\to\Isom_+(\H)$ is a representation, we write
\[\rho_{\vert S'}=(\rho_{\vert S_i'})_{1\leq i\leq I}.\]
 \end{itemize}

\subsection{Topological classification}

A connected oriented surface (with or without boundary) is said to be \emph{of finite topological type} if $\pi_1(S)$ is finitely generated.  An oriented surface without boundary is of finite type if and only if it is homeomorphic to the interior of a compact oriented surface of genus $g\geq 0$ with $K \geq 0$ boundary components or, equivalently, to the complement of $K$ points in a closed surface of genus $g$.

The Euler characteristic of a compact surface $S$ of genus $g$ with $K$ boundary components (or its interior) is
\[\chi(S) = 2-2g-K~.\]
It is negative unless $S$ is a closed sphere, a disc, an annulus or a closed torus.

Given a compact oriented surface $S$ with boundary components $(c_k)_{1\leq k \leq K}$, we denote by $S^\bullet$ the surface obtained by collapsing each boundary component to a point. It is a closed oriented surface with a continuous surjective map $\pi:S \to S^\bullet$ and a finite collection of marked points $c^\bullet_k, 1\leq k \leq K$.

Conversely, if $S$ is a closed oriented surface and $P = \{p_1,\ldots p_K\} \subset S$ a finite set of marked points, we set
\[(S\setminus P)^\circ = S\setminus  \bigcup_{1\leq k \leq K} D_k\]
where the $D_k$ are pairwise disjoint open discs containing $p_k$. The boundary of $D_k$ is a boundary component of $(S\setminus P)^\circ$, which we denote by $p_k^\circ$.

\subsection{Decompositions of a surface $S$}

\label{ss:Decompositions surface}

Let $S$ be a compact oriented surface, possibly disconnected, possibly with boundary. We denote its boundary components by $(c_k)_{1\leq k \leq K}$. Recall that these boundary curves have a natural orientation compatible with the orientation of $S$. We assume that each connected component of $S$ has negative Euler characteristic.

Recall that a \emph{multicurve} on $S$ is a finite collection $(b_l)_{1\leq l \leq L}$ of oriented disjoint smooth simple closed curves. It is \emph{essential} if the following conditions hold:
\begin{enumerate}
\item for all $l\leq L$, the curve $b_l$ is not homotopically trivial,
\item for all $k\leq K$ and all $l \leq L$ the curve $b_l$ is homotopic neither to $c_k$ nor to $c_k^{-1}$,
\item for all $l\neq l' \leq L$, the curve $b_l$ is homotopic neither to $b_{l'}$ nor to $b_{l'}^{-1}$.
\end{enumerate}

\begin{defi}\label{d: decomposition}
    A \emph{decomposition} of $S$ is a partition of $S$ as
    \[S = \bigsqcup_{1\leq i \leq I} S_i \sqcup \bigsqcup_{1\leq l \leq L} b_l\]
    where the $S_i$ are connected incompressible subsurfaces and $(b_l)_{1\leq l \leq L}$ is an essential multicurve.
\end{defi}

Though a decomposition is really just given by an essential multicurve (since the $S_i$ are just the connected components of the complement), for practical purposes, it will be convenient to include in these data the labeling of both the curves along which we cut and the components of the complement.

\begin{rmk}\label{r: essential curve}
    The fact that the multicurve $(b_l)_{1\leq l \leq L}$ is essential is equivalent to the fact that the connected components $S_i$ have negative Euler characteristic.
\end{rmk}

\begin{figure}[htbp]
\centering
\includegraphics[width=0.98\textwidth]{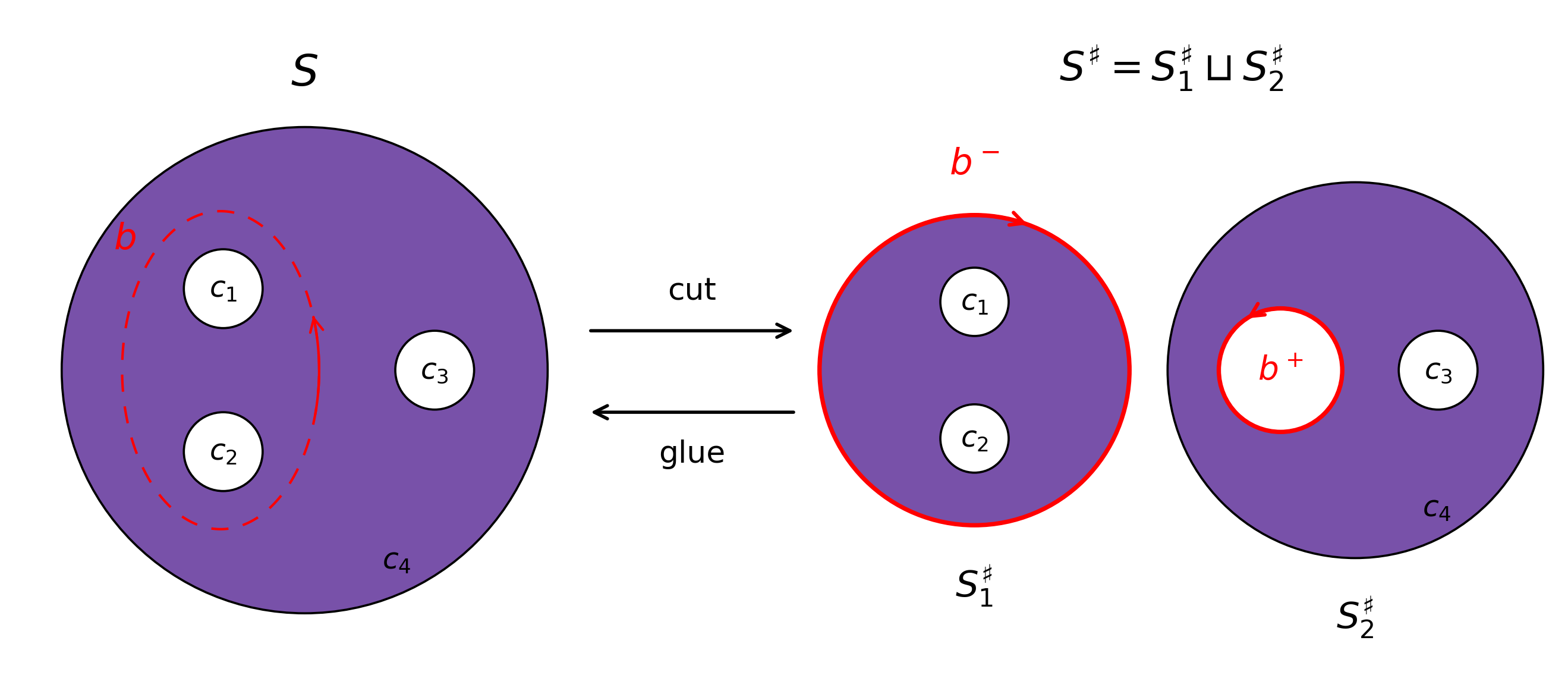}
\caption{Cutting along $b$ and gluing the two new boundary components.}
\label{fig:new-decomposition}
\end{figure}

Given a partition
\[S = \bigsqcup_{1\leq i \leq I} S_i \sqcup \bigsqcup_{1\leq l \leq L} b_l\]
we denote by $S_i^\sharp$ the completion of $S_i$ for the restriction of some Riemannian metric on $S$, and write 
\[S^\sharp = \bigsqcup_{1\leq i \leq I} S_i^\sharp~.\]
It is a compact oriented surface with boundary (typically disconnected) with a natural map $\iota: S^\sharp \to S$ that sends the interior of $S^\sharp$ homeomorphically to
\[S \setminus \left(\bigsqcup_{1\leq l \leq L} b_l \sqcup \bigsqcup_{1\leq k \leq K} c_k\right)~.\]
The boundary components of $S^\sharp$ are of two kinds: 
\begin{enumerate}
    \item for each boundary component $c_k$ of $S$, there is a unique boundary component of $S^\sharp$ which is mapped homeomorphically to  $c_k$ by $\iota$, which we still denote by $c_k$,
    \item for each curve  $b_l$ of the partition, there are exactly two boundary components  $b_l^+$ and $b_l^-$ which are mapped homeomorphically to $b_l$ by $\iota$, and such that  $\iota : b_l^+ \to b_l$ (resp.  $\iota: b_l^- \to b_l$) preserves (resp. reverses) the orientation.
\end{enumerate}

The surface $S$ is thus obtained from $S^\sharp$ by gluing $b_l^+$ to $b_l^-$ by an orientation-reversing homeomorphism. See Figure~\ref{fig:new-decomposition}.

\subsection{Mapping class group}
 
We denote by $\MCG(S,*)$ the group of orientation-preserving diffeomorphisms of $S$ equipped with a lift to $\tilde S$, modulo isotopy. The group $\MCG(S,*)$ acts on $\pi_1(S)$ and the morphism $\MCG(S,*) \to \Aut(\pi_1(S))$ is injective as soon as $S$ has negative Euler characteristic. From now on, we identify $\MCG(S,*)$ with its image in $\Aut(\pi_1(S))$.

The group $\MCG(S,*)$ contains as a normal subgroup the group $\pi_1(S)$, namely the lifts of the identity. The quotient $\MCG(S,*)/\pi_1(S)$ is the usual mapping class group $\MCG(S)$ of $S$, i.e. the group of isotopy classes of (unlifted) orientation-preserving diffeomorphisms of $S$.

For any group $G$, the group $\MCG(S,*)$ acts on $\Hom(\pi_1(S), G)$ by precomposition, and this action descends to an action of $\MCG(S)$ on
the quotient $\Rep(S,G)$ of $\Hom(\pi_1(S),G)$ by the conjugation action of $G$.

\subsection{Volume of representations}

Let $S$ be a surface of finite type and $\rho: \pi_1(S) \to \Isom_+(\H)$ a representation. A notion of \emph{volume} of the representation $\rho$ was developed independently by Koziarz--Maubon \cite{KM} and Burger--Iozzi--Wienhard \cite{BIW}.\footnote{They actually define the  \emph{Toledo invariant} $\Tol(\rho)$ of any representation of a surface group into a semisimple Lie group of Hermitian type. When the target group is $\Isom_+(\H)$, the volume introduced here equals $-2\pi\,\Tol(\rho)$.} Here we denote this volume by $\Vol(\rho)$ and normalize it so that, when the surface is closed, we have
\[\Vol(\rho) = \int_S f^*\vol_\H = -2\pi\, \euler(\rho)~,\]
where $f: \tilde S \to \H$ is any smooth $\rho$-equivariant map and $\euler(\rho)$ is the Euler class of the circle bundle
\[\tilde S \times \partial_\infty \H/(x,p)\sim (\gamma\cdot x, \rho(\gamma)p)~.\]
Goldman \cite{Goldman88} and Hitchin \cite{Hitchin87} independently proved that the Euler class classifies connected components of $\Hom(\pi_1(S),\Isom_+(\H))$.

When $S$ is not closed, the definition is a bit more involved and one must either impose a ``normalized behaviour'' in the ends of $S$ or add a correction term for each boundary component. The volume is not an integral multiple of $2\pi$ anymore. Rather, we have the following:

\begin{prop}[{\cite[Proposition 3.10, Lemma 8.2]{BIW}}] \label{prop: Volume continuous}
    Let $S$ be a compact connected oriented surface with boundary components $c_k, 1\leq k\leq K$. Then 
    \[\Vol: \Hom(\pi_1(S), \Isom_+(\H)) \to \R\]
    is continuous and, for any $\rho \in \Hom(\pi_1(S), \Isom_+(\H))$,   
    \[\Vol(\rho) \equiv -\sum_{k=1}^K r(\rho(c_k)) \pmod{2\pi}~,\]
    where $r: \Isom_+(\H) \to \R/2\pi \Z$ denotes the rotation number.
\end{prop}

In \cite[Proposition-definition 4.1]{KM}, the volume of a representation $\rho$ is defined as the integral of the pull-back of the volume form of $\H$ by any smooth $\rho$-equivariant map $f: \tilde S \to \H$, corrected by some boundary terms. One can verify that these terms are preserved by postcomposing $f$ by an isometry and precomposing it by a diffeomorphism of $S$. It follows that the volume is invariant under conjugation and under the action of $\MCG(S,*)$:
\begin{prop} \label{prop: volume MCG invariant}
    Let $\rho:\pi_1(S)\to \Isom_+(\H)$ be a representation. Then, for every $\tau \in \MCG(S,*)$ and every $g\in \Isom_+(\H)$, we have
    \[\Vol( g(\rho\circ \tau) g^{-1}) = \Vol(\rho)~.\]
\end{prop}

Another crucial property of the volume is its additivity:
\begin{prop}[{\cite[Proposition 3.2]{BIW}}]
\label{prop: Volume Additive}
    Given a decomposition 
    \[S = \bigsqcup_{1\leq i \leq I} S_i \sqcup \bigsqcup_{1\leq l \leq L} b_l~,\]
    for any representation $\rho: \pi_1(S) \to \Isom_+(\H)$, we have
    \[\Vol(\rho) = \sum_{1\leq i \leq I} \Vol(\rho_{\vert S_i})~.\]
\end{prop}

Finally, the volume satisfies the Milnor--Wood inequality:
\begin{prop}[{\cite[Corollary 3.4]{BIW}}]
    For every representation $\rho:\pi_1(S) \to \Isom_+(\H)$, we have
    \[\vert \Vol(\rho)\vert  \leq 2\pi \vert \chi(S)\vert~. \]
\end{prop}

\section{Singular hyperbolic surfaces} \label{s: hyperbolic surfaces}

\subsection{The hyperbolic plane and its isometries}

Recall that ``the'' hyperbolic plane is the (unique up to isometry) complete simply connected Riemannian surface of constant Gauss curvature $-1$. It admits several models, among which:
\begin{itemize}
    \item the upper half-plane model \[\H = \{z\in \C \mid \im(z)>0\}\]
    equipped with the conformal metric  \[\sigma_\H = \frac{\vert \d z\vert^2}{\im(z)^2}~,\]
    \item the Poincar\'e disc
    \[\D= \{z\in \C \mid \vert z \vert <1 \} \]
    equipped with the conformal metric \[\sigma_\D = \frac{4 \vert \d z \vert^2}{(1-\vert z \vert^2)^2}~.\]
\end{itemize}

Their groups of orientation-preserving isometries are respectively
\[\Isom_+(\H) = \PSL(2,\R)\]
and
\[\Isom_+(\D) = \PU(1,1)~,\]
where $\PU(1,1)$ and $\PSL(2,\R)\subset \PSL(2,\C)$ act on $\C\cup \{\infty\}$ by homographies.\\

Recall also that conjugacy classes of isometries of $\H$ fall into two one-parameter families and 3 additional classes:
\begin{itemize}
    \item An \emph{elliptic isometry} fixes a unique point in $\H$ and rotates about it with an angle $\alpha \in \R/2\pi\Z$. We simply call it a \emph{rotation of angle $\alpha$}. It is conjugate to
    \[z\mapsto e^{i\alpha} z\]
    acting on $\D$.

    \item A \emph{hyperbolic isometry} preserves a geodesic called its \emph{axis} and acts along it as a translation of length $\ell$. It is conjugate to
    $z\mapsto e^\ell z$
    acting on $\H$.

    \item A \emph{parabolic isometry} fixes a unique point in $\partial_\infty \H$ and translates along the horocycles centered at this point. Parabolic isometries come in two conjugacy classes: \emph{positive parabolics}, which are conjugate to
    \[z\mapsto z+1\]
    acting on $\H$, and \emph{negative parabolics}, which are conjugate to
    \[z\mapsto z-1~.\]

    \item The identity is its own conjugacy class.
\end{itemize}

\begin{rmk}
    Though positive and negative parabolics are not conjugate in $\Isom_+(\H)$, they are conjugate by an \emph{orientation-reversing} isometry of $\H$. The same holds for rotations of angles $\alpha$ and $2\pi-\alpha$.
\end{rmk}

The \emph{translation length} of an isometry $g\in \Isom_+(\H)$ is the number
\[\ell(g) = \inf_{x\in \H} d(x,g\cdot x)~.\]
It is $0$ if $g$ is elliptic, parabolic or the identity. When $g$ is hyperbolic, it is a translation of length $\ell(g)$ along its axis. 

Given an element $g\in \Isom_+(\H) \simeq \PSL(2,\R) \simeq \PU(1,1)$ and a lift $\hat g$ of $g$ to $\SL(2,\C)$, the number $\vert \Tr(g) \vert \coloneqq \vert \Tr(\hat{g})\vert$ does not depend on the choice of the lift $\hat{g}$. It is a conjugacy invariant of $g$ and is given in the following way:
\begin{itemize}
    \item If $g$ is elliptic with rotation angle $\alpha$, then
    \[\vert\Tr(g)\vert = \vert 2 \cos(\frac{\alpha}{2}) \vert < 2~.\]
    \item If $g$ is hyperbolic with translation length $\ell$, then
    \[\vert\Tr(g)\vert = 2\cosh(\ell/2)>2~.\]
    \item If $g$ is parabolic or the identity, then
    \[\vert\Tr(g)\vert = 2~.\]
\end{itemize}

An element $g\in \Isom_+(\H)$ acts on the circle $\partial_\infty \H$ with a \emph{rotation number} $r(g) \in \R/2\pi \mathbb Z$. The rotation number is a continuous conjugacy invariant. It takes the value $\alpha$ if $g$ is an elliptic element of rotation angle $\alpha$ and $0$ otherwise.

Given a representation \(\rho: G\rightarrow \Isom_+(\H)\) of a finitely generated group \(G\) equipped with a finite generating set \(F\), one defines, for every \(x\in \mathbb H\), the \emph{displacement}
\[ \delta_{\rho,F} (x) \coloneqq \max_{f\in F} d(x, \rho(f)x)~, \]
and the \emph{minimal displacement} of the pair \( (\rho, F) \) as the number 
\[ \delta_{\rho,F}= \inf_{x\in \mathbb H} \delta_{\rho,F}(x). \]

It vanishes if and only if \(\rho \) is conjugate to a  representation with values either in the group \(\text{PSO}(2)\) or in the group of unipotent upper triangular matrices. A fundamental fact is that the minimal displacement of a representation depends continuously on the representation\footnote{While upper semicontinuity is immediate, lower semicontinuity is a non-trivial fact.} (see \cite[Chapters II.6 to II.8]{BridsonHaefliger13}).

\subsection{Hyperbolic cones and cusps}
Via the exponential map at a point, the hyperbolic plane identifies with $\R^2$ equipped with the metric
\[\d r^2 + \sinh^2(r) \d \theta^2\]
in polar coordinates $(r,\theta)$. This expression can be generalized to give a model for a hyperbolic metric with cone singularities.

\begin{defi}
For $\alpha >0$, we define the \emph{standard hyperbolic cone} $C_\alpha$ of angle $\alpha$ as the plane $\R^2$ equipped with the singular Riemannian metric given in polar coordinates $(r,\theta)$, with $\theta$ of period $2\pi$, by
\[\d r^2 + \left(\frac{\alpha}{2\pi}\right)^2 \sinh^2(r) \d \theta^2 ~.\]
\end{defi}

The cone $C_{2\pi}$ is isometric to the hyperbolic plane, and the metric on $C_{\alpha}$ is (locally) the pull-back of the metric of $C_{2\pi}$ by the map $r e^{i\theta} \mapsto  r e^{i\alpha\theta/(2\pi)}$. In particular, it is smooth and hyperbolic away from $0$.

\begin{figure}[htbp]
\centering
\includegraphics[width=0.98\textwidth]{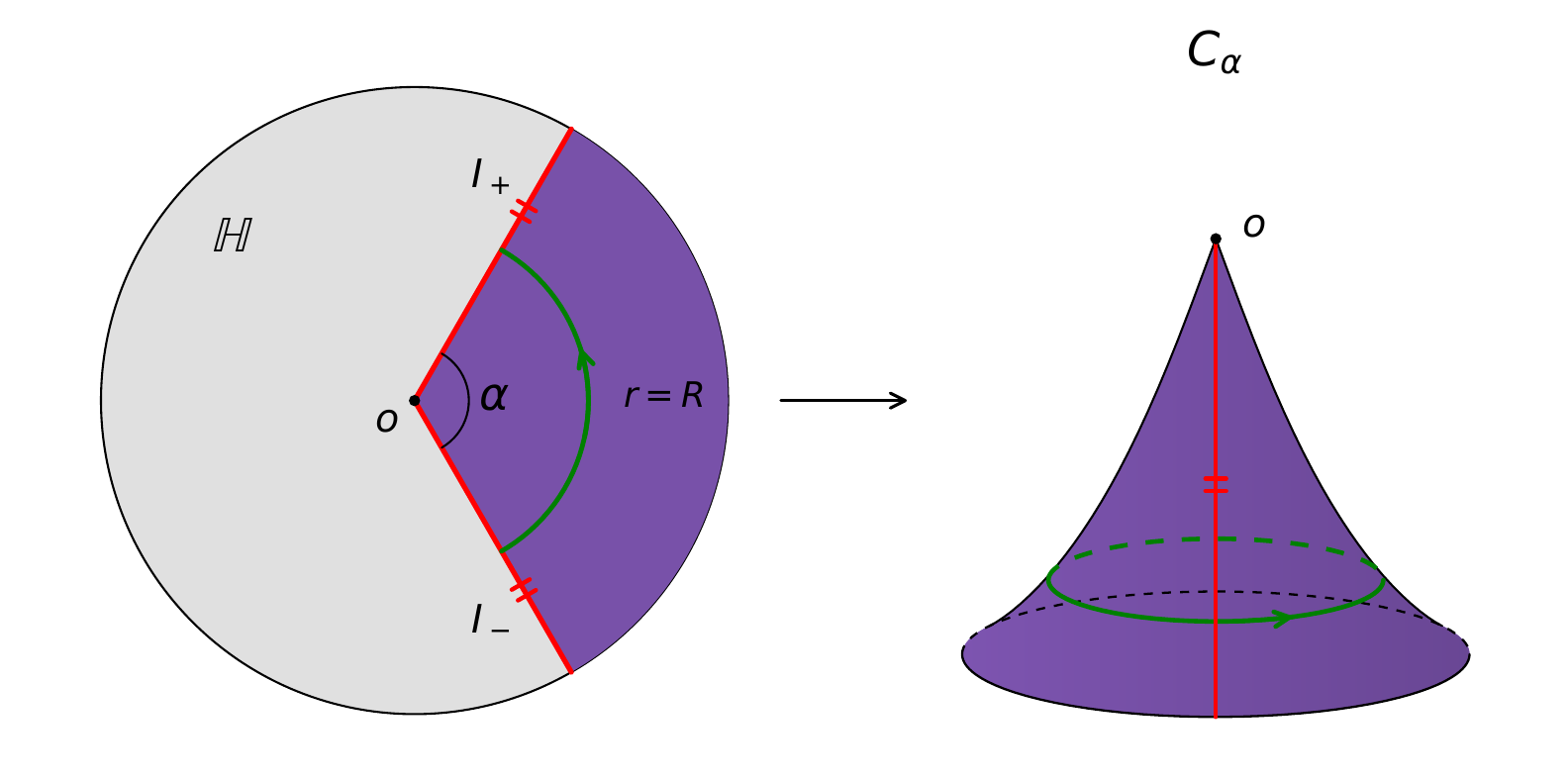}
\caption{A hyperbolic sector and a local surface-of-revolution model of the cone obtained by identifying its red sides. The green curves correspond to the same metric circle. The three-dimensional view shows a neighbourhood of the vertex, for $0<\alpha<2\pi$.}
\label{fig:new-cone}
\end{figure}

The cone $C_\alpha$ also admits a conformal model:
\begin{prop}
The cone $C_\alpha$ is isometric to the unit complex disc $\D$ equipped with the conformal metric
\[\frac{(\alpha/\pi)^2\vert z\vert^{\alpha/\pi-2}}{(1-\vert z\vert^{\alpha/\pi})^2}\,\vert\d z\vert^2~.\]
\end{prop}

When $\alpha$ goes to $0$, the cone $C_\alpha$ ``converges'' to a hyperbolic cusp:

\begin{defi}
We define the \emph{standard hyperbolic cusp} $C_0$ as the plane $\R^2$ equipped, on the complement of the origin, with the metric given in polar coordinates $(r,\theta)$ by
\[
\frac{\d r^2}{r^2}+r^2\d\theta^2.
\]
\end{defi}

\begin{rmk}
    While the origin is at infinite distance, we include it as a point in $C_0$, which will prove useful when discussing singular metrics on a closed surface.
\end{rmk}

The hyperbolic cusp also admits a conformal model:
\begin{prop}
The cusp $C_0$ is isometric to the complex disc $\D$ equipped with the conformal metric on the complement of the origin given by
\[\frac{\vert d z\vert^2}{\vert z \vert^2 \log^2  \vert z \vert }~.\]
\end{prop}

\begin{figure}[htbp]
\centering
\includegraphics[width=0.98\textwidth]{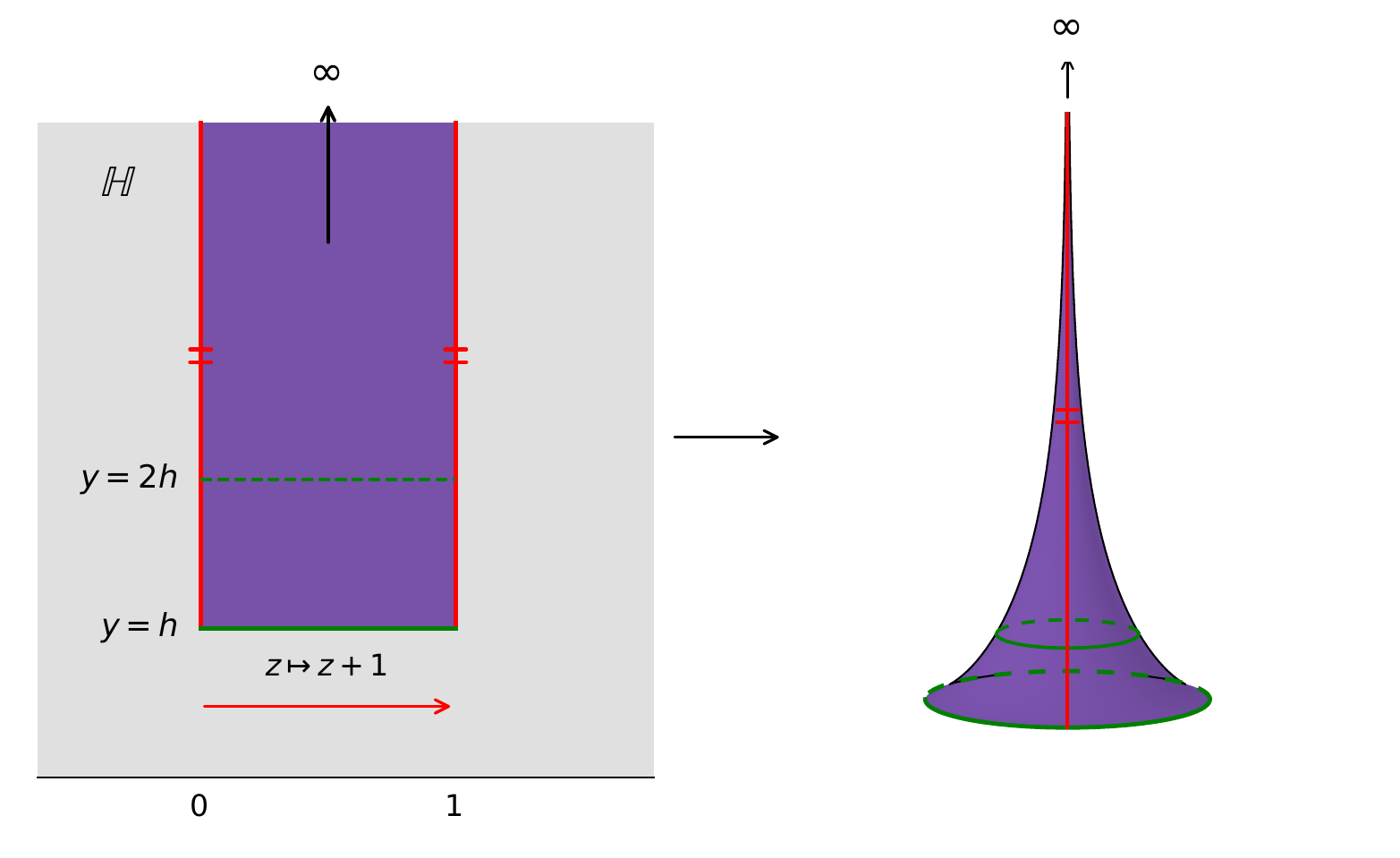}
\caption{A cusp end: the two red vertical sides are identified by $z\mapsto z+1$. On the right, a surface-of-revolution model of the end; its tip is at infinite distance. The two green horocycles have lengths in ratio $2:1$.}
\label{fig:new-cusp}
\end{figure}

\subsection{Singular hyperbolic structures on surfaces with boundary}

Recall that an oriented surface $S$ equipped with a Riemannian metric of constant curvature $-1$ is locally isometric to the hyperbolic plane, via an orientation-preserving map. These local isometries define an atlas on $S$ with values in $\H$ and transition maps in $\Isom_+(\H)$; that is,  a \emph{$(\Isom_+(\H),\H)$-structure} in  Thurston's terminology. Here we introduce a variation of this notion for surfaces locally modeled on hyperbolic cones.

Let $S$ be a compact oriented topological surface with boundary.

\begin{defi}
A singular hyperbolic structure on $S$ is given by a maximal atlas $(U_i,\alpha_i, \phi_i)$, where $(U_i)_{i\in I}$ is an open covering of $S$ and $\phi_i:  U_i \to C_{\alpha_i}$ is an orientation-preserving homeomorphism onto its image such that
\begin{itemize}
\item $\phi_i(U_i)$ is locally closed with piecewise $\mathcal C^2$ boundary,
\item $\phi_j\circ \phi_i^{-1}: \phi_i(U_i\cap U_j) \to \phi_j(U_i\cap U_j)$ is a local isometry.
\end{itemize}
\end{defi}

If $x$ is an interior point of $S$, there exists a unique $\alpha(x)$ for which there is an open neighbourhood $U_x$ of $x$ and a chart $\phi_x : U_x \to C_{\alpha(x)}$ mapping $x$ to $0$. We call $\alpha(x)$ the \textit{angle} of the hyperbolic structure at $x$. The point $x$ is called \emph{regular} if $\alpha(x) = 2\pi$, \emph{branched} if $\alpha(x) = 2k\pi$, $k\geq 2$ and \emph{singular} otherwise. A singular point $x$ is a \emph{cone singularity of angle $\alpha(x)$} if $\alpha(x)> 0$ and a \emph{cusp} if $\alpha(x)=0$. Note that branched and singular points are isolated.

Now let $x$ be a boundary point of $S$. Let $(U_x, \beta, \phi_x)$ be a local chart in a neighbourhood of $x$ such that $\phi_x(x) =0$. By definition, the image of $U_x$ (intersected with a sufficiently small ball centered at $0$) is a domain bounded by two $\mathcal C^1$ curves starting at $0$. 

Let $\gamma$ be the angle (for the standard Euclidean metric) between the tangents to these two curves at $0$. Then the quantity $\alpha(x) = \frac{\beta \gamma}{2\pi}$ does not depend on the choice of the local chart. We call it the \textit{angle} of the hyperbolic structure at $x$. We call $x$ \emph{regular} if $\alpha(x)= \pi$ and \emph{a corner} otherwise.

A singular hyperbolic structure defines a path metric on the complement
of the cusps. At an interior point $x$ it is locally isometric to
$C_{\alpha(x)}$ near its vertex; at a boundary point it is locally
isometric to a domain in a cone as in the preceding definition.
Conversely, these local models, together with the specified cusp ends,
define a singular hyperbolic structure. We call the resulting path metric
a \emph{singular hyperbolic metric}. In the interior, away from singular
and branch points, it is a smooth Riemannian metric of curvature $-1$.

Assume now that $S$ is connected, and let $P\subset S$ denote the set of singular points of $S$. Then there exists a representation $\rho: \pi_1(S\setminus P) \to \Isom_+(\H)$ called a \emph{holonomy representation} of the hyperbolic structure and a $\rho$-equivariant, orientation-preserving, locally isometric map\footnote{Here, \emph{locally isometric} means that it preserves the length of a path.} $\dev: \widetilde{S\backslash P} \to \H$ called a \emph{developing map}. We call the pair $(\dev,\rho)$ a \emph{developing pair} for the hyperbolic structure. If $(\dev', \rho')$ is another developing pair, then there exists a unique $g\in \Isom_+(\H)$ such that $\dev' = g\circ \dev$ and $\rho' = g\rho g^{-1}$. More informally, the developing map and holonomy representation are respectively unique up to postcomposition and conjugation by an isometry. When this causes no confusion, we will simply refer to them as \emph{the} developing map and holonomy representation. Finally, for a surface with multiple connected components, a \emph{developing map} and \emph{holonomy representation} are understood as a collection of developing maps and holonomy representations for each component.

\begin{prop}
    Let $S$ be a singular hyperbolic surface, $p\in P$ a singular point of $S$ and $p^\circ$ a small simple loop encircling the point $p$.
    \begin{itemize}
        \item if $p$ is a cusp, then $\rho(p^\circ)$ is a positive parabolic,
        \item if $p$ is a cone singularity of angle $\alpha$, then $\rho(p^\circ)$ is a rotation of angle $\alpha \mod 2\pi$.
    \end{itemize}
\end{prop}

\begin{proof}
Since this is a local property around the point $p$, it is enough to compute the holonomy of the hyperbolic cone $C_\alpha$. Let $\widetilde{C_\alpha\setminus \{0\}}$ denote $\R_{>0} \times \R$ equipped with the hyperbolic metric pulled back from $C_\alpha \setminus \{0\}$ by the universal covering map \[(r,\theta) \mapsto (r\cos(\theta),r\sin(\theta))~.\]

If $\alpha >0$, one can verify that the map
\[\dev: \tilde{C_\alpha \setminus \{0\}} \to C_{2\pi} \simeq \H\]
given by 
\[\dev(r,\theta) = \left(r \cos \left(\frac{\alpha\theta}{2\pi}\right), r\sin\left(\frac{\alpha\theta}{2\pi}\right)\right)\]
is a local isometry and satisfies
\[\dev(r,\theta+2\pi) = r_\alpha \dev(r,\theta)~,\]
where $r_\alpha$ is the rotation of angle $\alpha$ about $0$. This gives the holonomy of $C_\alpha$ at $0$.

If $\alpha = 0$, we similarly verify that the map
\[\dev: \widetilde{C_0\setminus \{0\}} \to \H\]
given by
\[(r,\theta) \mapsto (\theta+ \frac{i}{r})\]
is a local isometry and satisfies
\[\dev(r,\theta+2\pi) = \dev(r,\theta) +2\pi~,\]
showing that the holonomy of $C_0$ at $0$ is a positive parabolic. 
\end{proof}

The image of the boundary of $S$ by $\dev$ is a piecewise $\mathcal C^2$ curve. Except at finitely many points, it therefore has a well-defined \emph{curvature}. To be more precise, let $\sigma$ be a singular hyperbolic metric on $S$ and $\dev$ a developing map for $\sigma$. For every $x \in \partial S$, let $c: (-\epsilon, \epsilon) \to \partial S$ be an oriented isometric local parametrization of $\partial S$ such that $c(0) = x$, and $\tilde c$ a lift of $c$ to $\tilde S$. Then, for all but finitely many points $x$, the parametrized curve $\gamma \coloneqq \dev \circ \tilde c$ is of class $\mathcal C^2$ near $0$, and we can define
\[
\kappa(x)\coloneqq
-\vol_{\H}\bigl(\gamma'(0),
(\nabla_{\gamma'}\gamma')(0)\bigr),
\]
where $\nabla$ is the Levi-Civita connection of $\H$.
The minus sign reflects our convention that the surface lies
to the right of its oriented boundary. Thus the boundary of
a hyperbolic disc of radius $r$, oriented clockwise, has
curvature $\coth r>0$.

We call the boundary \emph{locally convex} (resp. \emph{locally concave}) at $x$ if $\kappa(x) >0$ (resp. $\kappa(x)<0$).

\subsection{Singular hyperbolic structures and equivariant holomorphic maps} \label{ss:hyperbolic structures and holomorphic maps}

Given that the standard cones admit conformal models, we see that a singular hyperbolic metric on $S$ defines a complex structure on $\INT(S)$ (including the singular points), for which the local charts (with values in the conformal models of the cones) are holomorphic. In particular, the developing map is holomorphic. We call this complex structure the \emph{conformal class} of the singular hyperbolic metric.

Conversely, let $(S,J)$ be a closed Riemann surface, $P\subset S$ a finite set of points, $\rho$ a representation of $\pi_1(S\backslash P)$ into $\Isom_+(\H)$ and $\dev: \tilde{S\backslash P} \to \H$ a \emph{non-constant} $\rho$-equivariant holomorphic map.

\begin{prop}\label{prop: holomorphic map -> singular metric}
The pull-back metric $\dev^*\sigma_\H$ descends to a singular hyperbolic metric on $S$ such that every $x\in S\setminus P$ is either regular (if any lift is a regular point of $\dev$) or a branch point (if any lift is a critical point of $\dev$).
\end{prop}

\begin{proof}
For any $x\in S\setminus P$, there is a local coordinate $z$ centered at $x$ and an isometry $g\in \Isom_+(\H)$ such that $g\circ \dev$ has the form $z\mapsto z^k$ for some positive integer $k$ near some lift $\tilde{x}$ of $x$, with $k=1$ if and only if $\dev'(\tilde{x}) \neq 0$. Hence the pulled-back metric $\dev^*\sigma_\H$ has a cone singularity of angle $2k\pi$ at $x$, which is thus regular if $k=1$ and a branch point otherwise.

Fix $p\in P$ and a coordinate $z$ centered at $p$. Write $F$ for the local developing map in these coordinates.
Let $D$ be a small disc centred at $p$, and let $g$ denote the holonomy around a positively oriented loop about $p$.
By the Schwarz--Pick lemma $\dev: \widetilde{D\setminus \{p\}} \to \H$ is contracting for the uniformizing metric of $D\setminus \{p\}$. Since this uniformizing metric has a cusp at $p$, we get that $d_\H(g\,\dev(x), \dev(x))\underset{x\to p}{\longrightarrow} 0$. Hence $g$ has translation length $0$ and $g$ is elliptic, parabolic, or the identity.

If $g$ is elliptic, use the disc model and conjugate $g$ to
$w\mapsto e^{2\pi ia}w$, with $0<a<1$.
Equivariance implies that
\[
    h(z)=z^{-a}F(z)
\]
is single-valued and holomorphic on the punctured disc.
Since $|h(z)|\leq |z|^{-a}$ and $a<1$, the function $zh(z)$ extends
holomorphically with value $0$ at $0$, and hence $h$ extends as well.
Write $h(z)=z^m u(z)$, where $m\geq0$ and $u(0)\ne0$.
Choosing a holomorphic determination of $u^{1/(a+m)}$ and setting
$\zeta=z\,u(z)^{1/(a+m)}$, we obtain $F=\zeta^{a+m}$.
The metric therefore has cone angle $2\pi(a+m)$.

If $g$ is the identity, $F$ is single-valued and bounded in the disc
model, hence extends holomorphically across $0$ with value in the disc.
After a target isometry and a change of coordinate, $F=\zeta^k$
for some $k\geq1$, giving angle $2k\pi$.

Finally, if $g$ is parabolic, conjugate it in the upper half-plane
model to $w\mapsto w+\varepsilon$, where $\varepsilon\in\{1,-1\}$.
The function $Q(z)=\exp(2\pi iF(z))$ is single-valued, nonvanishing
and bounded by $1$ on the punctured disc, and extends across $0$.
Its winding number about $0$ along a positive small circle is
$\varepsilon$. This equals its order of vanishing at $0$, so
$\varepsilon=1$ and $Q$ has a simple zero.
Thus $Q$ is a local coordinate, in which
\[
     F^*\sigma_\H = \frac{|dQ|^2}{|Q|^2\log^2|Q|},
\]
the standard cusp metric.
\end{proof}

For holonomies of singular hyperbolic metrics the volume has the following natural interpretation:

\begin{prop} \label{prop: volume singular metric = volume holonomy}
    Let $S$ be a closed surface, $\sigma$ be a singular hyperbolic metric on $S$, with singular set $P$, and $\rho: \pi_1(S\setminus P) \to \Isom_+(\H)$ its holonomy representation. Then
    \[\Vol(\rho) = \int_{S\setminus P} \vol_\sigma~.\]
\end{prop}

\begin{proof}
Let $j: \pi_1(S\setminus P) \to \Isom_+(\H)$ be the representation uniformizing $(S\setminus P, J)$, i.e. such that there exists a $j$-equivariant biholomorphism $u: \tilde{S\setminus P} \to \H$. Then the map $\dev \circ u^{-1}: \H \to \H$ is $(j,\rho)$-equivariant and holomorphic. In particular, it is Lipschitz, and we thus have
\[\Vol(\rho)=\int_{j(\pi_1(S\setminus P))\backslash \H} (\dev\circ u^{-1})^*\vol_\H~,\]
see for instance \cite[Proposition 3.7]{TholozanVolAdS}.
Finally, by construction, we have $(\dev\circ u^{-1})^*\vol_\H = u_*\dev^*\vol_\H = u_* \vol_\sigma$, and the conclusion follows since $u$ factors to a global diffeomorphism $S\setminus P \overset{\sim}{\to} j(\pi_1(S\setminus P))\backslash \H$.
\end{proof}

\subsection{Cutting and gluing} \label{ss: cutting and pasting}

A way to build (possibly singular) hyperbolic surfaces is to glue together pieces of hyperbolic surfaces along their boundaries. While this kind of operation is geometrically intuitive, it is not easy to formalize. Since we will make extensive use of gluing constructions, we sketch in this section the basics of a sufficiently general theory of gluing.\\

Let $S'$ be a compact oriented surface with boundary. Let $I'$ be either a closed interval or a circle contained in $\partial S'$ and let $f: I' \to \partial S'$ be a continuous injective orientation-reversing map such that $f(\INT{I'}) \cap I' = \emptyset$. Then the quotient space
\[S = S' /\{ f(x) \sim x, x\in I'\}\]
is a topological surface with boundary.

\begin{prop} \label{p: paste}
Assume that $S'$ carries a singular hyperbolic structure, and that $f$ is an isometry that satisfies $$\kappa(f(x)) = -\kappa(x)\text{ for all }x\in I'$$ where the curvatures are defined. Then there exists a unique singular hyperbolic structure on $S$ such that the projection $p: S' \to S$ is a local isometry. Moreover, for every $y\in S$, we have
\begin{equation}\label{eq: angles after gluing} \alpha(y) = \sum_{x\in p^{-1}(y)} \alpha(x)~.\end{equation}
\end{prop}

In order to prove this proposition, we will use the following fact:

\begin{fact} \label{fact: curve characterized by curvature}
    Fix some $\alpha>0$. Let $u, v: [0,t_0) \to C_\alpha$ be two curves of class $\mathcal C^2$ parametrized at unit speed such that:
    \begin{itemize}
        \item $u(0) = v(0)$ and $u(t), v(t)\neq 0$ for $t>0$,
        \item $u'(0)= v'(0)$,
        \item $\kappa(u(t)) = \kappa(v(t))$ for all $t\in [0,t_0)$.
    \end{itemize}
    Then $u(t) = v(t)$ for all $t\in [0,t_0)$.
\end{fact}

\begin{proof}[Proof of Proposition \ref{p: paste}]
Clearly, the singular hyperbolic structure on $S'\setminus (I'\cup f(I'))$ induces a singular hyperbolic structure on $S\setminus p(I')$. The goal is to extend this hyperbolic structure to the gluing points. 

Given a piecewise $\mathcal C^2$ curve $c$ from $(-\epsilon,\epsilon)$ to a manifold $M$, let us denote by $d_+(c)$ and $d_-(c)$ the derivatives at $0$ in positive time of $t\mapsto c(t)$ and $t\mapsto c(-t)$.

Let $x$ be a point in $\INT(I')$ with angle $\alpha$ and let $\beta$ be the angle of $f(x)$. Let $u: (-\epsilon, \epsilon) \to I'$ be a local parametrization of $I'$ with $u(0)= x$. By definition of our singular hyperbolic structures, there exist open neighbourhoods $U\ni x$, $V\ni f(x)$ and local coordinates $\dev_U: U \to C_{\alpha+\beta}$, $\dev_V: V\to C_{\alpha+\beta}$ such that:
\begin{itemize}
    \item $\dev_U(x) = \dev_V(f(x)) = 0~,$
    \item $\dev_U\circ u$ is piecewise $\mathcal C^2$, $d_-(\dev_U\circ u)$ and $d_+(\dev_U\circ u)$ form an angle $\frac{2\pi \alpha}{\alpha+\beta}$ and $\dev_U(U)$ lies to the right of $\dev_U\circ u(-\epsilon,\epsilon)$,
    \item $\dev_V\circ f\circ u$ is piecewise $\mathcal C^2$, $d_-(\dev_V\circ f\circ u)$ and $d_+(\dev_V \circ f\circ u)$ form an angle $\frac{2\pi \alpha}{\alpha+\beta}$ and $\dev_V(V)$ lies to the left of $\dev_V\circ f\circ u(-\epsilon,\epsilon)$.
\end{itemize}
Up to composing $\dev_V$ by a rotation, we can further assume that $d_-(\dev_U\circ u) = d_-(\dev_V\circ f\circ u)$, and the equality of angles implies that we also have $d_+(\dev_U\circ u) = d_+(\dev_V\circ f\circ u)$. By hypothesis on $f$, we have $\kappa(\dev_V \circ f\circ u(t)) = \kappa (\dev_U \circ u(t))$ for all $t$, and Fact \ref{fact: curve characterized by curvature} gives $\dev_V\circ f \circ u = \dev_U \circ u$. We conclude that $\dev_U$ and $\dev_V$ patch together into a local homeomorphism from $S$ to $C_{\alpha+\beta}$ near $p(x)$. Hence the quotient metric has a cone singularity of angle $\alpha+\beta$ at $p(x)$.

Similarly, if $I'$ is an interval and $x$ is one of its endpoints, one can show that we are in one of two cases:
\begin{itemize}
    \item if $f(x)\neq x$, then $p(x)$ is a boundary point of $S$ of angle $\alpha+\beta$,
    \item if $f(x) =x$, then $f$ is gluing the two sides of $x \in \partial S$ and $p(x)$ becomes a singularity of angle $\alpha$.
\end{itemize}
\end{proof}

The reverse operation can be described as follows. Let $S$ be an oriented topological surface with boundary. Let $I$ be either an embedded interval or  circle in $S$ with interior contained in $\INT(S)$. Then there exists a surface with boundary $S'$ equipped with an interval or circle $I'$ in its boundary, a continuous injective map $f:I' \to \partial S'$ such that $f(\INT(I')) \cap I' = \emptyset$ and a projection map $ p : S' \rightarrow S $ such that $p$ and $p\circ f$ are respectively orientation-preserving and reversing homeomorphisms from $I'$ to $I$, and the restriction $p: S'\backslash (I'\cup f(I')) \to S\backslash I$ is a homeomorphism.

\begin{prop}\label{p: cut}
Assume that $S$ is endowed with a singular hyperbolic structure, and that~$I$ can be arclength parametrized by a piecewise $\mathcal C^2$ curve in $S$. Then there exists a unique singular hyperbolic structure on $S'$ such that the projection $p$ is a local isometry. Moreover, for every $y\in S$, equation~\eqref{eq: angles after gluing} holds.
\end{prop}

\begin{proof}
    We equip $S'$ with the pull-back under $p$ of the singular hyperbolic metric on $S$. Again, near the preimage of a point in $S\setminus I$, $p$ is a local homeomorphism and the pulled-back metric is thus a singular hyperbolic metric.

    Let now $x$ be a point in $\INT(I)$. Let $I_+, I_-$ be the two preimages of $I$ in $S'$, and $x_+\in I_+$ and $x_-\in  I_-$ the two preimages of $x$. Let $U$ be a small neighbourhood of $x\in S$, $u:(-\epsilon, \epsilon) \to U\cap I$ a parametrization of $I$ at unit speed such that $u(0) = x$, and $\dev: U \to C_\alpha$ a local chart for the singular hyperbolic structure. Then $\dev\circ p$ is a local chart from a neighbourhood of $x_-$ (resp. $x_+$) to the left (resp. right) of $\dev\circ u(-\epsilon,\epsilon)$ in $\dev(U)$. Since $\dev\circ u$ is piecewise $\mathcal C^2$, we get that the pulled-back metric is a singular metric near $x_-$ and $x_+$, with $x_-$ and $x_+$ boundary points of respective angles $\alpha_-$ and $\alpha_+$ such that $\alpha_-+\alpha_+ = \alpha$.\\

    Similarly, if $I$ is an interval and $x$ is an endpoint of $I$, then we are in one of two cases:
    \begin{itemize}
        \item  either $x\in \partial S$, in which case $x$ has two  preimages $x_-$ and $x_+$ in $S'$. Again, precomposing a local chart with $p$, we get that $x_-$ and $x_+$ are corners for the pulled-back metric, of respective angles $\alpha_-$ and $\alpha_+$ such that $\alpha_-+\alpha_+$ is the angle of $x$.
        \item or $x\in \INT(S)$. Then the two preimages $I_-$ and $I_+$ have as a common  endpoint the unique  preimage $p^{-1}(x)$ of $x$. The map $\dev\circ p$ near $p^{-1}(x)$ can be lifted to the double cover $C_{2\alpha}$ of $C_\alpha$, defining a local chart near $p^{-1}(x)$, and we get that $p^{-1}(x)$ is a corner of angle $\alpha$.
    \end{itemize}
\end{proof}

We now discuss the relationship between the holonomies of the singular hyperbolic structures arising in Propositions \ref{p: paste} and \ref{p: cut}.

Let $S$, $S'$, $I$, $I'$ and $f$ be as above.
For the following two holonomy lemmas, assume that $S$ is connected
and that $I$ is either an embedded circle in $\INT(S)$ or a properly
embedded compact interval, meaning that
\[
I\cap\partial S=\partial I.
\]
Assume also that $I$ contains neither a singular point nor a branched
point. Let $P$ be the set of singular points of $S$ and set
$P'=p^{-1}(P)$.
The two copies $I'$ and $f(I')$ of the cutting curve are then
disjoint, including their endpoints when $I$ is an interval.
Fix a lift $\tilde I$ of $I$ to $\widetilde{S\setminus P}$.

We have two separate situations, depending on whether $I$ disconnects $S$ or not.

\vspace{0.2cm}

\subsubsection{{\bf Separating case.}} \label{sss:Van Kampen Separating}
If \(I\) disconnects \(S\), $S'$ has two connected components $S_1'$ and $S_2'$. We may assume that \( I'\) is contained in \(S_1'\). By assumption,  the set of singular points $P'$ of \(S'\) does not intersect $I '$. So it splits as
\[P' = P_1' \sqcup P_2' ~,\]
with $P_i' \subset S_i'$. Since singular points of \(S\) do not lie in \( I\), the projection from \(S'\) to \(S\) induces a bijection \(P' \rightarrow P\). 

Let $\tilde{I}'$ be a lift of $I'$ to $\tilde{S_1'\setminus P_1'}$, and $\tilde f: \tilde{I'}\to \tilde{S_2'\setminus P_2'}$ a lift of $f$. We can lift the inclusions
\[\iota_i: S_i'\setminus P_i' \to S\setminus P\]
to 
\[\tilde{\iota}_i: \tilde{S_i'\setminus P_i'} \to \tilde{S\setminus P}\]
in such a way that $\tilde{\iota}_1(\tilde I') = \tilde{\iota}_2(\tilde f(\tilde I')) = \tilde I$. By Van Kampen's theorem, the homomorphisms ${\iota_i}_*: \pi_1(S_i'\setminus P_i') \to \pi_1(S\setminus P)$ induce an isomorphism:
\[\pi_1(S_1'\setminus P_1') \underset{\pi_1(I)}{*} \pi_1(S_2'\setminus P_2') \overset{\sim}{\to} \pi_1(S\setminus P)~.\]
(When $I$ is a properly embedded interval, this is simply a free product.)

Let $(\dev_1,\rho_1)$ and $(\dev_2, \rho_2)$ denote developing pairs for the hyperbolic structures on $S_1'\backslash P_1'$ and $S_2 '\backslash P_2 '$.

\begin{lem} \label{lem: holonomy gluing separating case}
There exists a unique isometry $g\in\Isom_+(\H)$ such that
\[
\dev_2\circ\tilde f=g\circ\dev_1|_{\tilde I'},
\]
and a unique developing pair $(\dev,\rho)$ for the hyperbolic
structure on $S\setminus P$ such that
\begin{itemize}
\item $\dev\circ\tilde\iota_1=\dev_1$ and
$\dev\circ\tilde\iota_2=g^{-1}\circ\dev_2$;
\item $\rho\circ{\iota_1}_*=\rho_1$ and
$\rho\circ{\iota_2}_*=\Ad_g^{-1}\circ \rho_2$.
\end{itemize}
\end{lem}

\begin{proof}
Choose a developing pair $(\dev_0,\rho_0)$ on $S\setminus P$.
By uniqueness of developing pairs, there are
$g_1,g_2\in\Isom_+(\H)$ such that
\[
\dev_i=g_i\circ\dev_0\circ\tilde\iota_i,
\qquad
\rho_i=\Ad_{g_i}\circ\rho_0\circ(\iota_i)_*.
\]
Set $(\dev,\rho)=(g_1\circ\dev_0,\Ad_{g_1}\circ\rho_0)$
and $g=g_2g_1^{-1}$.
The required restrictions follow immediately, and
$\tilde\iota_2\circ\tilde f=\tilde\iota_1|_{\tilde I'}$ gives
\[
\dev_2\circ\tilde f
=g\circ\dev\circ\tilde\iota_2\circ\tilde f
=g\circ\dev_1|_{\tilde I'}.
\]
The isometry $g$ is unique, since an orientation-preserving isometry
fixing pointwise the non-constant curve $\dev_1(\tilde I')$
is the identity. The normalization on the first component uniquely
determines the developing pair.
\end{proof}

\subsubsection{{\bf Non-separating case.}} \label{sss: Van Kampen non separating}

Assume now that $S'$ is connected.
Choose a lift $\tilde I'$ of $I'$ to $\widetilde{S'\setminus P'}$
and a lift $\tilde f:\tilde I'\to\widetilde{S'\setminus P'}$ of $f$.
Let $\iota:S'\setminus P'\to S\setminus P$ be the restriction of $p$,
and choose a lift
$\tilde\iota:\widetilde{S'\setminus P'}\to\widetilde{S\setminus P}$
such that $\tilde\iota(\tilde I')=\tilde I$.
There is a unique $\gamma\in\pi_1(S\setminus P)$ such that
\[
\tilde\iota\circ\tilde f
=\gamma\circ\tilde\iota|_{\tilde I'}.
\]
By Van Kampen's theorem, we have an isomorphism

\[\pi_1(S'\setminus P') *_f \langle \gamma\rangle \overset{\sim}{\to} \pi_1(S\setminus P)\]
where $\pi_1(S'\setminus P') *_f \langle \gamma\rangle$ denotes the \emph{HNN extension}:
\[\pi_1(S'\setminus P') * \langle \gamma\rangle /\{\gamma \mu \gamma^{-1} = f_*\mu, \mu \in \pi_1(I')\}~.\]
(When $I$ is a properly embedded interval, this is simply a free product.)

Let $(\dev',\rho')$ denote a developing pair for the hyperbolic structure on $S'\backslash P'$.

\begin{lem} \label{lem: holonomy gluing non-separating}
    There exists a unique isometry $g\in \Isom_+(\H)$ such that
    \[\dev'\circ \tilde{f} = g\circ \dev'_{\vert \tilde I'}\]
    and a unique developing pair $(\dev, \rho)$ of the hyperbolic structure on $S\backslash P$ such that
    \begin{itemize}
        \item $\dev \circ \tilde{\iota} = \dev'$,
        \item $\rho \circ \iota_* = \rho'$ and $\rho(\gamma) = g$.
    \end{itemize}
\end{lem}

\begin{proof}
Choose a developing pair $(\dev,\rho)$ for the hyperbolic structure
on $S\setminus P$. By uniqueness of developing pairs, we may
normalize it so that
\[
\dev\circ\tilde\iota=\dev',
\qquad
\rho\circ\iota_*=\rho'.
\]
Using
$\tilde\iota\circ\tilde f
=\gamma\circ\tilde\iota|_{\tilde I'}$
and equivariance, we obtain
\[
\dev'\circ\tilde f
=\dev\circ\tilde\iota\circ\tilde f
=\dev\circ\gamma\circ\tilde\iota|_{\tilde I'}
=\rho(\gamma)\circ\dev'|_{\tilde I'}.
\]
Hence $g=\rho(\gamma)$. This isometry is uniquely determined,
since an orientation-preserving isometry fixing pointwise the
non-constant curve $\dev'(\tilde I')$ is the identity.
The normalization also uniquely determines the developing pair.
\end{proof}

\begin{rmk}
When one endpoint of the cutting interval $I$ is a singular point $p$, one can remove a small open disc $D$ around $p$ to obtain a singular hyperbolic metric on the complement $S\setminus D$, in which the interval $I \setminus (I\cap D)$ is now properly embedded. Lemmas \ref{lem: holonomy gluing separating case} and \ref{lem: holonomy gluing non-separating} then describe the holonomy representation of the hyperbolic structure on $S \setminus (P\cup D)$, hence on $S\setminus P$ (since the inclusion $S \setminus (P\cup D) \hookrightarrow S\setminus P$ is a homotopy equivalence). The same applies when the two endpoints of $I$ are singular points. (One needs to remove a small disc around each endpoint).

When one endpoint of the cutting interval $I$ is regular or a branch point $p$, one can similarly remove a small open disc $D$ around it. The argument as above will describe the holonomy $\rho$ of the hyperbolic structure on $S\setminus (D\cup P)$. Now, because $p$ is regular or branched, we have $\rho'(\partial D) = \Id$ and we deduce that $\rho'$ factors to a representation $\bar{\rho}: \pi_1(S\setminus P)\to \Isom_+(\H)$. This is the holonomy representation of the singular hyperbolic metric we started with.

Finally, if $\INT(I)$ contains singular points, one can decompose $I$ into finitely many intervals whose interiors do not contain singular points, and apply interatively the above arguments to infer the holonomy of $S$ from that of $S'$.
\end{rmk}

\section{Singular metrics and Teichm\"uller theory}

In this section, $S$ denotes a closed oriented connected topological surface of negative Euler characteristic and $\rho : \pi_1(S) \to \Isom_+(\H)$ an irreducible representation. As seen in Section~\ref{ss:hyperbolic structures and holomorphic maps}, finding a branched hyperbolic metric on $S$ with holonomy $\rho$ is equivalent to finding a complex structure $J$ on $S$ for which there exists a non-constant $\rho$-equivariant holomorphic map $\dev: (\tilde S, J) \to \H$.

In this section, we develop a classical approach to this problem, dating back at least to the works of Sacks--Uhlenbeck \cite{SacksUhlenbeck82} and Schoen--Yau \cite{SchoenYau79}. For every complex structure $J$ on $S$, there is a unique $\rho$-equivariant \emph{harmonic map} $f_J: (\tilde S, J) \to \H$, which is conformal exactly when $J$ is a critical point of some \emph{energy functional} $\EE_\rho$ on the Teichm\"uller space of $S$.

\subsection{Harmonic maps and the energy functional on Teichm\"uller space}

Recall that the Teichm\"uller space of $S$, denoted $\Teich(S)$, is the space of isotopy classes of complex structures on $S$ compatible with the orientation. It is a complex manifold of dimension $3g-3$, biholomorphic to a contractible domain in $\C^{3g-3}$ (see \cite[Chapter VI]{Ahlfors06}).

Given a complex structure $J$ on $S$, let us denote by $\sigma_J$ the uniformizing metric (i.e. the unique $J$-conformal metric of curvature $-1$) and $\vol_J$ the associated volume form. The hyperbolic surface $(S,\sigma_J)$ is isometric to $j(\pi_1(S))\backslash \H$, where $j: \pi_1(S) \to \Isom_+(\H)$ is a \emph{Fuchsian} representation. In the sequel, we will denote points in $\Teich(S)$ alternatively by $J$ or $j$ depending on whether we want to see them as isotopy classes of complex structures or conjugacy classes of Fuchsian representations.\\

Let $\sigma$ be a conformal metric on $(S,J)$ and $f:\tilde S \to \H$ be a smooth $\rho$-equivariant map. Then the pull-back of $\sigma_\H$ by $f$ factors to a symmetric $2$-tensor on $S$, which decomposes as
\[ f^*\sigma_\H = \phi_f + e_\sigma(f) \sigma + \bar \phi_f~,\]
where $e_\sigma(f)$ is a smooth function called the \emph{energy density of $f$} (with respect to the metric $\sigma$), and $\phi_f$ is a $(2,0)$-form on $(S,J)$ called the \emph{Hopf differential}.

If $\sigma= e^u \sigma_J$, then we have $e_\sigma(f) = e^{-u} e_{\sigma_J}(f)$ and $\vol_\sigma = e^{u} \vol_J$. The product $e_\sigma(f) \vol_\sigma$ thus only depends on the conformal structure of $\sigma$.
\\

The following theorem summarizes the fundamental properties of harmonic maps to a Riemann surface.

\begin{thm}[\cite{EellsSampson64,Corlette88}]
    For every complex structure $J$ on $S$ there exists a unique smooth $\rho$-equivariant map $f_{J,\rho}: \tilde S \to \H$ which minimizes the functional
    \[\int_S e_{\sigma_J}(f) \vol_J\]
    over all smooth equivariant maps. Moreover, its Hopf differential $\phi_{J,\rho} \coloneqq \phi_{f_{J,\rho}}$ is holomorphic.
\end{thm}

The map $f_{J,\rho}$ is called \emph{harmonic}.
\begin{defi}
    The \emph{energy functional} of the representation $\rho$ is the functional 
    \[\function{\EE_\rho}{\Teich(S)}{\R_+}{J}{\int_S  e_{\sigma_J}(f_{J,\rho}) \vol_J~.}
    \]
\end{defi}

\begin{rmk}
    If $j$ is the Fuchsian representation associated to the complex structure $J$, we will set
    \[f_{j,\rho} \coloneqq f_{J,\rho} \circ u^{-1}~,\]
    where $u = f_{J,j}: \tilde S \to \H$ is the unique $j$-equivariant biholomorphism. The map $f_{j,\rho}$ is a $(j,\rho)$-equivariant harmonic map from $\H$ to $\H$. When convenient, we will also write $\phi_{j,\rho}$ for $u_*\phi_{J,\rho}$, $e_{j,\rho}$ for $e_{\sigma_\H}(f_{j,\rho})$ and $\EE_\rho(j)$ for $\EE_\rho(J)$.
\end{rmk}

Recall that, since $S$ is closed, for any complex structure $J$ on $S$ we have
\[\Vol(\rho) = \int_S f_{J,\rho}^* \vol_\H~.\]
An elementary computation then shows that
\[\EE_\rho \geq  \vert \Vol(\rho)\vert~.\]

Recall that the holomorphic cotangent space to $\Teich(S)$ at a point $J$ naturally identifies with the space $\mathcal Q(S,J)$ of holomorphic quadratic differentials on $(S,J)$. In particular, one can see the map $J\mapsto \phi_{J,\rho}$ as a section of the holomorphic cotangent bundle to $\Teich(S)$. A classical result, dating back to Uhlenbeck and Yau (see \cite{Wentworth07}), states that $\phi_{J,\rho}$ is precisely (up to a multiplicative constant), the $(1,0)$ part of the differential of $\EE_\rho$ at $J$. As a consequence, we get the following:

\begin{thm}[see \cite{Wentworth07}] \label{t: Critical points energy}
    For $J \in \Teich(S)$, the following are equivalent:
    \begin{itemize}
    \item the harmonic map $f_{J,\rho}$ is conformal (hence holomorphic or antiholomorphic),
    \item the Hopf differential $\phi_{J,\rho}$ vanishes identically,
    \item $J$ is a critical point of the functional $\EE_\rho$,
    \item $\EE_\rho(J) = \vert \Vol(\rho)\vert$.
    \end{itemize}
\end{thm}

Hence, to prove that $\rho$ is the holonomy of a branched hyperbolic metric, it suffices to find a critical point of $\EE_\rho$ (which will necessarily be a global minimum). This is however far from easy because $\EE_\rho$ is never proper unless $\rho$ is Fuchsian. 

\begin{prop} \label{prop: Energy not proper}
    Let $\rho:\pi_1(S)\to \Isom_+(\H)$ be an irreducible representation of a closed surface group. Then $\EE_\rho$ is proper if and only if $\rho$ is Fuchsian.
\end{prop}

\begin{proof}
    Properness of the energy functional of a Fuchsian representation is classical (see for instance \cite{Tromba92}). Conversely, assume that $\EE_\rho$ is proper. Then the critical locus of $\EE_\rho$ is non-empty.
    
    By Theorem \ref{t: Critical points energy}, this critical locus is the set of complex structures $J$ for which $f_{J,\rho}$ is holomorphic or antiholomorphic. This is known to be a complex subvariety of $\Teich(S)$ of dimension $2g-2 - \vert \euler(\rho)\vert$ (see for instance the appendix of \cite{CDF14}). Since $\EE_\rho = \vert \Vol(\rho)$ on this locus, properness of $\EE_\rho$ implies that this locus is compact. Since $\Teich(S)$ embeds as a bounded domain in $\C^{3g-3}$, this locus is finite. Hence $\vert \euler(\rho)\vert = 2g-2$ and $\rho$ is Fuchsian.
\end{proof}

\begin{rmk}
    Since reducible representations have Euler class $0$, we don't need to consider them in this paper. However, it is still possible to define their energy functional, and prove that it is non-proper: at a minimum, the energy would vanish, implying that the semisimple part of the representation fixes a point in $\H$. This in turn implies that the energy functional is identically $0$.
\end{rmk}

\subsection{Stable nodal curves and the compactified moduli space}

Since a minimizing sequence for $\EE_\rho$ may diverge, we need a geometric understanding of how a sequence of Riemann surfaces degenerates. 
Such a description is provided by the Deligne--Mumford compactification of the moduli space of curves.

A \textit{stable curve} of genus $g$ is a compact connected complex curve $C$ of
arithmetic genus $g$ whose only singularities are nodes, and such that every
connected component of $C\setminus \text{Sing}(C)$ has negative Euler
characteristic \cite{DM69}, \cite[Chap.~X, Section 3]{ACG11}. By the uniformization
theorem, $C\setminus\text{Sing}(C)$ then carries a unique complete
hyperbolic metric of finite area compatible with its complex structure, which
we denote by $\text{hyp}_C$; each node corresponds to a pair of cusps
of this hyperbolic surface. 

\begin{thm}[Deligne-Mumford \cite{DM69}]\label{thm:DM}
    The moduli space $\overline{\mathcal M}_g$  of
stable curves of genus $g$ can be endowed with a complex orbifold structure, which contains the moduli space
$\mathcal M_g=\mathcal T(S)/\mathrm{MCG}(S)$ of smooth curves as an open
dense subset. Moreover,  $\overline{\mathcal M}_g$ is compact. 
\end{thm}

We refer to \cite[Chap. XII, Section 2]{ACG11} for the construction of
$\overline{\mathcal M}_g$ as an analytic space; compactness is a consequence
of the stable reduction theorem \cite[Chap. X, Theorem (4.11)]{ACG11}.

We will use the following local description of $\overline{\mathcal M}_g$ near
a stable curve $C_\infty$ with nodes $y_1,\dots,y_K$. There is a holomorphic
family $\pi:\mathcal C\to B$ of stable curves over a polydisc $B\ni 0$ (a
Kuranishi family) with $\pi^{-1}(0)=C_\infty$, such that
$B/\mathrm{Aut}(C_\infty)$ is homeomorphic to a neighbourhood of $[C_\infty]$
in $\overline{\mathcal M}_g$, the point $t$ corresponding to the fibre
$C_t=\pi^{-1}(t)$ \cite[Chap. XI and Chap.~XII, Section 2]{ACG11}. Near each node
$y_k$, the family is given in suitable coordinates by $\{xy=t_k\}$. 
On compact subsets away from the nodes of $C_\infty$, the family admits
smooth local trivializations for small $t$. There are also
\emph{collapse maps} $c_t:C_t\to C_\infty$
\cite[Chap. X, Section 9]{ACG11} with the following properties:
\begin{itemize}
\item $c_t$ is continuous and surjective, and restricts to a diffeomorphism
from $C_t\setminus c_t^{-1}(\{y_1,\dots,y_K\})$ onto
$C_\infty\setminus \mathrm{Sing}(C_\infty)$,
\item if $C_t$ is smooth, each $c_t^{-1}(y_k)$ is a simple closed curve (the
\emph{vanishing cycle} associated to $y_k$),
\item as $t\to0$, the maps $c_t^{-1}$ converge to the identity in the
$\mathcal C^\infty$ topology on compact subsets of
$C_\infty\setminus\mathrm{Sing}(C_\infty)$.
\end{itemize}
The key analytic input is the continuity of the hyperbolic metric in such a
family, due to Wolpert.

\begin{thm}[Wolpert \cite{Wol90}]\label{thm:Wolpert}
The hyperbolic metrics $\mathrm{hyp}_{C_t}$ define a continuous metric
on the vertical tangent bundle $T_{\mathcal C^{\mathrm{sm}}/B}$, where
\[
\mathcal C^{\mathrm{sm}}
=\{x\in\mathcal C:x\text{ is a smooth point of }C_{\pi(x)}\}
\]
is the relative smooth locus of the family. In particular,
$(c_t^{-1})^*\mathrm{hyp}_{C_t}$ converges to $\mathrm{hyp}_{C_\infty}$
uniformly on compact subsets of $C_\infty\setminus\mathrm{Sing}(C_\infty)$
as $t\to0$.
\end{thm}

See also \cite[Chap. XVIII, Lemma (6.9) and Section 11]{ACG11} and \cite{HK14}.
Wolpert moreover gives precise expansions of $\mathrm{hyp}_{C_t}$ near the
nodes. In particular the hyperbolic lengths of the vanishing cycles tend to
$0$ as $t\to0$.

To interpret this convergence in terms of holonomies, we will rely on the following elementary lemma:

\begin{lem}\label{lem:metrics-to-holonomies}
Let $S$ be a connected surface of finite type and let $(\sigma_n)_{n\in \N}$ be a
sequence of (possibly incomplete) hyperbolic metrics on $S$ converging
uniformly on compact subsets of $S$ to a hyperbolic metric $\sigma_\infty$.
Denote by $\rho_n$ and $\rho_\infty$ their holonomies. Then, up to extracting
a subsequence, there exist $(g_n)_{n\in \N}\in\mathrm{Isom}^+(\mathbb H)^\N$ such that
$\mathrm{Ad}_{g_n}\circ\rho_n$ converges to a representation conjugate to
$\rho_\infty$.
\end{lem}

\begin{proof}
Let $\dev_n$ be the $\rho_n$-equivariant developing map for $\sigma_n$ and $g_n$ an isometry such that $g_n\circ \dev_n(\tilde x_0)=o$.
Uniform convergence of the metrics on compact sets gives local uniform
Lipschitz bounds of $g_n \circ \dev_n$ with respect to a fixed metric and bounds the images of compact sets. By Ascoli's theorem, up to extracting a subsequence, $g_n \circ \dev_n$ converges locally uniformly to a map $\dev'_\infty$.

For each element $\gamma$ of a finite generating set of  $\pi_1(S)$,
equivariance gives
$\mathrm{Ad}_{g_n}\rho_n(\gamma)o=g_n \circ \dev_n(\gamma\tilde x_0)$.
Since these points remain bounded, after further extracting a subsequence, we can assume that $\Ad_{g_n} \circ \rho_n$ converges to a representation $\rho'_\infty$. Passing equivariance to the limit, we get that $\dev'_\infty$ is $\rho'_\infty$-equivariant.

We claim that $\rho'_\infty$ is conjugate to $\rho_\infty$. Indeed, $g_n\circ \dev_n :(\tilde{S},\sigma_n) \to \H$ is locally isometric by the definition of a developing map, and this property passes to local uniform limits. Hence $\dev'_\infty: (\tilde{S},\sigma_\infty)\to \H$ is locally isometric, and $(\dev'_\infty,\rho'_\infty)$ is a developing pair for $\sigma_\infty$. By essential uniqueness of the holonomy of a hyperbolic structure, we deduce that $\rho'_\infty$ is conjugate to $\rho_\infty$.
\end{proof}

We deduce the following reformulation of Deligne--Mumford's compactification in terms of surface group representations.
\begin{cor}\label{coro: Fuchsian decompoisition Deligne Mumford}
Let $(j_n)_{n\in\mathbb N}$ be a sequence of Fuchsian representations of
$\pi_1(S)$. Then, up to extracting a subsequence, there exist a decomposition
\[
S=\bigsqcup_{1\le i\le I}S_i\sqcup\bigsqcup_{1\le l\le L}b_l ,
\]
a sequence $(\tau_n)_{n\in\mathbb N}$ in $\mathrm{MCG}(S,*)$, sequences
$(g^i_n)_{1\le i\le I,\,n\in\mathbb N}$ in $\mathrm{Isom}^+(\mathbb H)$ and
representations $j^i_\infty:\pi_1(S_i)\to\mathrm{Isom}^+(\mathbb H)$ such that
for all $i\le I$,
\[
(\mathrm{Ad}_{g^i_n}\circ j_n\circ\tau_n)|_{\pi_1(S_i)}
\xrightarrow[n\to+\infty]{} j^i_\infty ,
\]
where $j^i_\infty$ is the holonomy of a complete hyperbolic metric of finite
area on $S_i$.

Moreover, for every $R>0$ there exists $N$ such that, for every $\gamma \in \pi_1(S)$, and every $n\geq N$, we have
\[d(o,\Ad_{g_n^i}\circ j_n \circ \tau_n(\gamma)\cdot o) \leq R \Rightarrow \gamma \in \pi_1(S_i)~.\]
\end{cor}

\begin{proof}
Let $X_n=(S,J_n)$ be the Riemann surface uniformized by $j_n$, and set
$\sigma_{j_n}=\mathrm{hyp}_{X_n}$. By Theorem~\ref{thm:DM}, up to extracting a
subsequence, $X_n$ converges in $\overline{\mathcal M}_g$ to the class of a
stable curve $C_\infty$, with nodes $y_1,\dots,y_L$. Let
$\pi:\mathcal C\to B$ be a Kuranishi family of $C_\infty$ as above. Up to
extracting again, there are $t_n\to0$ in $B$ and biholomorphisms
$h_n:X_n\to C_{t_n}$.

 Fix $n_0$ and put $c=c_{t_{n_0}}\circ h_{n_0}:S\to
C_\infty$. Let $b_l=c^{-1}(y_l)$, and let $S_i$, $1\le i\le I$, be the
connected components of $S\setminus\bigcup_l b_l$. The map $c$ sends each
$S_i$ diffeomorphically onto a connected component $C^i_\infty$ of
$C_\infty\setminus\mathrm{Sing}(C_\infty)$. By stability,
$\chi(S_i)=\chi(C^i_\infty)<0$, so $S=\bigsqcup S_i\sqcup\bigsqcup b_l$ is a
decomposition in the sense of Definition \ref{d: decomposition} (see Remark \ref{r: essential curve}).

 The map $c_{t_n}^{-1}\circ c$ is a
diffeomorphism from $S\setminus\bigcup_l b_l$ onto the complement of the
vanishing cycles in $C_{t_n}$. Modifying it in collars of the curves $b_l$
which shrink as $n\to\infty$, we obtain orientation-preserving
diffeomorphisms $\psi_n:S\to C_{t_n}$ which coincide with $c_{t_n}^{-1}\circ c$
outside these collars. Choosing a lift of the diffeomorphism $h_n^{-1}\circ\psi_n$ of $S$ to $\widetilde S$, we obtain an element
$\tau_n\in\mathrm{MCG}(S,*)$. Since $h_n$ is an isometry for the
hyperbolic metrics, $j_n\circ\tau_n$ is the holonomy of
$\psi_n^*\mathrm{hyp}_{C_{t_n}}$.

By Theorem~\ref{thm:Wolpert}, the metrics
$\psi_n^*\mathrm{hyp}_{C_{t_n}}$ converge uniformly on compact subsets of
$S\setminus\bigcup_lb_l$ to $c^*\mathrm{hyp}_{C_\infty}$, whose restriction to
each $S_i$ is the complete finite area hyperbolic metric
$\sigma^i_\infty=(c|_{S_i})^*\mathrm{hyp}_{C^i_\infty}$. The holonomy of the
restriction of $\psi_n^*\mathrm{hyp}_{C_{t_n}}$ to $S_i$ is
$(j_n\circ\tau_n)|_{\pi_1(S_i)}$. Applying
Lemma~\ref{lem:metrics-to-holonomies} to $S_i$ gives, after extraction,
elements $g^i_n$ such that $(\mathrm{Ad}_{g^i_n}\circ j_n\circ\tau_n)|_{\pi_1(S_i)}$
converges to a holonomy representation $j^i_\infty$ of $\sigma^i_\infty$.

To prove the last part of Corollary \ref{coro: Fuchsian decompoisition Deligne Mumford}, let $x_0$ be the projection to $S_i$ of $\dev_\infty^{-1}(o)$. Since the metric $\sigma_\infty$ on $S_i$ is complete, we get that
\[d_{\sigma_n}(x_0, \partial S_i) \underset{n\to +\infty}{\longrightarrow} +\infty~.\]
Hence, for $n$ large enough, if $d(o, \Ad_{g_n^i}\circ j_n \circ \tau_n(\gamma)\cdot o) \leq R$, then the geodesic loop for $\sigma_n$ based at $x_0$ and homotopic to $\gamma$ remains in $S_i$, hence $\gamma \in \pi_1(S_i)$.
\end{proof}

\subsection{Harmonic maps along a minimizing sequence}

Throughout this section, we fix a closed oriented connected surface $S$ of genus $g \geq 2$ and a representation $\rho: \pi_1(S) \to \Isom_+(\H)$ such that
\[ 0 < \Vol(\rho) < 2\pi(2g-2)~.\]
We see $\Teich(S)$ as the space of conjugacy classes of Fuchsian representations.

Our general goal is to find a minimizer of $\EE_\rho$ (see Theorem \ref{t: Critical points energy}). Since $\EE_\rho$ is not proper (Proposition \ref{prop: Energy not proper}), a minimizing sequence $(j_n)_{n\in \N}$ for $\EE_\rho$ may not converge. However, it will converge -- up to the action of the mapping class group of $S$ -- to a point in the Deligne--Mumford compactification of the moduli space of $S$. Keeping track of the behaviour of the harmonic maps $f_{j_n,\rho}$ along this sequence, we will prove the following lemma:

\begin{lem} \label{lem: Conformal limits harmonic maps minimizing sequence}
Up to extracting a subsequence, there exists a decomposition
\[S= \bigsqcup_{1 \leq i \leq I} S_i \sqcup \bigsqcup_{1\leq l \leq L} b_l~,\]
a sequence $(\tau_n)_{n\in \N}\in \MCG(S,*)^\N$, sequences $(g^i_n)_{1\leq i \leq I, n\in \N}$ in $\Isom_+(\H)$ and a family $(j^i_\infty)_{1\leq i \leq I}$ of finite volume Fuchsian representations of $\pi_1(S_i)$ such that for each $i$,
\begin{itemize}
\item $\rho^i_n \coloneqq {\Ad_{g^i_n} \circ \rho \circ \tau_n}_{|\pi_1(S_i)}$
converges to a representation
\[\rho^i_\infty: \pi_1(S_i) \to \Isom_+(\H)~,\]
\item there exists a $(j^i_\infty, \rho^i_\infty)$-equivariant holomorphic or antiholomorphic map from $\H$ to $\H$.
\end{itemize}
\end{lem}

To prove this lemma, let us consider a sequence $j_n$ of points in $\Teich(S)$ such that 
\[\EE_\rho(j_n) \underset{n\to + \infty}{\to} \inf_{j\in \Teich(S)} \EE_\rho(j)~.\]
Then, by Corollary \ref{coro: Fuchsian decompoisition Deligne Mumford}, up to taking a subsequence, there exists a decomposition \[S= \bigsqcup_{1\leq i \leq I} S_i \sqcup \bigsqcup_{1\leq l \leq L} b_l~,\]
a sequence $(\tau_n)_{n\in \N}$ of elements of $\MCG(S,*)$, sequences $(h^i_n)_{1\leq i\leq I,\,n\in\N}$ in $\Isom_+(\H)$ and a family $(j^i_\infty)_{1\leq i\leq I}$ of finite-covolume Fuchsian representations of $\pi_1(S_i)$ such that for each $i$,
\[{\Ad_{h^i_n} \circ j_n \circ \tau_n}_{|\pi_1(S_i)} \underset{n\to +\infty}{\longrightarrow} j^i_\infty~.\]
To simplify notation, let us set
\[j^i_n \coloneqq \Ad_{h^i_n} \circ j_n \circ \tau_n~.\]

Fix a base point $o$ in $\H$. Define 
\[f^i_n = g^i_n \circ f_{j_n,\rho} \circ {h^i_n}^{-1}~,\]
where $g_n^i$ is chosen so that $f_n^i(o)=o$. Set
\[
\rho_n^i\coloneqq
\Ad_{g_n^i}\circ\rho\circ\tau_n~.
\]
Then $f_n^i$ is $(j_n^i,\rho_n^i)$-equivariant and harmonic.

\begin{prop} \label{p:HolomorphicLimit}
Up to extracting a subsequence, each $\rho_n^i|_{\pi_1(S_i)}$ converges to a representation $\rho_\infty^i$ and each $f^i_n$ converges in the $\mathcal C^2$ topology on compact subsets to a harmonic map $f^i_\infty$ which is $(j^i_\infty, \rho^i_\infty)$-equivariant.
\end{prop}

The proof will make use of the following regularity estimates for harmonic maps:

\begin{prop}\label{prop: Uniform estimates harmonic maps}
For any relatively compact domains $V\Subset U\Subset\H$, there exists
$C(V,U)>0$ such that every harmonic map $f:\H\to\H$ satisfies
\[
    \sup_V e_{\sigma_\H}(f)\leq C(V,U)\EE_U(f)~,
\textrm{ where }
    \EE_U(f)=\int_U e_{\sigma_\H}(f)\,\vol_\H~.
\]
Consequently, any sequence of harmonic maps $f_n:\H\to\H$ with
$f_n(o)=o$ and uniformly bounded energy on each relatively compact domain has a subsequence converging in $\mathcal C^\infty$ topology on compact subsets to a harmonic map.
\end{prop}
While this result is classical, we give the outline of the proof for completeness.

\begin{proof}
The Bochner formula gives $\Delta e(f)\geq-2e(f)$, so the interior mean-value estimate gives the stated inequality. It yields local
uniform bounds on $\d f_n$; together with $f_n(o)=o$, these keep the
images of each compact set in a fixed compact subset of $\H$.
Interior elliptic regularity for the harmonic map equation then bounds all derivatives on smaller domains, and a diagonal extraction gives the conclusion.
\end{proof}

\begin{proof}[Proof of Proposition \ref{p:HolomorphicLimit}]Let $U$ be a bounded open domain in $\H$. Let $k$ be the cardinality of the set\[\{\gamma \in \pi_1(S_i) \mid j^i_\infty(\gamma) \cdot \overline{U} \cap \overline{U} \neq \emptyset\}~.\]Then, for $n$ large enough, one still has\[\sharp \{\gamma \in \pi_1(S) \mid j^i_n(\gamma) \cdot \overline{U} \cap \overline{U} \neq \emptyset\} = \sharp \{\gamma \in \pi_1(S_i) \mid j^i_n(\gamma) \cdot \overline{U} \cap \overline{U} \neq \emptyset\} \leq k~.\]
(The first equality comes from the last statement of Corollary \ref{coro: Fuchsian decompoisition Deligne Mumford}, and the second from the geometric convergence of the Fuchsian groups $j^i_n(\pi_1(S_i))$ to $j^i_\infty(\pi_1(S_i))$).

Thus, for all $n$ large enough, \[\EE_U(f^i_n) \leq k \EE(j^i_n, \rho^i_n) = k \EE(j_n,\rho) \leq k (\inf_{j \in \Teich(S)} \EE(j,\rho) +1)~~\]
which is bounded. Since, moreover, $f^i_n(o) = o$ for all $n$, by Proposition \ref{prop: Uniform estimates harmonic maps} and Ascoli's theorem, one can extract a subsequence of the $(f^i_n)_{n\in \N}$ that converges uniformly on every compact subset $K\subset\H$ in $\mathcal{C}^2$ topology. With a diagonal extraction, one can thus assume that $f^i_n$ converges uniformly on every compact set in $\mathcal{C}^2$ topology, to some map $f^i_\infty: \H \to \H$ which must then be harmonic. In particular, for every $\gamma \in \pi_1(S_i)$, the sequence $\rho^i_n (\gamma) \cdot o = f^i_n(j^i_n(\gamma) \cdot o)$ is bounded. Upon extracting a further subsequence, we can thus assume that $\rho^i_n|_{\pi_1(S_i)}$ converges to some representation $\rho^i_\infty: \pi_1(S_i) \to \Isom_+(\H)$ such that $f^i_\infty$ is $(j^i_\infty, \rho^i_\infty)$-equivariant.\end{proof}

It remains to prove that the maps $f^i_\infty$ are in fact conformal, i.e. have vanishing Hopf differential. To do so, we will prove the next lemma, which needs the following notation. For a quadratic differential $\phi=q\,dz^2$ and a conformal metric
$\sigma=\lambda^2|dz|^2$, we set
\[
    \norm{\phi}_\sigma=\frac{|q|}{\lambda^2}.
\]
On a hyperbolic surface, the subscript $\H$ refers to its hyperbolic
metric, and $\vol_\H$ denotes the associated area form.

\begin{lem}\label{lem: Decreasing the energy with Hopf}
Let $j:\pi_1(S)\to\Isom_+(\H)$ be a Fuchsian representation and let
$X=j(\pi_1(S))\backslash\H$. Suppose that $D(x,2R)$ is an embedded
hyperbolic disc in $X$ and that
$\sup_{D(x,2R)}e_{j,\rho}\leq m$ for some $m>0$.
Then there exists a Fuchsian representation $j'$ such that
\[
    \EE_\rho(j')\leq\EE_\rho(j)
    -\frac{1}{m}\int_{D(x,R)}
        \norm{\phi_{j,\rho}}_\H^2\,\vol_\H.
\]
\end{lem}

\begin{proof}
To simplify notation, we write $\sigma$ for the hyperbolic metric on $X$, and set
$f=f_{j,\rho}$, $e=e_{j,\rho}$, $\phi=\phi_{j,\rho}$ and
$p=\norm{\phi}_\sigma$.
The identity $f^*\sigma_\H=e\sigma+\phi+\bar\phi$ implies $2p\leq e$.
Choose a smooth function $\beta$, supported in $D(x,2R)$, such that
$0\leq\beta\leq \frac{1}{2m}$ and $\beta=\frac{1}{2m}$ on $D(x,R)$.
Since $2\beta p\leq\beta e\leq\frac{1}{2}$, the tensor
\[
    \sigma'=\sigma+\beta(\phi+\bar\phi)
\]
is a positive definite metric. A direct computation gives
\[
    e_{\sigma'}(f)\,\vol_{\sigma'}
    =\frac{e-4\beta p^2}{s}\,\vol_\sigma~,
    \textrm{ where }  s=\sqrt{1-4\beta^2p^2}\in(0,1].
\]
Using $\beta e\leq 1$, we obtain
\[
\begin{aligned}
    e-\frac{e-4\beta p^2}{s}
    &=\frac{4\beta p^2}{s}
        \left(1-\frac{\beta e}{1+s}\right)\\
    &\geq\frac{4\beta p^2}{1+s}
    \geq2\beta p^2.
\end{aligned}
\]
Let $j'$ uniformize the marked conformal structure determined by $\sigma'$.
By conformal invariance of energy and the minimizing property of the
$\rho$-equivariant harmonic map, we have
\[
\begin{aligned}
    \EE_\rho(j')
    &\leq\int_X e_{\sigma'}(f)\,\vol_{\sigma'}\\
    &\leq\EE_\rho(j)-2\int_X\beta p^2\,\vol_\sigma\\
    &\leq\EE_\rho(j)
        -\frac{1}{m}\int_{D(x,R)}p^2\,\vol_\sigma.
\end{aligned}
\]
\end{proof}

We can now conclude the proof of
Lemma~\ref{lem: Conformal limits harmonic maps minimizing sequence}.

\begin{proof}[Proof of Lemma \ref{lem: Conformal limits harmonic maps minimizing sequence}]To conclude the proof, it remains to show that the maps $f^i_\infty$ constructed above are conformal. Let $\phi^i_\infty$ denote the Hopf differential of $f^i_\infty$. 

Let $x$ be a point in $\H$. By convergence of $j_n^i|_{\pi_1(S_i)}$ to $j_\infty^i$ and the last statement of Corollary \ref{coro: Fuchsian decompoisition Deligne Mumford}, there exists $R>0$ such that, for all sufficiently large $n$ and all $\gamma \in \pi_1(S) \setminus \{\Id\}$, we have\[j_n^i\gamma \cdot D(x,2R)\cap D(x,2R) = \emptyset~.\]
In other words, the injectivity radius of $x$ in $j_n^i(\pi_1(S))\backslash \H$ is uniformly bounded from below by $2R$. Since the maps $f^i_n$ converge in the $\mathcal C^2$ topology to $f^i_\infty$ uniformly on $D(x,2R)$, there exists a constant $m$ such that \[e_\H(f^i_n) \leq m\]for all $n$. By Lemma \ref{lem: Decreasing the energy with Hopf}, we can find a sequence of Fuchsian representations $j'_n$ such that
\[\int_{D(x,R)} \norm{\phi^i_n}^2_\H\,\vol_\H\leq m (\EE_{\rho_n^i}(j_n^i) - \EE_{\rho_n^i}(j'_n \circ \tau_n))= m (\EE_\rho(j_n)- \EE_\rho(j_n'))~.\]
On the other hand, we have\[\EE_\rho(j_n)- \EE_\rho(j'_n) \leq \EE_\rho(j_n) - \inf_j \EE_\rho(j) \underset{n\to +\infty}{\longrightarrow} 0\]since $(j_n)$ is a minimizing sequence for $\EE_\rho$. Passing to the limit, we obtain that
\[\int_{D(x,R)} \norm{\phi^i_\infty}^2_\H\,\vol_\H= 0~.\]
Hence $\phi^i_\infty$ vanishes in a neighbourhood of every point, and we conclude that $f^i_\infty$ is conformal.\end{proof}

We conclude this section with the following property of the volume along this minimizing sequence:

\begin{prop} \label{prop: Additivity volume limit}
    With the notation above, we have
    \[\Vol(\rho) = \sum_{1\leq i \leq I}  \Vol(\rho^i_\infty)~.\]
\end{prop}

\begin{proof}
    Since the volume is additive (Proposition \ref{prop: Volume Additive}) and invariant under conjugation and the mapping class group (Proposition \ref{prop: volume MCG invariant}), we have
    \begin{eqnarray*}
        \Vol(\rho) & = & \Vol(\rho\circ \tau_n) \\
        & = & \sum_{1\leq i \leq I} \Vol({\rho\circ \tau_n}_{\vert S_i})\\
        &=& \sum_{1\leq i \leq I} \Vol(\rho_n^i)~.
    \end{eqnarray*}
    Finally, the continuity of the volume (Proposition \ref{prop: Volume continuous}) gives
    \[\lim_{n\to+\infty} \Vol(\rho_n^i) = \Vol(\rho_\infty^i)~.\]
\end{proof}

\section{Building hyperbolic surfaces from hyperbolic blocks}

In this section we introduce a notion of \emph{hyperbolic block}, which is well-adapted to cutting and gluing constructions: informally, these are branched hyperbolic surfaces with standardized behaviour on the boundary. The main results of the section are a deformation theorem \`a la Ehresmann--Thurston for hyperbolic blocks (Lemma \ref{l:perturbation of hyperbolic block}) and a result stating that, under suitable conditions, hyperbolic blocks can be pieced together to form a hyperbolic surface with prescribed holonomy.

\subsection{An Ehresmann--Thurston principle for oriented branched hyperbolic surfaces with boundary}

The so-called Ehresmann--Thurston principle (see for instance \cite{BergeronGelander04}) states that, given a $(G,X)$-structure on a closed manifold $M$ with holonomy $\rho$, any representation $\rho':\pi_1(M) \to G$ sufficiently close to $\rho$ is again the holonomy of another (essentially unique) $(G,X)$-structure on $M$. Several extensions of this principle are known for singular geometric structures and non-compact manifolds. In particular, the Ehresmann--Thurston principle holds for hyperbolic metrics with branch points and cone singularities (though the uniqueness part fails in the presence of branch points).

Before introducing hyperbolic blocks, let us state a general Ehresmann--Thurston principle for singular hyperbolic structures  with smooth boundary.

\begin{lem}\label{l:ehresman-thurston}
Let $S$ be a compact connected oriented surface, possibly with boundary components $(c_k)_{1\leq k\leq K}$, and let $\sigma_0$ be a branched
hyperbolic metric on $S$, with smooth regular boundary and all its
branch points in the interior. Fix a developing pair
$(\dev_0,\rho_0)$ and lifts $\tilde c_k$ of the boundary components.
There exist a neighbourhood $U$ of $\rho_0$ in
$\Hom(\pi_1(S),\Isom_+(\H))$ and $\epsilon>0$ such that, for every
$\rho\in U$ and every collection of smooth $\rho(c_k)$-equivariant maps
$d_k:\tilde c_k\to\H$ that are $\epsilon$-close to
$\dev_0|_{\tilde c_k}$ in the $\mathcal C^2$ topology on fixed compact
fundamental intervals, there is a branched hyperbolic structure on $S$
with developing pair $(\dev,\rho)$ satisfying
\[
\dev|_{\tilde c_k}=d_k\qquad(1\leq k\leq K).
\]
The branch points and their orders can be kept unchanged.
\end{lem}

\begin{proof}
We begin by reviewing the Ehresmann--Thurston transversality argument, based on the following observation: if $\rho:\pi_1(S)\to\text{Isom}_+(\mathbb H)$ is any representation, then a $\rho$-equivariant map defines a section $\sigma$ of the flat $\mathbb H$-bundle $E_\rho$ having monodromy $\rho$, and vice-versa. 

Note that the representation variety \(\text{Hom} (\pi_1(S) , \text{Isom}_+ (\mathbb H))\) is smooth near $\rho_0$. Indeed, in the case $S$ is a surface with boundary, its fundamental group is a free group on $s$ generators, so that $\text{Hom}(\pi_1(S), \text{Isom}_+(\mathbb H)) \simeq \text{Isom}_+(\mathbb H)^s$. If \(S\) is closed, then \(\rho_0\) is non-elementary (see for instance \cite{GKM}) so the representation variety $\text{Hom}(\pi_1(S), \text{Isom}_+(\mathbb H))$ is smooth near $\rho_0$ too. 

Let $M$ be the quotient of $U \times \tilde{S} \times \mathbb H$ by the action of $\pi_1(S)$ given by $\gamma \cdot (\rho, x, z )= (\rho, \gamma (x), \rho (\gamma) (z))$. It has a structure of an $\mathbb H$-bundle over $U\times S$ given by $(\rho, x,z) \mapsto (\rho, x)$. Smooth sections of this $\mathbb H$-bundle are in one-to-one correspondence with smoothly parameterized families $\{ \dev_\rho \}_{\rho \in U }$ where $\dev_\rho$ is a $\rho$-equivariant smooth map from $\tilde{S}$ to $\mathbb H$. 

This $\mathbb H$-bundle comes naturally with a connection defined on the slices $\rho \times S$, on which it is flat (this is the trivial connection at the level of the universal covering). The bundle over any of these slices is then exactly the flat bundle $E_\rho$.  Extend the connection to one which is defined on the whole $T(U\times S)$. Any smooth path $\{\rho_t \}_{t\in [0,1]}$ connecting $\rho_0$ to $\rho_1=\rho$ then gives a smooth identification between $E_{\rho_0}$ and $E_\rho$, sending fibers of $E_{\rho_0}$ to corresponding fibers of $E_\rho$. The section corresponding to $\dev_0$ can be transported under these identifications to define a family of sections $\sigma_t$ of $E_{\rho_t}$ which depends smoothly on $t$.  Under the above dictionary, this family of sections can be thought of as a smooth family of (smooth) maps $\dev_t: \tilde{S} \rightarrow \mathbb H$ that are $\rho_t$-equivariant. 

If $\rho$ is sufficiently close to $\rho_0$, the smooth path $\{\rho_t \}_{t\in [0,1]}$ can be chosen close to a constant path in smooth topology, and the equivariant maps $\dev_t$ are close to each other in the smooth topology on compact subsets. In particular, if $\dev_0$ is a local diffeomorphism, then $\dev_1$ remains a local diffeomorphism, since this property can be detected on a fundamental domain, which is compact.

 Finally, to prescribe the boundary behaviour of the developing map, note that the maps $d_k: \tilde c_k \to \H$ correspond to a section $\sigma^{\partial S}$ of $E_\rho$ restricted to the boundary of $S$, which is close to the section $\sigma_1$ that we just constructed. We can then write $\sigma^{\partial S} = \exp (v) \sigma_1$ where $v$ is a section of class $C^2$ of the tangent bundle of $E_\rho$ restricted to the section $\sigma_1$, $v$ being close to $0$ in the $C^2$ topology. This section can be extended to a section of class $C^2$ of the same bundle over the whole of $S$, keeping the property that it is close to $0$ in the $C^2$ topology. In particular, we get a new $C^2$ section of $E_\rho$ which is close to $\sigma_1$ (hence close to $\sigma_0$) and which coincides with $\sigma^{\partial S}$ on $\partial S$. The corresponding $\rho$-equivariant map $\dev : \tilde S \rightarrow \mathbb H $ is a new developing map with holonomy $\rho$ which satisfies $\dev|_{\tilde c_k} = d_k$ for all $1\leq k\leq K$.

These considerations establish Lemma \ref{l:ehresman-thurston} when there are no branch points. It remains to indicate how to proceed when there are some branch points in the interior of $S$. For each branch point $b$ of $\sigma_0$, let $D_b$ be a small open disc containing $b$ on which the developing map $\dev_0$ is a ramified cyclic covering of order $k_b$ onto a hyperbolic disc in $\mathbb H$. Also, choose a preimage $\tilde{b}$ of $b$ in $\tilde{S}$ and denote by $\tilde D_b$ the lift of $D_b$ containing~$\tilde{b}$. 

Now consider the surface $S' = S\setminus \cup_b D_b'$, let $\iota$ be the inclusion $S' \hookrightarrow S$ and $\tilde \iota: \tilde S' \to \tilde S$ its lift to the universal cover. Note that $\tilde \iota$ is not injective. More precisely, the kernel of $\iota_*$ is the normal subgroup generated by the loops $\partial D_b$, and the map $\tilde \iota$ factors to an embedding
\[\hat \iota: \hat S' \hookrightarrow \tilde S~,\]
where $\hat S' = \tilde S'/\ker \iota_*$.

Let $\tilde{\partial D}_b$ denote the lift of $\partial D_b$ in $\tilde S'$ such that $\tilde \iota(\tilde{\partial D}_b) = \partial \tilde D_b$. 

Let $\rho:\pi_1(S)\to \Isom_+(\H)$ be close to $\rho_0$, and $d_k: \tilde c_k \to \H$ be $\rho_{\vert c_k}$-equivariant maps close to ${\dev_0}_{\vert \tilde c_k}$. Applying the first part of the proof to $S'$, we can thus find a $\rho\circ \iota_*$-equivariant developing map
\[\dev': \tilde S' \to \H~,\]
close to $\dev_0 \circ \tilde \iota$, such that
\[\dev'_{\vert \tilde{\iota}^{-1}(\tilde c_k)} = d_k\circ \tilde\iota\]
and
\[\dev'_{\vert \tilde{\partial D}_b} = \dev_0 \circ \tilde \iota~.\]
The map $\dev'$ factors to a $\rho$-equivariant map $\dev: \hat S' \to \H$ such that
\[\dev_{\vert \gamma \cdot \partial \tilde D_b} = \rho(\gamma)\cdot {\dev_0}_{\vert \partial \tilde D_b}~.\]
We can thus extend $\dev$ to $\tilde S$ by setting
\[\dev_{\vert \gamma \cdot \tilde D_b} = \rho(\gamma)\cdot {\dev_0}_{\vert \tilde D_b}~.\]
The map $\dev$ is thus the developing map of a hyperbolic structure on $S$ with holonomy $\rho$ and with the same branch points as $\dev_0$. Moreover, it is close to $\dev_0$ and has the required behaviour on the boundary. This concludes the proof of Lemma \ref{l:ehresman-thurston}.
\end{proof}

\subsection{Hyperbolic blocks} \label{ss: Hyperbolic blocks}

We now introduce a special kind of branched hyperbolic surface with boundary, which we call \emph{hyperbolic blocks}.

\begin{defi}
Let $S$ be a compact oriented surface with boundary components $(c_k)_{1\leq k \leq K}$. A branched hyperbolic metric on $S$ with developing pair $(\dev,\rho)$ will be called a \emph{hyperbolic block} if for each $1\leq k \leq K$,
\begin{enumerate}
    \item the boundary component $c_k$ has constant curvature,
    \item $\dev_{\vert \tilde c_k}$ is  injective on a fundamental domain of $\tilde c_k$,
    \item if $\rho(c_k) = \Id$, then $c_k$ is concave.
\end{enumerate}
\end{defi}

\begin{figure}[htbp]
\centering
\includegraphics[width=0.43\textwidth]{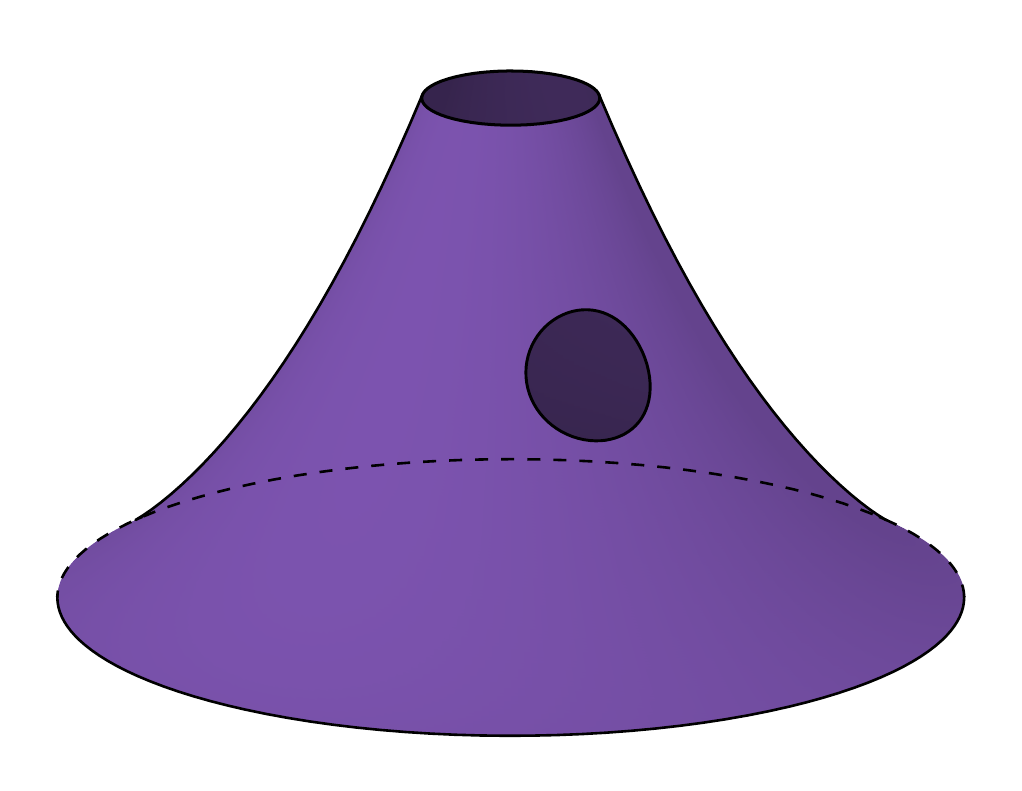}
\caption{A hyperbolic block on a punctured sphere }
\end{figure}

To clarify this definition, let us recall what curves of constant curvature in $\H$ look like:
\begin{itemize}
    \item If $\rho(c_k)$ is hyperbolic, $\dev(\tilde c_k)$ is the curve at distance $r$ from the translation axis of $\rho(c_k)$. The curve $c_k$ will be convex (resp. concave) when $\dev(\tilde c_k)$ is to the left (resp. to the right) of the axis of $\rho(c_k)$. In that case $\dev_{\vert \tilde c_k}$ is injective, so condition (2) is automatically satisfied.\\

    \item If $\rho(c_k)$ is parabolic, then $\dev(\tilde c_k)$ is a horocycle centered at the fixed point at infinity of $\rho(c_k)$. In that case, the curve $c_k$ is convex (resp. concave) when $\rho(c_k)$ is a negative (resp. positive) parabolic. Again, in that case, $\dev_{\vert \tilde c_k}$ is injective and condition (2) is automatically satisfied.\\

    \item If $\rho(c_k)$ is elliptic of rotation angle $\alpha \in (0,2\pi)$, or the identity, $\dev(\tilde {c_k})$ is a circle centered at a fixed point of $\rho(c_k)$. Condition (2) requires that, on a fundamental domain $I$ of $\tilde c_k$, $\dev$ does not wind around this circle more than once. Then, if $c_k$ is concave, $\dev(I)$ runs counterclockwise through a portion $\frac{\alpha}{2\pi}$ of the circle, and if $c_k$ is convex, $\dev(I)$ runs through a portion $\frac{2\pi-\alpha}{2\pi}$ of this circle clockwise.
\end{itemize}

As a consequence of the definition, the curvatures of the boundary components of a hyperbolic block are determined up to sign by their length and holonomy. To be more precise:\footnote{These identities follow from the classical relations between length, curvature and radius of circles and hypercircles in the hyperbolic plane.}
\begin{itemize}
    \item If $\rho(c_k)$ is elliptic with rotation angle $\alpha$, then 
    \[\kappa(c_k) =  \sqrt{\left(\frac{\alpha}{\ell(c_k)}\right)^2 +1} \textrm{ or } \kappa(c_k) = -\sqrt{\left(\frac{2\pi-\alpha}{\ell(c_k)}\right)^2+1}~.\]
    \item If $\rho(c_k)$ is hyperbolic with translation length $a$, then 
    \[\kappa(c_k) = \pm \sqrt{1- \left(\frac{a}{\ell(c_k)}\right)^2}~.\]
    \item If $\rho(c_k)$ is a positive parabolic, then \[\kappa(c_k)= 1~.\]
    \item If $\rho(c_k)$ is a negative parabolic, then
    \[\kappa(c_k)= -1~.\]
    \item If $\rho(c_k) = \Id$, then 
    \[\kappa(c_k) = \sqrt{\left(\frac{2\pi}{\ell(c_k)}\right)^2 +1}~.\]
\end{itemize}

Hyperbolic blocks will serve as building blocks that we will glue together in order to get branched hyperbolic structures on closed surfaces.\\

Let us start by noting that singular hyperbolic metrics on a closed surface lead to hyperbolic blocks on the complement of small discs around singular points.

\begin{prop} \label{prop: Cone metrics -> hyperbolic blocks}
    Let $S$ be a compact connected surface with boundary curves $(c_k)_{1 \leq k \leq K}$ and $\rho$ a representation of $\pi_1(S)$ such that $\rho(c_k)$ is elliptic, parabolic or the identity for all $k$. Then the following are equivalent:
    \begin{itemize}
        \item[(i)] $\rho$ is the holonomy of a singular hyperbolic structure on $S^\bullet$ with cusps and cone singularities at the $c_k^\bullet$ (and possibly branch points elsewhere),

        \item[(ii)] $\rho$ is the holonomy of a singular hyperbolic structure on $S^\bullet$ with cusps and cone singularities of angle in $(0,2\pi]$ at the $c_k^\bullet$ (and possibly branch points elsewhere),
        
        \item[(iii)] $\rho$ is the holonomy of a hyperbolic block on $S$ with concave boundary, 

        \item[(iv)] There exists $\epsilon_0>0$ such that, for $(\epsilon_k)_{1\leq k\leq K}\in(0,\epsilon_0)^K$, $\rho$ is the holonomy of a hyperbolic block on $S$ such that $c_k$ is concave of length $\epsilon_k$.
    \end{itemize}
\end{prop}

\begin{proof}
The implications (ii) $\Rightarrow$ (i) and
(iv) $\Rightarrow$ (iii) are clear.

To prove (ii) $\Rightarrow$ (iv), remove a small open metric disc
around each marked point of positive cone angle, and truncate
each cusp along a horocycle, choosing the resulting boundary
lengths to be $\epsilon_k$. These boundaries are concave and
satisfy the defining conditions of a hyperbolic block.

To prove (iii) $\Rightarrow$ (ii), start with a hyperbolic block on $S$ with concave boundary. If $\rho(c_k)$ is elliptic, let $\alpha_k\in(0,2\pi)$ be its rotation angle; if $\rho(c_k)=\Id$, set $\alpha_k=2\pi$. Let $D_k$ be the disc centered at the origin in $C_{\alpha_k}$ whose
boundary has length $\ell(c_k)$. The disc $D_k$ has convex boundary with holonomy a rotation of angle $2\pi - \alpha_k$. By the relations between length, holonomy and curvature given above, the curvature of $\partial D_k$ is opposite to the curvature of $c_k$. We can thus glue $D_k$ to $S$ along $c_k$ as in Proposition \ref{p: paste}.

Similarly, if $\rho(c_k)$ is parabolic, $c_k$ has curvature $-1$ by concavity and we can glue to it a disc centered at the origin in the standard cusp $C_0$ with the same boundary length.

This operation ``fills'' each boundary component with a singular disc, and thus defines a singular hyperbolic structure on $S^\bullet$ with the same holonomy, and with a cusp or a cone of angle in $(0,2\pi]$
at each marked point.

Finally, to prove (i) $\Rightarrow$ (ii), we describe a local
surgery that replaces a cone of angle $\alpha+2\pi$, with
$\alpha>0$, by a branch point of angle $4\pi$ and a cone
of angle $\alpha$.

Let $p$ be a cone point of angle $\alpha+2\pi$.
Choose two embedded geodesic segments $I_1,I_2$ of the same sufficiently small length, meeting only at $p$, and dividing its angle into sectors of angles $2\pi$ and $\alpha$. Cut along these segments using Proposition~\ref{p: cut}.
The two copies $p_+,p_-$ of $p$ have angles $2\pi$ and $\alpha$, while the other endpoints $q_1,q_2$ each have angle $2\pi$.

Glue together the two sides bounding the $2\pi$ sector with each other, and the two remaining sides of the slit with each other, matching
points at equal distance from $p$.
By Proposition~\ref{p: paste}, the resulting surface carries
a singular hyperbolic metric.
The point $p_+$ becomes regular, $p_-$ has angle $\alpha$,
and $q_1,q_2$ merge into a branch point of angle $4\pi$.
The surface is homeomorphic to the original one, and its
hyperbolic structure is unchanged outside the surgery disc.
The identifications agree with the local developing map,
since the two segments are separated by an angle $2\pi$;
hence the holonomy is unchanged. See Figure~\ref{fig:cone-surgery}.
\begin{figure}[htbp]
\centering
\includegraphics[width=\textwidth]{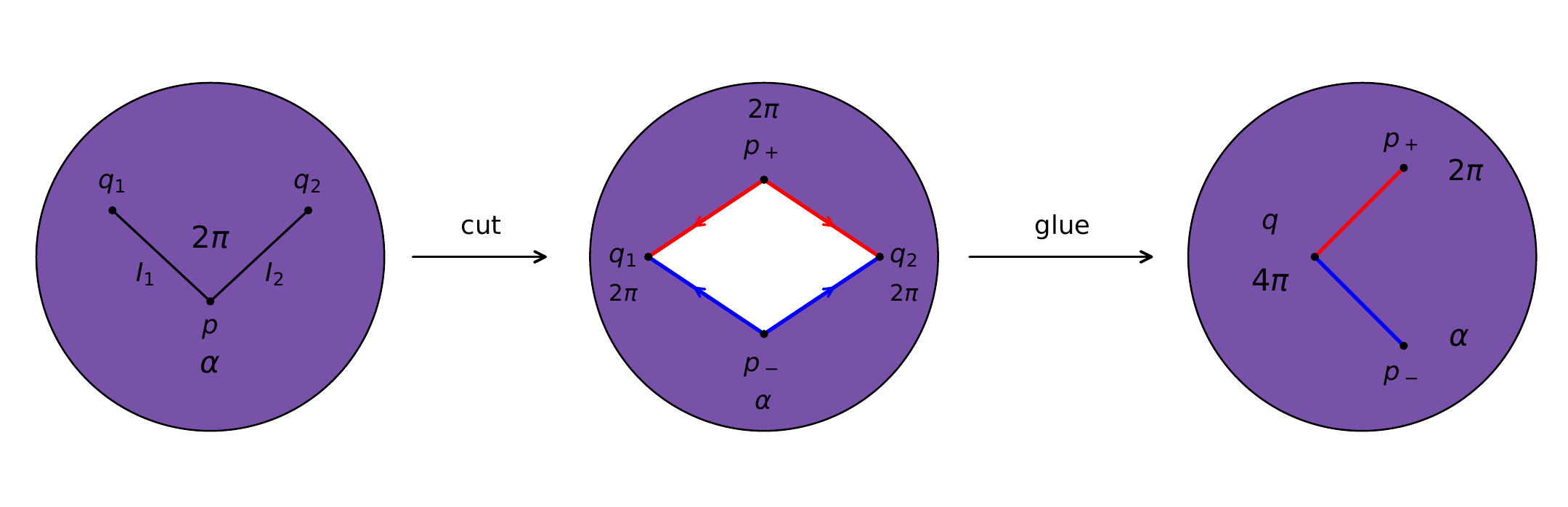}
\caption{Local surgery at a cone point of angle $\alpha+2\pi$.
Cut along the black segments, then glue the two red sides and the two blue sides.
The points $q_1$ and $q_2$ merge into a point of angle $4\pi$, while
$p_+$ becomes regular and $p_-$ has angle $\alpha$.
This is a topological diagram; the labels give the intrinsic angles.}
\label{fig:cone-surgery}
\end{figure}

Applying this construction iteratively at each cone singularity,
one can transform a singular hyperbolic metric with cone
singularities of angles $\alpha_k$ (possibly larger than $2\pi$)
into a singular hyperbolic metric with cone singularities of
angles $\bar\alpha_k\in(0,2\pi]$ and additional branch points
without changing the holonomy.
\end{proof}

We now refine Lemma \ref{l:ehresman-thurston} to get an Ehresmann--Thurston principle for hyperbolic blocks:

\begin{lem} \label{l:perturbation of hyperbolic block}
Let $S$ be a compact surface with boundary components $c_k, 1\leq k \leq K$. Let $\sigma_0$ be a hyperbolic block on $S$ with holonomy $\rho_0$. Then, for every $\epsilon>0$, there exists a neighbourhood $U$ of $\rho_0$ in $\Hom(\pi_1(S), \PSL(2,\R))$ such that every $\rho\in U$ is the holonomy of a hyperbolic block $\sigma$ on $S$ such that, for every boundary curve $c_k$:
\begin{itemize}
    \item if $\rho_0(c_k) = \Id$, then $\ell_\sigma(c_k) < \epsilon$,
    \item otherwise, $\vert \ell_\sigma(c_k)-\ell_{\sigma_0}(c_k)\vert < \epsilon$.
\end{itemize}
\end{lem}

\begin{rmk}
    When $\rho_0(c_k) \neq \Id$ for all $k$, the lemma is a rather straightforward application of Lemma~\ref{l:ehresman-thurston}, and the hyperbolic block $\sigma$ is a small perturbation of $\sigma_0$. The case where $\rho_0(c_k) = \Id$ for some $k$ requires more care, and the hyperbolic block $\sigma$ constructed in Lemma \ref{l:perturbation of hyperbolic block} cannot be chosen arbitrarily close to $\sigma_0$, even when $\rho$ converges to $\rho_0$.
\end{rmk}

\begin{proof}
We will construct the hyperbolic block $\sigma$ near the boundary of $S$, then extend it to the whole surface using Lemma \ref{l:ehresman-thurston}. In order to do so, let us fix a component $c_k$ of $\partial S$ and a lift $\tilde{c}_k$ of $c_k$ to $\tilde S$. Denote by $\gamma_k$ the generator of the stabilizer of $\tilde{c}_k$ in $\pi_1(S)$ which acts on $\tilde{c}_k$ by a positive translation, with respect to the orientation of $c_k$.

Assume first that $\rho_0(\gamma_k) \neq \Id$. There is a neighbourhood $U$ of $\rho_0$ such that $\rho(\gamma_k)\neq \Id$ for every $\rho \in U$. We can then write $\rho_0 (\gamma_k)  = \exp (a_0)$ and $\rho (\gamma_k)= \exp (a)$ with $a$ close to $a_0$ in the Lie algebra $\mathfrak{isom}(\H)$ of $\Isom_+(\H)$. The orbits of the subgroup $\exp (\mathbb R a_0)$ are curves of constant curvature in $\H$, or fixed points of $\rho_0(\gamma_k)$. 

By definition of a hyperbolic block, the developing map $\dev_0$ sends the curve  $\tilde{c}_k$ to one of these curves (but not to a fixed point). More precisely, we can take a parametrization of $\tilde{c}_k$ by $\mathbb R$ so that the restriction of $\dev_0$ to $\tilde{c}_k$ is given by $t\mapsto \exp (t a_0) q$, for some point $q\in \H$ which is not fixed by $\rho_0(\gamma_k)$. We now define $d_k(t) = \exp (t a) q$ for $t\in \mathbb R \simeq \tilde{c}_k$. This map is close to ${\dev_0}_{\vert \tilde c_k}$. In particular, for $\rho$ close enough to $\rho_0$, it is still injective on a fundamental domain of $\tilde{c}_k$. By Lemma \ref{l:ehresman-thurston}, we can thus find a branched hyperbolic structure on $S$ with holonomy $\rho$ and developing map $\dev$ such that $\dev|_{\tilde c_k} = d_k$ for each $k$. This gives the hyperbolic block with holonomy $\rho$.\\

The situation where the holonomy of some peripheral curve is the identity is more delicate, and the structure $\sigma$ cannot be close to $\sigma_0$. The construction needs two steps in that case: we first perturb the hyperbolic block structure $\sigma_0$ into a singular hyperbolic structure $\sigma'$, and then glue to $\sigma'$ some piece of hyperbolic surface to get the desired $\sigma$.

We begin by constructing the structure $\sigma'$ on the boundary of $S$. We proceed exactly as before for the peripheral curves that have a $\rho_0$-holonomy different from the identity, namely we define $d_k(t)= \exp(ta)q$ on $\tilde{c}_k$, for $a\in \mathfrak{isom}(\H)$ such that $\exp(a) = \rho(\gamma_k)$. 

For a boundary component $c_k$ with $\rho_0(\gamma_k)=\Id$,
the map $\dev_0$ descends to an orientation-preserving diffeomorphism
from $\langle\gamma_k\rangle\backslash\tilde c_k$ onto a round circle
oriented counterclockwise. If $\rho(\gamma_k)=\Id$ as well, set
$d_k=\dev_0|_{\tilde c_k}$; we will shrink this boundary after extending
the boundary data. In the remaining cases, write
$\rho(\gamma_k)=\exp(a)\ne\Id$, with $a$ small, and proceed as follows.

First assume that $\rho(\gamma_k)$ is elliptic with centre close to that of the circle $\dev_0(\tilde{c}_k)$. We can then find a small perturbation $d_k$ of ${\dev_0}_{\vert \tilde c_k}$ such that $\tilde{c}_k$ is mapped to a circle invariant under $\rho(\gamma_k)$ and oriented positively. 

Otherwise, the element \(a\in \mathfrak{isom}(\H)\), viewed as a vector field on \(\mathbb \H\), is transverse to the circle \( \dev_0 (\tilde{c}_k)\) at some point. In particular, there is an interval in \( \dev_0( \tilde{c}_k)\) whose image under \(\exp(ta) \) is disjoint from itself for all sufficiently small nonzero \(t\). Choosing a point \(q\) in this interval,  we parametrize $\tilde{c}_k$ by $\mathbb R$ equivariantly, in such a way that $\dev_0 (0)= q$, and we define as before  $d_k(t) = \exp (ta ) \dev_0(t)$. The image of $d_k$ is a helicoidal curve in $\mathbb H$, invariant under the transformation $\rho(\gamma_k)$ (this is the green curve in Figure \ref{fig: Escargot}).

\begin{figure}[h]
\centerline{\includegraphics[scale=0.6]{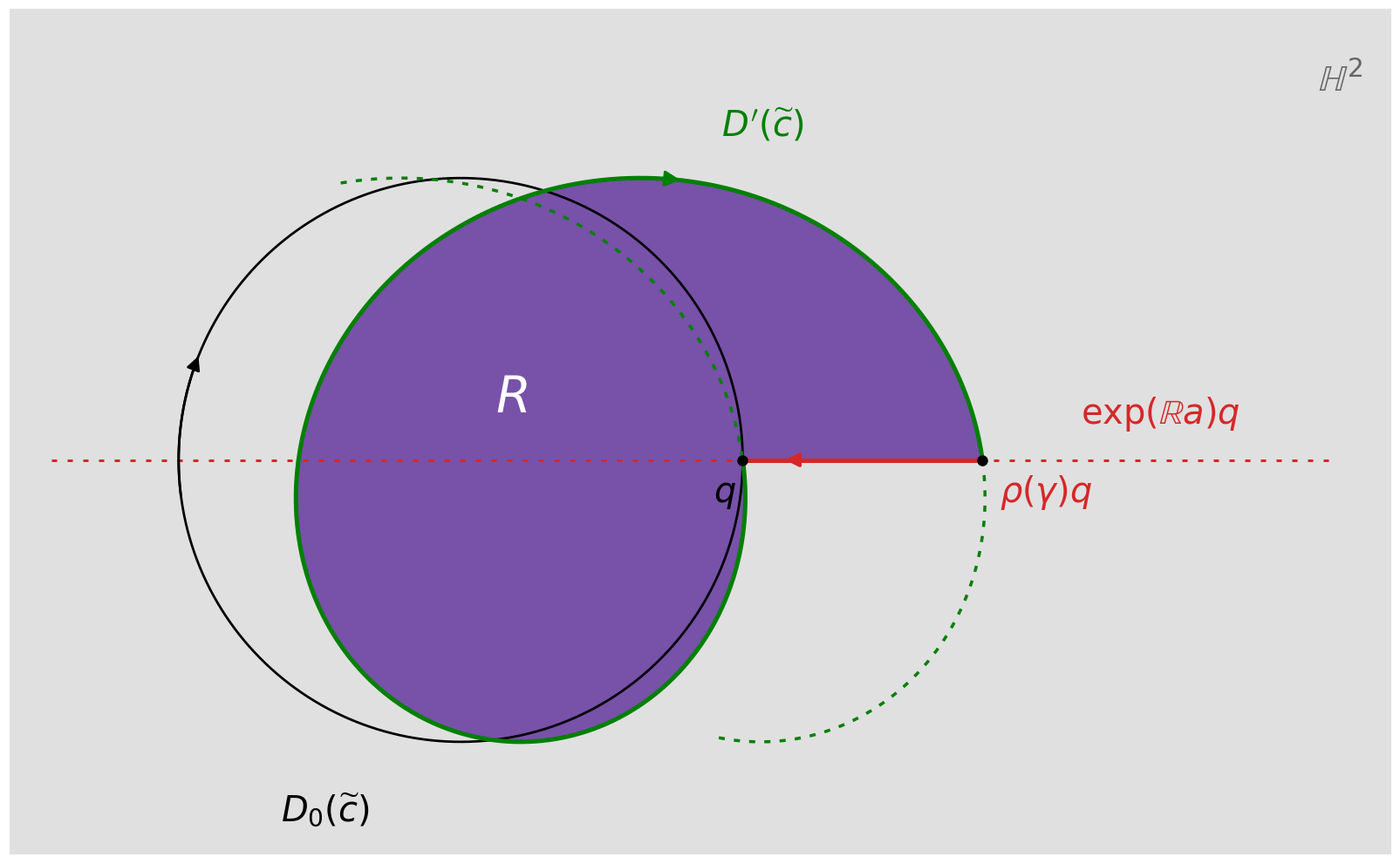}}
\caption{The curve $D'$.}
 \label{fig: Escargot}
\end{figure}

By Lemma~\ref{l:ehresman-thurston}, these boundary data extend to
a branched hyperbolic structure $\sigma'$ on $S$, with holonomy $\rho$
and developing map $\dev'$ satisfying $\dev'|_{\tilde c_k}=d_k$.
For every component with $\rho_0(\gamma_k)=\rho(\gamma_k)=\Id$,
fill its concave circular boundary with a regular hyperbolic disc,
then remove a smaller concentric open disc with boundary length
less than $\epsilon$. This preserves the holonomy and gives a regular
concave boundary satisfying the hyperbolic block conditions.

For the remaining components with $\rho_0(\gamma_k)=\Id$,
the structure $\sigma'$ may still fail to be a hyperbolic block:
a fundamental interval may develop more than once around a circle, or to a helicoidal curve with non-constant curvature.

In the circular case, $\dev'$ maps a fundamental domain of $\tilde c_k$ to a portion $\frac{\alpha_k}{2\pi}$ of a circle, where $\alpha_k$ is close to $2\pi$  and the  holonomy of $c_k$ is a rotation of angle $\alpha_k$. If $\alpha_k<2\pi$ then $\dev'$ is injective on a fundamental domain of $\tilde c_k$. Glue to $c_k$ a disc centreed at $0$ in the cone $C_{\alpha_k}$  to obtain a cone singularity of angle $\alpha_k$. If $\alpha_k > 2\pi$, replace this singularity by a singularity of angle $\alpha_k-2\pi$ and a branch point (as in the proof of Proposition~\ref{prop: Cone metrics -> hyperbolic blocks}, (i) $\Rightarrow$ (ii)). Finally remove a small disc around the remaining cone singularity (as in the proof of Proposition~\ref{prop: Cone metrics -> hyperbolic blocks}, (ii) $\Rightarrow$ (iv)) to obtain a hyperbolic structure which satisfies the conditions of a hyperbolic block near $c_k$ and for which $c_k$ is arbitrarily small.

Finally, suppose that $\dev'$ maps $\tilde c_k$ to a helicoidal curve.
Set $g=\rho(\gamma_k) =\exp(a)$, $d(t)=\dev'(t)$ and
$\eta(t)=\exp(ta)q$ for $t\in[0,1]$.
We glue to $\sigma'$ the domain $R$ bounded by the simple closed curve
obtained by following $d$ from $q$ to $gq$ and then $\eta$ in the reverse
direction, attaching $R$ along $d$.

This creates a corner of angle $3\pi$ at the boundary point $p$ obtained
by identifying the endpoints of $\eta$.
Indeed, equivariance identifies the tangent vectors of $d$ at its two
endpoints by $Dg$, and the same holds for $\eta$.
Consequently, the two interior angles of $R$ at these endpoints add up
to $2\pi$. Together with the angle $\pi$ contributed by the original
surface, they give an angle $3\pi$ at $p$.

We remove this corner by a boundary version of the surgery used in the
proof of Proposition~\ref{prop: Cone metrics -> hyperbolic blocks}.
Choose two sufficiently short embedded geodesic segments $J_1,J_2$ of
the same length, starting at $p$, disjoint except at $p$, and otherwise
contained in the interior of the surface. Choose their directions so
that the three successive sectors at $p$ have angles
\[
    \theta,\qquad 2\pi,\qquad \pi-\theta,
    \qquad 0<\theta<\pi.
\]
The two segments have the same developed image.
Cut along $J_1\cup J_2$. Identify the two sides adjacent to the middle
sector with each other, and identify the two remaining sides with each
other, in both cases matching points at the same distance from $p$.
These identifications agree with the developing map.

The copy of $p$ bounding the middle sector becomes a regular interior
point of angle $2\pi$. The two other copies of $p$ are identified to
form a regular boundary point of angle
$\theta+(\pi-\theta)=\pi$.
Finally, the other endpoints of $J_1$ and $J_2$, each of angle $2\pi$
after cutting, are identified to form an interior branch point of
angle $4\pi$.
The resulting neighbourhood is again a disc, and the hyperbolic
structure outside this neighbourhood is unchanged. Thus the surgery
preserves the topology and the holonomy.

The boundary arcs and their lengths are unchanged by this surgery.
Their developed images join smoothly along the orbit
$\exp(\mathbb R a)q$; since the boundary angle is now $\pi$, the new
boundary is regular and has constant curvature.
Its developing map is injective on a half-open fundamental domain,
and its length is the length of the arc $\eta$.
This length tends to zero as $a$ tends to zero, uniformly for $q$ on
the fixed circle $\dev_0(\tilde c_k)$.
Performing this construction at each boundary component of this type
gives the required hyperbolic block. The lemma follows.
\end{proof}

\subsection{A gluing lemma}

In this subsection, we explain how to glue together blocks with compatible boundaries to create larger blocks with prescribed holonomy. Though the construction is rather intuitive, it requires some care.

Let $S$ be a compact connected oriented surface with boundary components $c_k, 1\leq k \leq K$, and let 
\[S = \bigsqcup_{1\leq i \leq I} S_i \sqcup \bigsqcup_{1\leq l \leq L} b_l\]
be a decomposition of $S$ (see Section \ref{ss:Decompositions surface}). Recall that $S^\sharp$ denotes the metric completion of $\bigcup_{1\leq i \leq I} S_i$, and $b_l^+$ and $b_l^-$ denote the two copies of $b_l$ in $S^\sharp$.

Let $\rho: \pi_1(S) \to \Isom_+(\H)$ be a representation. The goal of this section is to prove the following gluing lemma for hyperbolic blocks:

\begin{lem}\label{l: gluing lemma}
Suppose that $\rho(b_l)\neq\Id$ for every cutting curve $b_l$.
Assume that $S^\sharp$ carries a hyperbolic block $\sigma^\sharp$
with holonomy $\rho|_{\pi_1(S_i)}$ on each component, and that for
each $l$ one of the following holds:
\begin{itemize}
\item[(a)] both $b_l^+$ and $b_l^-$ are concave;
\item[(b)] both $b_l^+$ and $b_l^-$ are geodesic;
\item[(c)] after exchanging the two sides if necessary, $b_l^+$ is
concave, $b_l^-$ is convex, and
$\ell_{\sigma^\sharp}(b_l^-)\leq\ell_{\sigma^\sharp}(b_l^+)$.
\end{itemize}
Then $S$ carries a hyperbolic block structure $\sigma$ with holonomy $\rho$
such that $\kappa_\sigma(c_k)=\kappa_{\sigma^\sharp}(c_k)$ for every
boundary component $c_k$ of $S$.
\end{lem}

The  main idea is that we can glue a hyperbolic cylinder with the appropriate boundaries to $b_l^-$ and $b_l^+$. Because of the constant curvature assumption, this gluing operation has a ``twist'' parameter that allows to recover all possible holonomies with a given restriction. Let us start by building these hyperbolic cylinders. Let $C$ be an
oriented cylinder with boundary components $c_-$ and $c_+$.

\begin{figure}[h] \label{fig: plug}
\centerline{\includegraphics[scale=0.4]{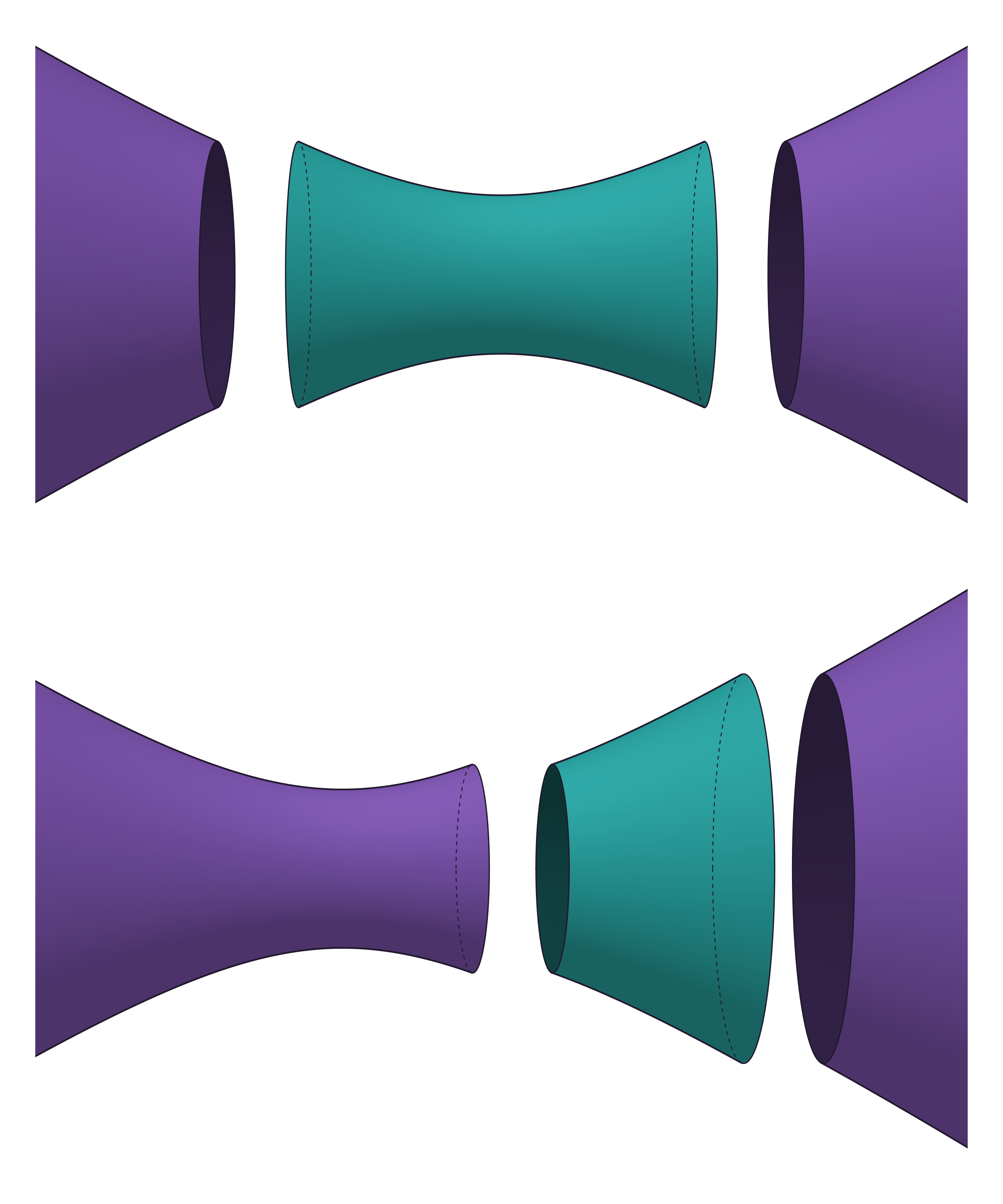}}
\caption{The plugs.}
\end{figure}

\begin{prop}\label{p: Cylinders}
\begin{itemize}
\item[(1)] Let $g$ be non-trivial elliptic, hyperbolic, or positive
parabolic. For every $\ell(g)<l_1<l_2$, there is a hyperbolic block
on $C$ whose boundary $c_-$ is concave, of length $l_1$ and holonomy
$g$, while $c_+$ is convex, of length $l_2$ and holonomy $g^{-1}$.
\item[(2)] Let $g$ be non-trivial elliptic or hyperbolic.
For every $l_1,l_2>\ell(g)$, there is a hyperbolic block on $C$ whose
two boundaries are convex, with respective lengths $l_1,l_2$ and
holonomies $g,g^{-1}$.
\end{itemize}
\end{prop}

\begin{proof}
In case (1), if $g$ is a rotation of angle $\alpha\in(0,2\pi)$,
take the annulus $r_1\leq r\leq r_2$ in the cone $C_\alpha$, where
$r_j=\sinh^{-1}(l_j/\alpha)$ and $c_-$ is the inner boundary.
Indeed, the circle of radius $r$ has length $\alpha\sinh r$.
If $g$ is a positive parabolic, conjugate it to $z\mapsto z+1$ and use
\[
    \{z\in\H:1/l_2\leq\operatorname{Im}z\leq1/l_1\}
    /\langle z\mapsto z+1\rangle,
\]
with $c_-$ the upper boundary.

If $g$ is hyperbolic, put $L=\ell(g)$. Let \(r\) be the signed distance to the axis of \(g\), positive on its right, and let \(s\) be the arclength coordinate of the orthogonal projection onto the axis, increasing in the direction of \(g\). Then $(r,s)$ give global coordinates on $\H$ with respect to which the hyperbolic metric and action of $g$ are given by
\[
    \sigma_\H = dr^2+\cosh^2 r\,ds^2,
    \qquad g(r,s)=(r,s+L).
\]
Set $r_j=\cosh^{-1}(l_j/L)>0$.
For (1), quotient the strip $r_1\leq r\leq r_2$ by $g$;
for (2), use $-r_1\leq r\leq r_2$.
The boundary $c_-$ is the image of $r=r_1$ in case (1),
and of $r=-r_1$ in case (2). In both cases, its orientation
corresponds to increasing $s$, and its holonomy is $g$.
The circle at signed distance $r$ has length $L\cosh r$, and the
required convexity properties follow from the position of the axis.

For (2) with $g$ a rotation of angle $\alpha$, take discs
$D_-\subset C_{2\pi-\alpha}$ and $D_+\subset C_\alpha$ of boundary
lengths $l_1$ and $l_2$. Cut each disc along a radial segment of the same
sufficiently small length, with its other endpoint in the interior.
Glue the two slit boundaries by an orientation-reversing isometry,
matching the cone vertices and matching the other endpoints.
The result is a cylinder: the cone vertices give a regular point of
angle $(2\pi-\alpha)+\alpha=2\pi$, while the other endpoints give
one branch point of angle $4\pi$.
The original disc boundaries remain convex and retain their lengths.
Their oriented holonomies are rotations of angles
$-(2\pi-\alpha)\equiv\alpha$ and $-\alpha$, respectively, as required.
\end{proof}

\begin{proof}[Proof of Lemma~\ref{l: gluing lemma}]
By an immediate induction, it suffices to treat one cutting curve at a time. Let $b$ be one curve of the decomposition and $b_+,b_-$ its two copies in $S^\sharp$. We start by gluing a cylinder to $b_+$ to obtain a new hyperbolic block with holonomy $\rho\vert_{S^\sharp}$ and for which $b_+$ and $b_-$ have opposite curvature.

If both $b_+, b_-$ are geodesic, their lengths equal the translation length of $\rho(b)$, and no cylinder is needed.
Otherwise label the sides so that $b_+$ is concave, and write
$L_\pm=\ell(b_\pm)$.
If $b_-$ is also concave, then $\rho(b_-)$ and $\rho(b_+)$ cannot be parabolics, since they should both be positive parabolics by concavity, while on the other hand $\rho(b_-)$ is conjugate to $\rho(b_+)^{-1}$. 
By Proposition~\ref{p: Cylinders}, part (2), we can find a hyperbolic cylindrical block $C$ whose boundary curves $b'_-$, $b'_+$ have respective lengths $L_+$, $L_-$ and holonomies $\rho(b_+)^{-1}$, $\rho(b_+)$. 
By the relations between holonomy, length and curvature recalled in Section \ref{ss: Hyperbolic blocks}, we can glue $b'_-$ to $b_+$ and obtain a new hyperbolic block structure on $S^\sharp$, with the same holonomy, such that the new boundary curves $b'_+$ and $b_-$ now have the same length and opposite curvatures. Similarly, if $b_-$ is convex and $L_-<L_+$, using Proposition~\ref{p: Cylinders}, part (1), we can again find a cylinder $C$ as above and glue it to $b_+$ to get a new hyperbolic block where $b_-$ and $b_+$ now have the same length and opposite curvature.

From now on, we thus assume that $b_-$ and $b_+$ have the same length and opposite curvatures. Assume first that $b$ separates $S$ into two connected components $S\pm$. Choose lifts $\widetilde S_\pm\hookrightarrow\widetilde S$ of the inclusions $S_\pm\hookrightarrow S$, and a lift $\tilde b\hookrightarrow\widetilde S$ of $b\hookrightarrow S$ as in Section \ref{sss:Van Kampen Separating} and identify correspondingly $\pi_1(S)$ with the amalgamated product
\[\pi_1(S_+) *_{\langle b\rangle} \pi_1(S_-)~.\]
Choose developing maps $\dev_\pm$ on the two components $\tilde S_\pm$ with holonomies equal to the restrictions of $\rho_\pm \coloneqq \rho \vert_{S_\pm}$.
Since $\rho_+(b_+)= \rho_-(b_-)^{-1} = \rho(b)$, the curves $\dev_-(\tilde b_-)$ and $\dev_+(\tilde b_+)$ are $\rho(b)$-invariant. Since they have opposite curvature, and since $\rho(b) \neq \Id$, there exists an orientation-reversing boundary isometry
$f:b_-\to b_+$ with a lift $\tilde f: \tilde b_- \to \tilde b_+$ satisfying
\[
    \dev_+\circ\widetilde f=\dev_-|_{\widetilde b_-}.
\]
Gluing $b_-$ to $b_+$ by $f$ yields a hyperbolic block structure on $S$ with a developing pair $(\dev,\rho')$ satisfying $\dev\vert_{\tilde S_\pm} = \dev_\pm$ and $\rho'\vert_{S\pm}= \rho_\pm$ by Lemma \ref{lem: holonomy gluing separating case}. Hence $\rho' = \rho$ and we have constructed the desired hyperbolic block on $S$.

Suppose now that $b$ does not separate $S$. Let $S'$ be the surface $S$ cut along $b$, and choose lifts $\tilde b_\pm \hookrightarrow\tilde{S'}\hookrightarrow \tilde S$, as in Section \ref{sss: Van Kampen non separating}, and identify correspondingly $\pi_1(S)$ with the HNN extension
\[\pi_1(S')*\langle \gamma \rangle /(\gamma \beta_- \gamma^{-1} = \beta_+^{-1})~,\]
where $\beta_\pm$ denote the positively oriented generators of $\pi_1(b_\pm) \subset \pi_1(S)$.
Let $\dev'$ be a developing map of the given hyperbolic block on $S'$ with holonomy $\rho'\coloneqq \rho|_{\pi_1(S')}$. The curves $\rho(\gamma)\dev'(\tilde b_-)$ and $\dev'(\tilde b_+)$ have opposite curvature and are both invariant under $\rho(\beta_+)^{-1}= \rho(\gamma) \rho(\beta_-)\rho(\gamma)^{-1}$.  Hence they coincide modulo orientation and we obtain an orientation-reversing isometry $f:b_-\to b_+$ and a lift $\tilde f: \tilde b_- \to \tilde b_+$ satisfying
\[
    \dev'\circ\widetilde f
    =\rho(\gamma)\circ\dev'|_{\widetilde b_+}~.
\]
Again, gluing $b_-$ to $b_+$ by $f$ yields a hyperbolic block structure on $S$ whose holonomy representation coincides with $\rho$ on $\pi_1(S')$ and on $\gamma$ by Lemma \ref{lem: holonomy gluing non-separating}, hence is equal to $\rho$.
All constructions leave the original boundary components of $S$ unchanged, which proves the conclusion.
\end{proof}

\section{Geometrizing small representations of the fundamental group of a pair of pants}

We recall that we adopt the unusual convention that each boundary component is oriented so that the surface lies to the right of its boundary.

Let $\Sigma_{0,3}$ denote a pair of pants, i.e. the surface of genus $0$ with three boundary components. We present its fundamental group as
\[ \pi_1 (\Sigma_{0,3} ) =\ \langle u,v,w\mid uvw=1\rangle, \]
where the elements \(u,v,w\) are represented by loops homotopic to the three boundary components of the surface.

Our goal is to establish the following result:

\begin{lem} \label{l: geometrization of representations with small translation length}
Let  $\rho_n : \pi_1(\Sigma_{0,3}) \to \Isom_+(\H)$ be a sequence of non-trivial representations whose minimal displacement \(\delta_{\rho_n}\) with respect to \(\{u,v,w\}\) tends to zero when \(n\) tends to infinity. Then, for every $\epsilon>0$ there exists $N$ such that $\rho_n$ is the holonomy of a hyperbolic block with boundary lengths $\leq \epsilon$ for all $n\geq N$.
\end{lem}

\subsection{Conjugating and taking limits}

We start by noticing that, up to conjugating and extracting a subsequence, the sequence $\rho_n$ converges to a non-trivial abelian representation.

\begin{prop}\label{p: bounded displacement}
Let \(G\) be a finitely generated group with finite generating set \(S\).
Let \(\rho_n : G \to \Isom_+(\H)\) be a sequence of non-trivial representations such that the minimal displacement sequence \((\delta_{\rho_n,S})_n\) is bounded.
Then, after conjugating each \(\rho_n\) by some \(g_n \in \Isom_+(\H)\) and passing to a subsequence, the sequence \((\rho_n)\) converges to a non-trivial representation \(\rho\).
\end{prop}

\begin{proof}
Let $C>0$ be such that \( \delta_{\rho_n} < C \) for all $n$. Since the displacement $\delta_g(x)$ of a non-trivial element \( g\in \Isom_+(\H) \) is an unbounded function of $x$, and since \( \rho_n \) is not the trivial representation, the function \( x\mapsto \delta_{\rho_n,S}(x) \) is unbounded as \(x\) ranges over \( \mathbb H \). By the intermediate value theorem, there thus exists a point \(x_n \in \mathbb H\) such that \( \delta_{\rho_n}(x_n)= C \).

Let \(g_n \in \Isom_+(\H)\) be such that \( g_n(x_n)=o \), where \(o \in \mathbb H\) is a fixed base point, and define \( \rho_n' := g_n \rho_n g_n^{-1} \). These representations satisfy \( \delta_{\rho_n'}(o)= C \). Hence, up to passing to a subsequence, we may assume that \( \rho_n' \) converges to a representation \( \rho' \) as \(n \to \infty\). Moreover, this limit representation is non-trivial, since it satisfies \( \delta_{\rho'}(o)=C >0 \).
\end{proof}

Applying Proposition~\ref{p: bounded displacement} to the sequence \((\rho_n)_n\) of representations given by the statement of Lemma~\ref{l: geometrization of representations with small translation length} (or to any of its subsequences), we find that, up to passing to a subsequence and conjugating by suitable elements of \(\Isom_+(\H)\), one may assume that \(\rho_n\) converges to a representation \(\rho\). Since minimal displacement depends continuously on the representation, the minimal displacement of \(\rho\) vanishes.

Lemma \ref{l: geometrization of representations with small translation length} will thus follow from Lemma~\ref{l:perturbation of hyperbolic block}, and the following lemma, which deals with representations with abelian image:

\begin{lem} \label{l:geometrization pair of pants abelian}
 Let \( \rho : \pi_1 (\Sigma_{0,3} ) \rightarrow \Isom_+(\H)\) be a non-trivial representation such that $\delta_{\rho,S} = 0$. Then, for every \(\varepsilon >0\), there exists a hyperbolic block structure on \(\Sigma_{0,3}\) whose boundary components have length bounded by \(\varepsilon\) and whose holonomy is \(\rho\). 
\end{lem}

\begin{proof}[Proof of Lemma \ref{l: geometrization of representations with small translation length} assuming Lemma \ref{l:geometrization pair of pants abelian}]

Fix some $\epsilon>0$. By Lemma \ref{l:perturbation of hyperbolic block}, the set $\Omega$ of holonomies of hyperbolic blocks on $\Sigma_{0,3}$ with boundary curves of length $<\epsilon$ is open in $\Hom(\pi_1(\Sigma_{0,3}), \Isom_+(\H))$. It is also invariant under conjugation. Finally, by Lemma \ref{l:geometrization pair of pants abelian}, it contains every non-trivial representation with vanishing minimal displacement. Let $(\rho_n)_{n\in \N}$ be as in Lemma \ref{l: geometrization of representations with small translation length}, and suppose for a contradiction that there exists a subsequence $(\rho_{n_k})$ that never lies in $\Omega$. By Proposition \ref{p: bounded displacement}, there exists a further subsequence $(\rho_{n_l})$ and conjugates $\rho'_{n_l}$ of $\rho_{n_l}$ that converge to a non-trivial representation $\rho_\infty$. By continuity of $\rho\mapsto \delta_{\rho,S}$, we have $\delta_{\rho_\infty,S}=0$. By Lemma \ref{l:geometrization pair of pants abelian}, $\rho_\infty$ belongs to $\Omega$. Since $\Omega$ is open we conclude that $\rho_{n_l}$ belongs to $\Omega$ for $n_l$ large enough, which contradicts the assumption on the sequence $(\rho_{n_k})$. We thus conclude that $\rho_n$ belongs to $\Omega$ for $n$ large enough.
\end{proof}

Representations as in Lemma \ref{l:geometrization pair of pants abelian} take values either in a one-parameter elliptic subgroup or in a one-parameter parabolic subgroup. We treat these two cases separately.

\subsection{The elliptic case}

Here we prove Lemma \ref{l:geometrization pair of pants abelian} when $\rho$ takes values in a one-parameter group of rotations. The image of $\rho$ consists of rotations about a common point $o$.\\

Suppose first that one peripheral holonomy is trivial. Since $\rho$
is non-trivial, after relabelling the boundary components we may write
\[
\rho(u)=r_\theta,\qquad
\rho(v)=r_{-\theta},\qquad
\rho(w)=\Id,
\qquad 0<\theta<2\pi,
\]
where $r_\theta$ denotes rotation of angle $\theta$ centreed at $o$.
Take an annulus $r_1\leq r\leq r_2$ in the hyperbolic cone
$C_\theta$, with $0<r_1<r_2$ sufficiently small, and remove a small
open metric disc whose closure lies in its interior.
The resulting pair of pants has the required holonomy: its inner
and outer circular boundaries have holonomies $r_\theta$ and
$r_{-\theta}$, respectively, while the new boundary is concave
and has trivial holonomy. All three boundary lengths can be made
smaller than $\varepsilon$. There are no branch points.

We may now assume that all three peripheral holonomies are non-trivial.
Write their rotation angles as $2\pi-\alpha$, $2\pi-\beta$,
and $2\pi-\gamma$, with $\alpha,\beta,\gamma\in(0,2\pi)$.
The relation $\rho(u)\rho(v)\rho(w)=\Id$ implies that
$\alpha+\beta+\gamma$ equals either $2\pi$ or $4\pi$.

\begin{figure}[h]
    \centering
    \includegraphics[width=0.5\textwidth]{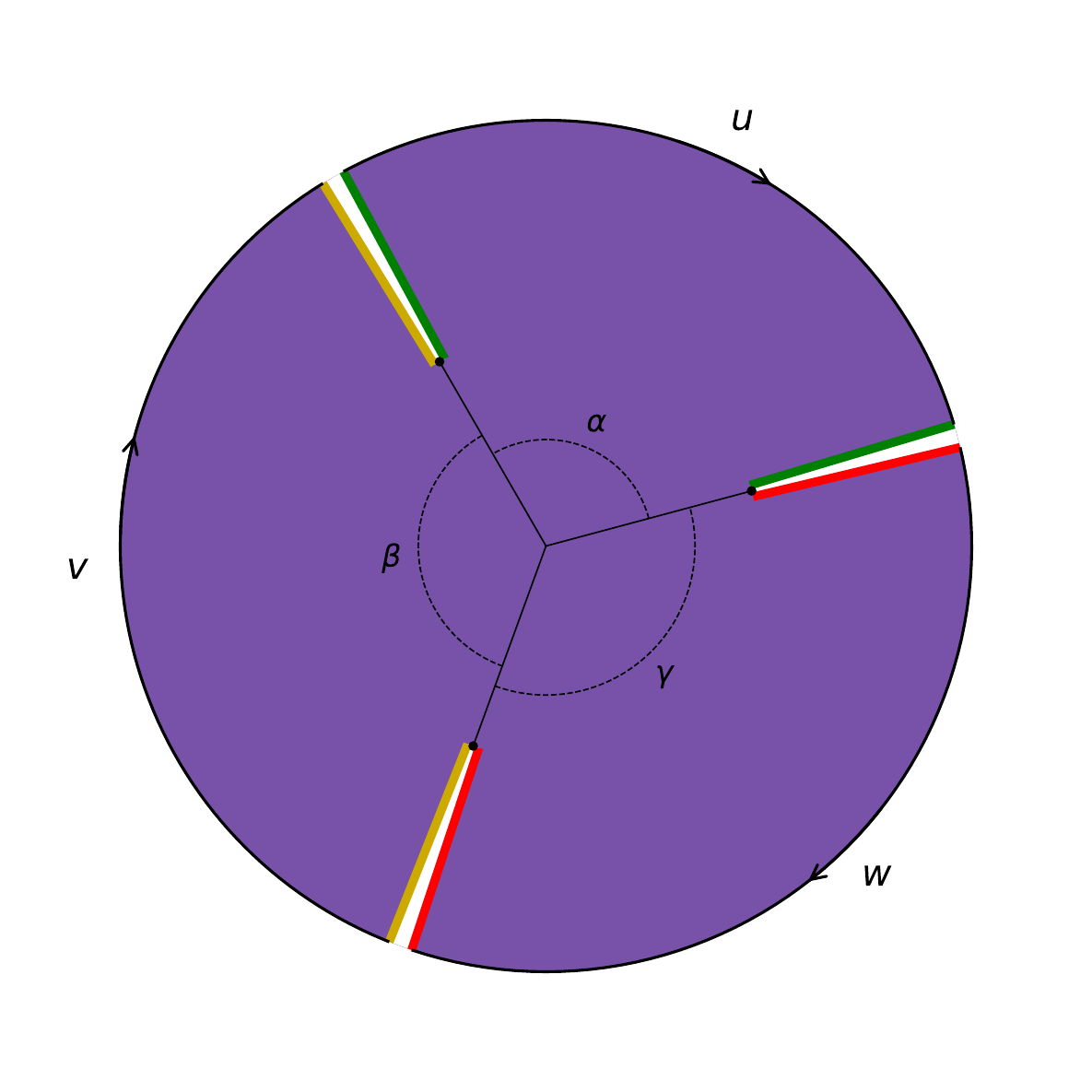}
\captionsetup{justification=centering}
\caption{Completion of a radially slit $\varepsilon$-hyperbolic ball with angles $\alpha$, $\beta$, and $\gamma$.}
\end{figure}

Assume first that $\alpha+\beta+\gamma=2\pi$.
Let $B$ be a hyperbolic disc of radius $r$ centreed at $o$.
Choose three radial segments dividing $B$ into sectors of angles
$\alpha,\beta,\gamma$, and fix $t\in(0,r)$.
Cut each radial segment from its point at distance $t$ from $o$
up to the boundary of $B$.
For each sector, glue the two slits bounding that sector
by the rotation about $o$ carrying one bank to the other.

The three inner endpoints of the slits become a single point
of angle $6\pi$. At the outer endpoints, two angles of $\pi/2$
are joined, so the boundary is smooth.
The quotient is a pair of pants whose three boundary components
are convex and have lengths
$\alpha\sinh r$, $\beta\sinh r$, and $\gamma\sinh r$.
With our boundary orientation, their holonomies are rotations
about $o$ through angles $-\alpha$, $-\beta$, and $-\gamma$,
respectively, hence are precisely $\rho(u)$, $\rho(v)$, and $\rho(w)$.
Taking $r$ sufficiently small proves the claim in this case.
If $\alpha+\beta+\gamma=4\pi$ the same construction works if we replace the hyperbolic disc $B$ with a disc in $C_{4\pi}$ centreed at the origin.

\subsection{The parabolic case}

Here we prove Lemma \ref{l:geometrization pair of pants abelian} when $\rho$ fixes a point in $\partial_\infty \H$.\\

Conjugating the representation by an element of \(\Isom_+(\H)\), one may assume that the representation fixes $\infty$. Since the minimal displacement vanishes, the three maps \( \rho(u), \rho(v), \rho(w)\) are translations of the form \(z\mapsto z+\tau_u\), \(z\mapsto z+\tau_v\), and \(z\mapsto z+\tau_w\), respectively.

Suppose first that one peripheral holonomy is trivial.
Since $\rho$ is non-trivial, after relabelling the boundary components
we may assume that
\[
\tau_u>0,\qquad \tau_v=-\tau_u,\qquad \tau_w=0.
\]
Choose $0<h_0<h_1$ with $\tau_u/h_0<\varepsilon$ and take the cylinder
\[
\{z\in\H:h_0\leq\operatorname{Im}z\leq h_1\}
/\langle z\mapsto z+\tau_u\rangle.
\]
Its upper and lower boundaries have holonomies $\rho(u)$ and
$\rho(v)$, respectively, and lengths smaller than $\varepsilon$.
Remove an embedded open metric disc whose closure lies in the
interior and whose boundary length is smaller than $\varepsilon$.
The resulting surface is a pair of pants with holonomy $\rho$.
The new boundary is concave and has trivial holonomy.
There are no branch points.

Let us now assume that none of \(\tau_u,\tau_v,\tau_w\) vanishes. The relation \( \rho(u)\rho(v)\rho(w)= \Id\) implies that the sum \(\tau_u + \tau_v + \tau_w \) vanishes. So among the numbers \( \tau_u, \tau_v , \tau_w\), two of them have the same sign, and the third the opposite sign.

Assume first that $\tau_u>0$ and $\tau_v,\tau_w<0$.
Since $\tau_u=-(\tau_v+\tau_w)$, choose
$0<h_0<h_1$ with $\tau_u/h_0<\varepsilon$.
Consider the rectangle
\[
R=\{x+iy:0\leq x\leq\tau_u,\ h_0\leq y\leq h_1\}.
\]
Its horizontal sides have hyperbolic lengths $\tau_u/h_0$
and $\tau_u/h_1$, both smaller than $\varepsilon$.

\begin{figure}[h]
    \centering
    \includegraphics[width=1\textwidth]{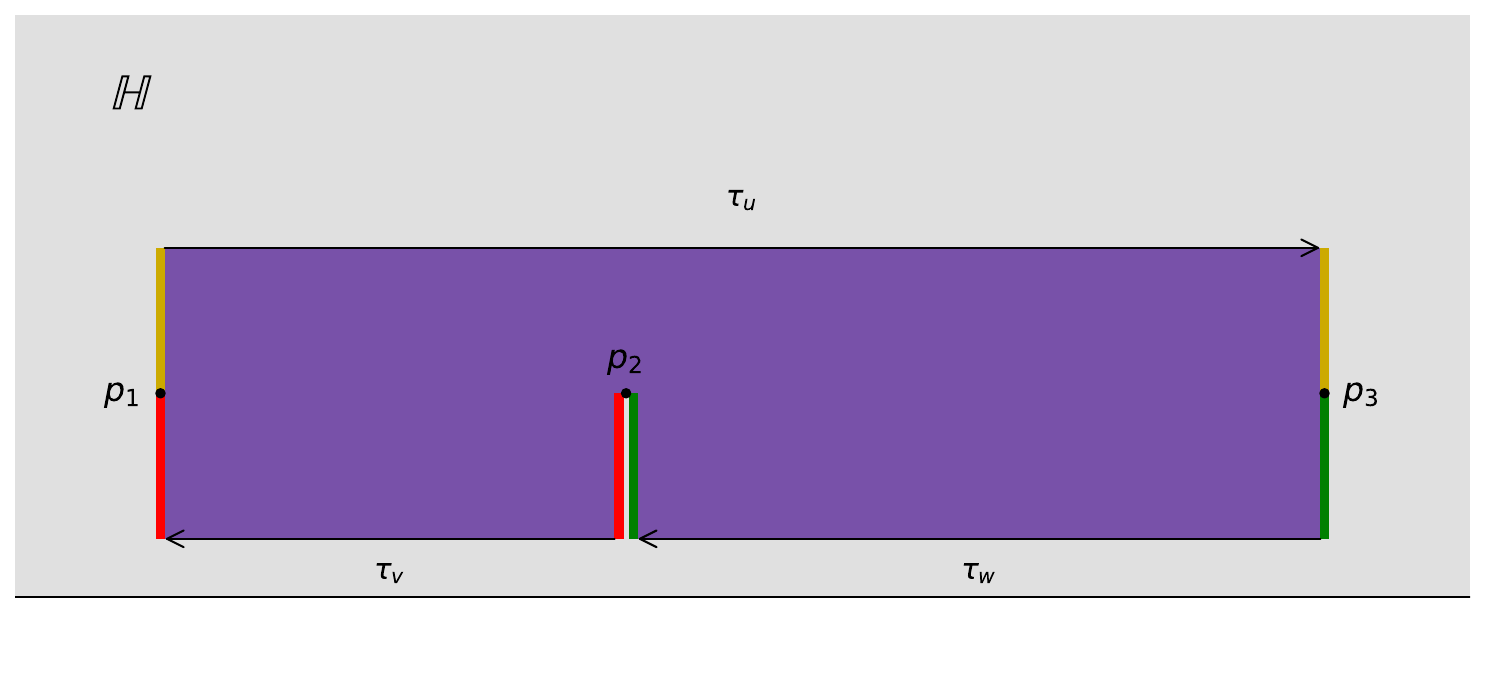}
    \caption{Completion of a rectangle with a vertical slit
    emanating from its lower side, for $\tau_u>0$ and $\tau_v,\tau_w<0$.}
    \label{f: slit rectangle}
\end{figure}

Choose $h_*\in(h_0,h_1)$ and cut $R$ along the vertical segment $(-\tau_v+ih_0, -\tau_v+ih_*)$, as in Figure~\ref{f: slit rectangle}.

Glue the right bank of the slit to the lower part of the right
vertical side by the translation $\rho(w)^{-1}$, and the left
bank to the lower part of the left vertical side by $\rho(v)$.
Glue the remaining upper parts of the vertical sides by $\rho(u)$.
The quotient is a pair of pants. Its upper boundary has holonomy
$\rho(u)$, while the lower left and lower right boundaries have
holonomies $\rho(v)$ and $\rho(w)$, respectively.
Their lengths are
\[
\frac{\tau_u}{h_1},\qquad
\frac{-\tau_v}{h_0},\qquad
\frac{-\tau_w}{h_0},
\]
all smaller than $\varepsilon$.
The tip of the slit and the two corresponding points on the vertical
sides are glued together to form a branch point of angle $2\pi+\pi+\pi=4\pi$.
All boundary points are regular.

\begin{figure}[h]
    \centering
    \includegraphics[width=1\textwidth]{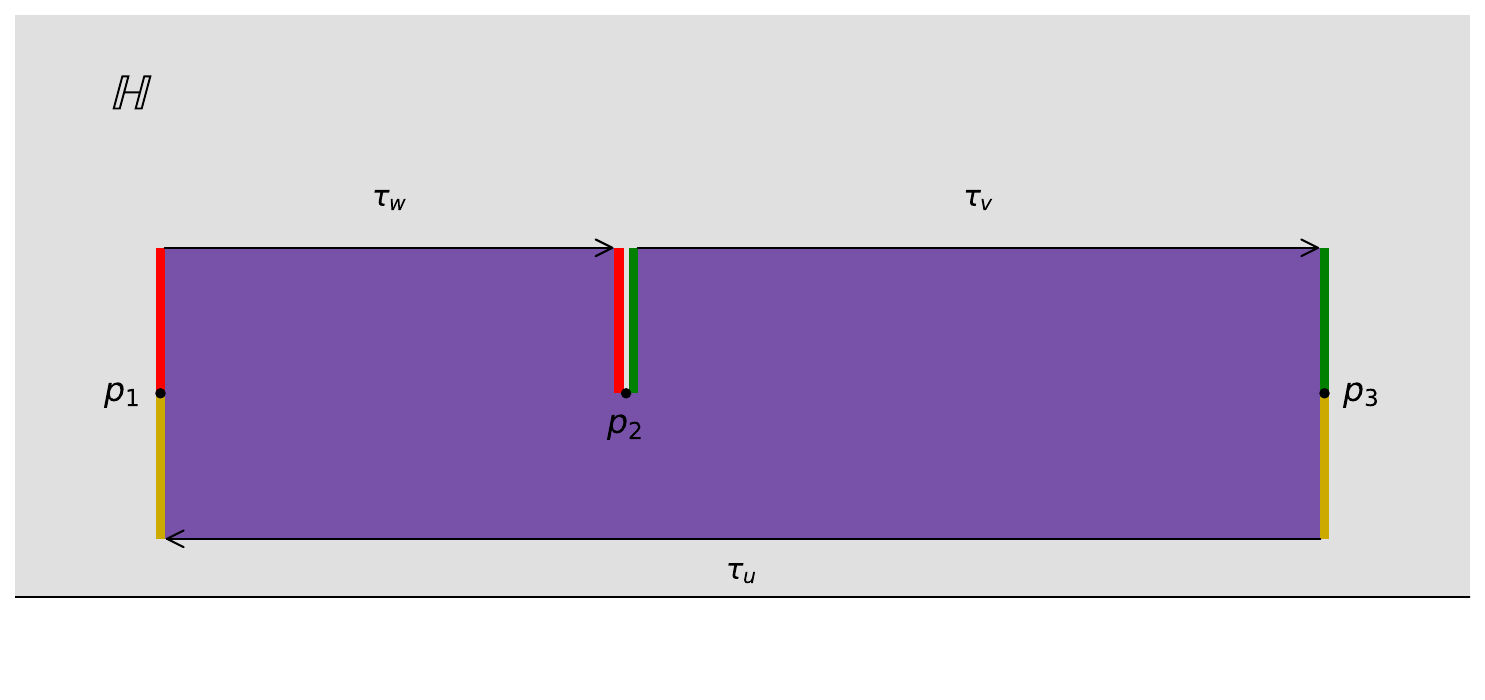}
    \caption{The symmetric configuration, with a slit emanating
    from the upper side, for  $\tau_u<0$ and $\tau_v,\tau_w>0$}.
    \label{f: reflected slit rectangle}
\end{figure}

The argument is identical in the remaining cases where a single one of the \(\tau_i\)'s is positive, since cyclically permuting the letters does not affect the reasoning. In the case where two of them are positive, one argues similarly by cutting a slit emanating from a point on the upper horizontal side of the rectangle, as in Figure~\ref{f: reflected slit rectangle}. \\

\section{Geometrizing representations of the fundamental group of a one-holed torus}

In this section, we prove a geometrization theorem for a large class of representations of the fundamental group of a one-holed torus, with control over the length of the boundary. Throughout this section, \(\Sigma_{1,1}\) denotes a surface of genus $1$ with $1$ boundary component. We write $\pi_1(\Sigma_{1,1})$ as the free group on two generators $a,b$ represented by simple closed curves in $\Sigma_{1,1}$ intersecting exactly once and such that the homological interesection number $a\wedge b$ is positive. Their commutator $c= [a,b]$ is then homotopic to the boundary curve.

Cutting $\Sigma_{1,1}$ along $a$, one obtains a pair of pants $\Sigma_{0,3}$. We can write $\pi_1(\Sigma_{0,3}) = \langle u,v,w\mid uvw= \Id\rangle$ as in the previous section. Choosing the generators $u,v,w$ and the lift $\tilde \Sigma_{0,3} \hookrightarrow \tilde \Sigma_{1,1}$ of the inclusion appropriately, the induced inclusion of fundamental groups is given by
\[u \mapsto  a^{-1},~v\mapsto c, w \mapsto bab^{-1}~.\]

\begin{figure}[h]
    \centering
\includegraphics[width=0.8\textwidth]{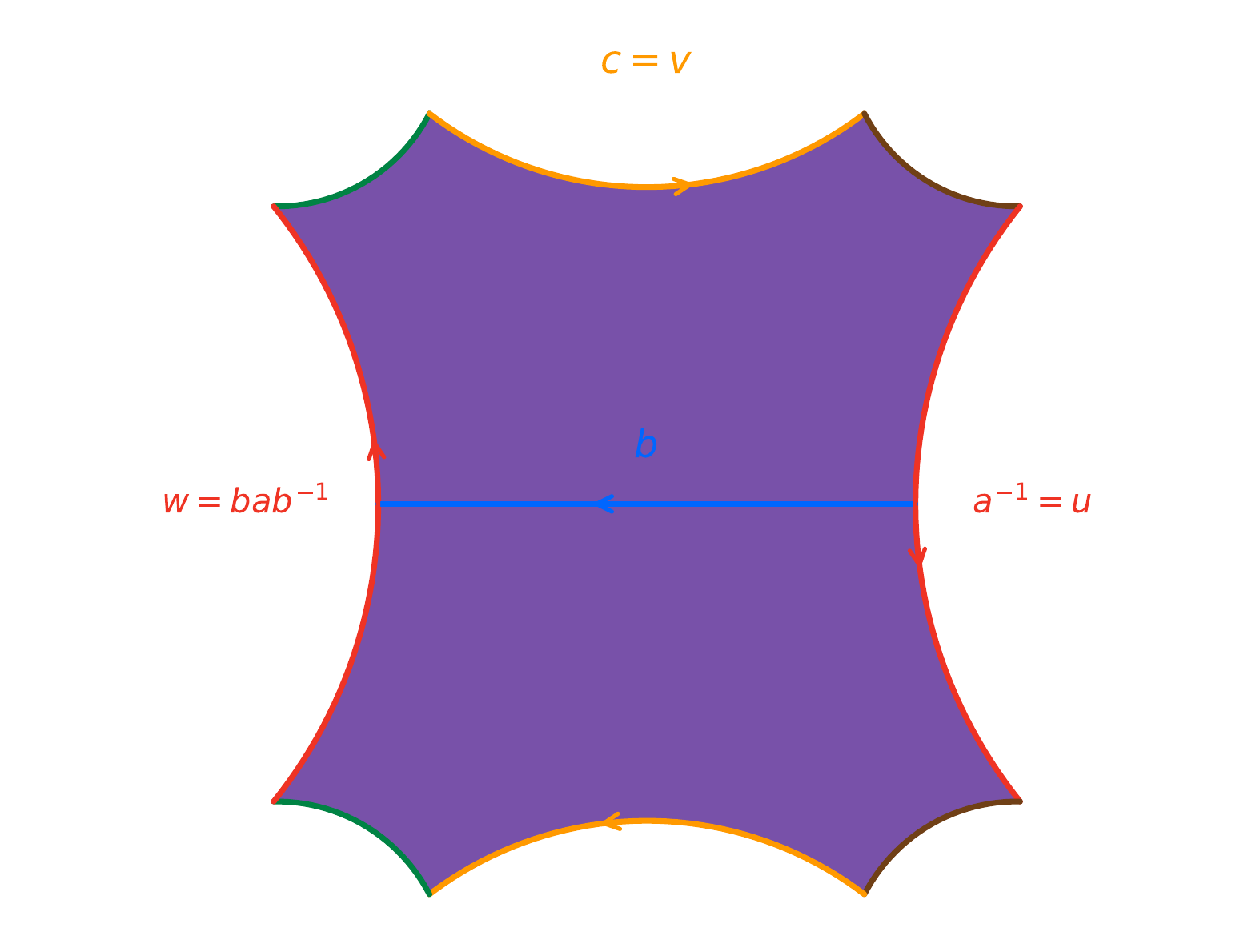}
    \caption{Gluing together the green and brown sides of this octagon produces a pair of pants $\Sigma_{0,3}$ with marked boundary curves $u,v,w$. Further gluing the red sides gives a one-holed torus with marked curves $a,b$. The region above the blue arc $b$ corresponds topologically to the plug in
Figure~\ref{f: domain for hyperbolic block Sigma_1_1}, used in the proof
of Lemma~\ref{l:Uniformization punctured torus elliptic}}.
    \label{f: Sigma_1_1}
\end{figure}

\subsection{Goldman's invariant}\label{ss: Goldman invariant}
Let $\rho: \pi_1(\Sigma_{1,1}) \to \Isom_+(\H) \simeq \PSL(2,\R)$ be a representation, and $\hat{\rho}: \pi_1(\Sigma_{1,1}) \to \SL(2,\R)$ be a lift of $\rho$. Then the number
\[\uptau(\rho) = \Tr(\hat{\rho}(c))\]
does not depend on the lift $\hat \rho$. We call it \emph{Goldman's invariant}, in reference to his investigation of the modular dynamics on the level sets of $\uptau$ in the character variety \cite{GoldmanPuncturedTorus}.

Note that conjugation by an orientation-reversing isometry of $\H$ preserves $\uptau(\rho)$. Goldman's invariant gives precise information on the geometric behaviour of the representation $\rho$. More precisely, we have the following relations (see \cite{GoldmanPuncturedTorus}):

\begin{itemize}
\item If $\uptau(\rho) < -2$, then $\Vol(\rho) = \pm 2\pi$ and up to conjugation by an orientation-reversing isometry, $\rho$ is the holonomy of a smooth hyperbolic structure on $\Sigma_{1,1}$ with totally geodesic boundary.

\item If $\uptau(\rho) = -2$, then $\Vol(\rho) = \pm 2\pi$. Up to conjugation by an orientation-reversing isometry, $\rho$ is the holonomy of a hyperbolic structure on $\Sigma_{1,1}^\bullet$ with a cusp at $c^\bullet$.

\item If $-2 < \uptau(\rho) < 2$, then $\Vol(\rho) = \pm (2\pi - \theta)$, where $\theta = 2\arccos(-\uptau(\rho)/2)$. Up to conjugation by an orientation-reversing isometry, $\rho$ is the holonomy of a hyperbolic structure on $\Sigma_{1,1}^\bullet$ with a cone singularity of angle $\theta$ at $c^\bullet$.

\item If $\uptau(\rho)=2$, then $\Vol(\rho)=0$ and $\rho$ has solvable image.

\item If $\uptau(\rho)>2$, then $\Vol(\rho)=0$.
\end{itemize}

In the range $2<\uptau(\rho)<18$, Goldman proves the following:

\begin{thm}[Goldman, {\cite[Theorem 5.2.1]{GoldmanPuncturedTorus}}] \label{thm: Goldman elliptic curve}
   If $2<\uptau(\rho)< 18$, then there exists a non-separating simple closed curve $a$ such that $\rho(a)$ is elliptic.
\end{thm}

\begin{rmk}
    Goldman actually shows the existence of a simple closed curve with elliptic or parabolic image, but one can verify that his trace-reduction algorithm proves this slightly stronger result.
\end{rmk}

\subsection{Geometrizing representations with an elliptic simple closed curve}

The goal of this section is to prove that representations of the one-holed in the range $2< \uptau < 18$ are still holonomies of a hyperbolic block and to control the boundary length of this hyperbolic block. Let $a$ be a simple closed curve in $\Sigma_{1,1}$ which is not homotopic to the boundary. We prove the following:

\begin{lem} \label{l:Uniformization punctured torus elliptic}
\label{l: geometrization pairs of pants in the elliptic case}

Let \(\rho: \pi_ 1 (\Sigma_{1,1}) \rightarrow \Isom_+(\H)\) be a representation with non-abelian image such that $\rho(a)$ is elliptic. Then, for all \(l> \ell(\rho(c))  \cosh(\ell(\rho(c))/4)\), there exists a hyperbolic block structure \(\sigma\) on \(\Sigma_{1,1}\) with holonomy \(\rho\) and boundary length less than \(l\).
\end{lem}

\begin{rmk}
    In the proof we will see that, under these assumptions, $\rho(c)$ is always hyperbolic.
\end{rmk}

\begin{rmk}
    Since all simple closed curves non-homotopic to $c$ are equivalent modulo the action of the mapping class group, Lemma \ref{l:Uniformization punctured torus elliptic} implies that a non-abelian representation is the holonomy of a hyperbolic block as soon as it sends one of these simple curves to an elliptic element. Combined with  Theorem~\ref{thm: Goldman elliptic curve}, we get that any representation $\rho$ with $2 < \uptau(\rho) < 18$ is the holonomy of a hyperbolic block, whose boundary has small length when $\uptau(\rho)$ is close to $2$.
\end{rmk}

\begin{proof}
Up to simultaneously replacing $a$ and $b$ by their inverses, we may assume that $\rho(a)$ is a rotation of angle $2\theta$ with $0<\theta\leq\pi/2$.
We denote by \( \rho ': \pi_1 (\Sigma_{0,3})\rightarrow \Isom_+(\H) \) the restriction of \(\rho\) to \(\pi_1(\Sigma_{0,3})\); we thus have \( \rho'(u) = \rho(a^{-1}) \), \( \rho'( v) = \rho([a,b])\), and \(\rho'(w) = \rho (bab^{-1})\).

\begin{figure}[htbp]
    \centering
    \includegraphics[width=0.8\textwidth]{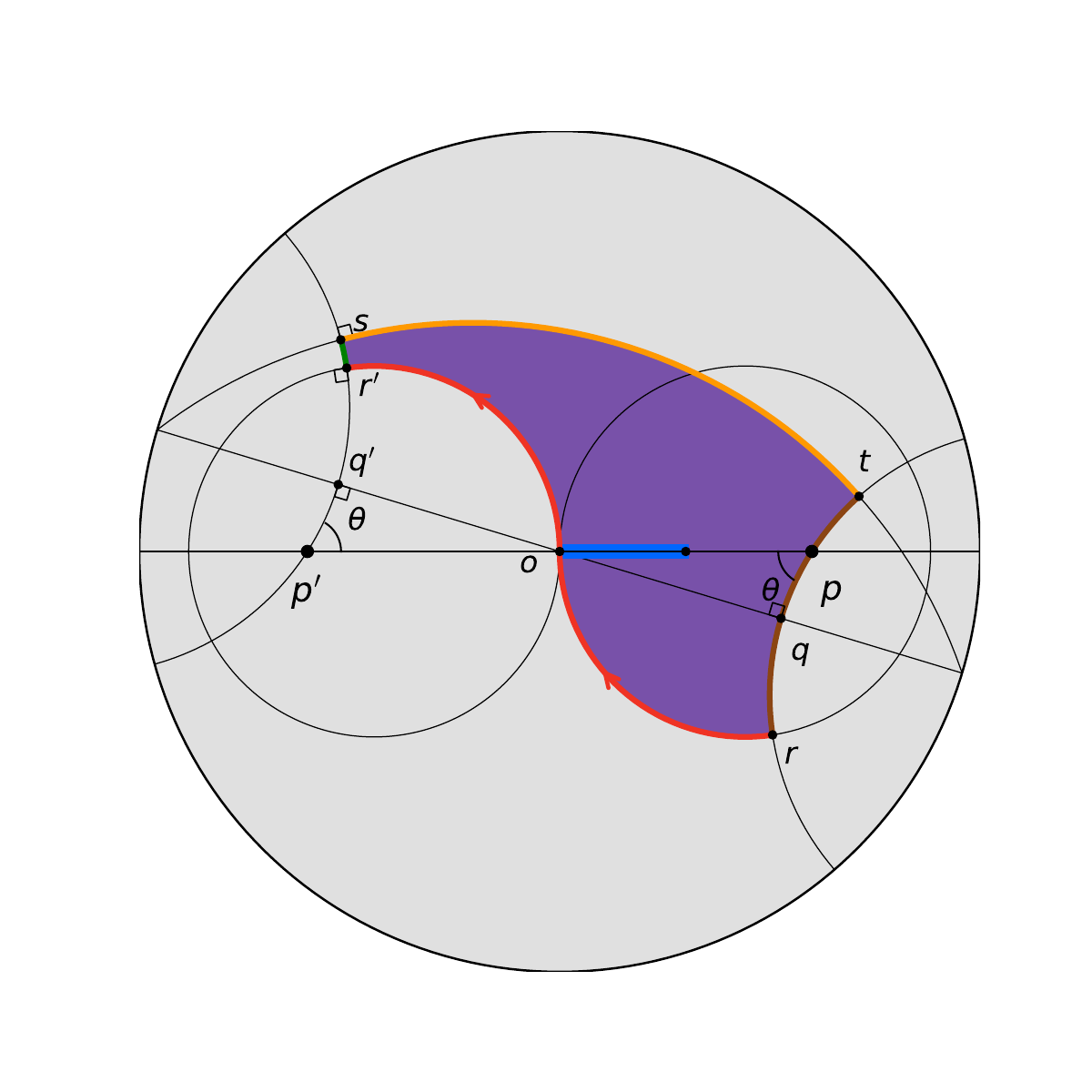}
    \caption{The plug: completion of a slit domain in the hyperbolic plane (unit disc model) with constant curvature sides.}
    \label{f: domain for hyperbolic block Sigma_1_1}
\end{figure}

Since \(\rho(a)\) is a non-trivial elliptic element, it has a unique fixed point \(p\in \mathbb H \). We define \(p'\coloneqq \rho(b) p\), the fixed point of the transformation \(\rho(b) \rho(a)\rho(b^{-1})\). 
The points \(p\) and \(p'\) are distinct. Indeed, otherwise, the maps \( \rho (a) \) and \(\rho(b) \rho(a)\rho(b^{-1})\)  would be two rotations with equal angles centreed at the same point, so they would be equal. But this would imply that \(\rho(a)\) and \(\rho(b)\) commute, and this is impossible since we assume that \(\rho\) has non-abelian image.

Let us describe Figure \ref{f: domain for hyperbolic block Sigma_1_1} in words. Let \(o\) be the midpoint of the segment \([p,p']\), and \(C\) (resp. \( C'\) ) the circle centred at \(p\) (resp. \(p'\) ) passing through \(o\). The geodesic ray based at \( p \) and forming an oriented angle \(\theta\) with the ray \( [p, p') \) intersects \(C\) at a point~\(r\). We let \(r'\) be the image of \(r\) under the half-turn about \(o\), so that the geodesic rays $[p',r')$ and $[p',p)$ form an angle $\theta$. We choose a point \(s\) slightly further than \(r'\) on the ray $[p',r')$.  The geodesics  \( (p,r) \) and \( (p',r') \) do not intersect at any point of \(\mathbb H \cup \partial_\infty \mathbb H\). Indeed, if this were the case, the intersection point would be the unique fixed point of the central symmetry with centre~\(o\) (since this symmetry exchanges the geodesics~\((p,r)\) and~\((p',r')\)), but~\(o\) does not belong to either of them. In particular, there is a unique geodesic perpendicular to both~\( (p,r) \) and~\( (p',r')\) at points denoted respectively~\(q\) and~\(q'\). Finally, we draw the curve passing through $s$ consisting of points whose distance to the geodesic \( (q,q') \) is equal to \(d=d(q',s)\); it intersects the geodesic \((p,r)\) at a point \(t\).  We consider the interior region of the closed curve defined by the concatenation of 
\begin{itemize} 
\item the arc between \(o\) and \(r\) in \(C\) (with the counterclockwise orientation on \(C\)),
\item the geodesic segment \([r,t]\), 
\item the arc between \(t\) and \(s\) of the curve at distance \(d\) from \((q,q')\), 
\item the geodesic segment \( [s,r']\), 

\item the arc between \(r'\) and \(o\) in \(C'\) (with the clockwise orientation on \(C'\)).
\end{itemize}
Finally, we cut a slit of this domain at the point \(o\), along an interval strictly contained in the geodesic segment \([o,p]\). We call the resulting hyperbolic surface (depicted in Figure~\ref{f: domain for hyperbolic block Sigma_1_1}) the \emph{plug}.

Let us now consider two copies of the plug, the first with its natural orientation, and the second with the opposite one. We glue these two copies by identifying the pairs of segments of the same color on each plug, as in Proposition~\ref{p: paste}.  The result is a hyperbolic block homeomorphic to a sphere minus three discs. After a suitable identification with \(\Sigma_{0,3}\) by an orientation-preserving diffeomorphism, this gives a hyperbolic block structure $\sigma'$ on \(\Sigma_{0,3}\) such that:
\begin{itemize}
    \item the curve $u$ is mapped to the double of the arc from $r$ to $o$ on $C$;
    \item the curve $v$ is mapped to the double of the arc from $s$ to $t$ at distance $d$ from $(q,q')$;
    \item the curve $w$ is mapped to the double of the arc $(r',o)$ on $C'$.
\end{itemize}
This identification with $\Sigma_{0,3}$ is illustrated by
Figure~\ref{f: Sigma_1_1}: the two copies of the plug in
Figure~\ref{f: domain for hyperbolic block Sigma_1_1} correspond,
respectively, to the closed regions above and below the blue arc $b$,
with boundary arcs matched by color. In each copy, the two slides
of the slit taken together correspond to the blue arc $b$.
The holonomy representation of $\sigma'$ then maps $u$ to the rotation of angle $-2\theta$ around $p$, i.e. $\rho(a)^{-1}$, and $w$ to the rotation of angle $2\theta$ around $p'$, i.e. $\rho(bab^{-1})$. Since $u$ and $w$ generate $\pi_1(\Sigma_{0,3})$, we conclude that $\sigma'$ has holonomy $\rho'$.

Finally, the curves $u$ and $w$ have the same length and opposite curvature for $\sigma'$. By Lemma~\ref{l: gluing lemma}, these boundaries can thus be glued together by an orientation-reversing isometry to form a hyperbolic block $\sigma$ homeomorphic to $\Sigma_{1,1}$ with holonomy $\rho$.

It remains to estimate the boundary length of this hyperbolic block. It is equal to twice the length of the arc between \(s\) and \(t\) on the curve at distance \(d\) from \((q,q')\), which by classical hyperbolic geometry gives:
\begin{equation} \label{eq: length boundary punctured torus}
l( \partial \Sigma_{1,1} ) = 2 d(q,q') \cosh(d(q',s))~.
\end{equation}

Denote by \( \sigma_l\) the symmetry with respect to a geodesic \(l\) of the hyperbolic plane. We have
\[ \rho (a) = \sigma_{(pq)} \circ \sigma_{(pp')} \text{ and } \rho (ba^{-1}b^{-1}) = \sigma_{(pp')} \circ \sigma_{(p' q')},\]
hence \( \rho ([a,b]) = \sigma_{(pq)} \circ \sigma_{(p' q')} \). Therefore, \(\rho ([a,b])\) is the hyperbolic isometry with axis \( (q' q ) \) and translation length
\[\ell(\rho(c)) = 2 d(q,q')~.\]

We have
\begin{eqnarray*}
    d(q',r')+d(p',q') &= & d(p',r') = d(p',o)\\
    &\leq & d(p',q') +d(q',o)~,
\end{eqnarray*}
from which we obtain
\[d(q',r')  \leq d(q',o) = \frac{1}{4} \ell(\rho(c))~.\]
Since $s$ can be chosen arbitrarily close to $r'$ and $l>\ell(\rho(c))\cosh(\ell(\rho(c))/4)$, we can choose $s$ so that \eqref{eq: length boundary punctured torus} gives the required bound
\[ l( \partial \Sigma_{1,1} ) < l~. \]
\end{proof}

\subsection{Geometrizing reducible representations}

We will also need to geometrize some representations $\rho$ for which $\uptau(\rho) =2$. These are exactly the representations that are reducible over $\C$, i.e. those fixing a point in $\H \sqcup \partial_\infty \H$ (see \cite[Section 2.3]{GoldmanPuncturedTorus}).

\begin{lem}
\label{l: geometrization pairs of pants in the parabolic case}
Let \(\rho: \pi_ 1 (\Sigma_{1,1}) \rightarrow \Isom_+(\H)\) be a reducible representation such that $\rho(c)$ is a negative parabolic.

Then, for all \(\epsilon > 0\), there exists a hyperbolic block structure \(\sigma\) on \(\Sigma_{1,1}\) with holonomy \(\rho\) and boundary length less than \(\epsilon\).
\end{lem}

\begin{rmk}
    If $\rho$ is reducible, $\rho(c)$ has to be parabolic or the identity. If $\rho(c)$ is a positive parabolic or the identity, then $\rho$ cannot be the holonomy of a hyperbolic block. Indeed, this block would have concave boundary, and both the volume of the hyperbolic structure and the curvature of the boundary would contribute positively to the volume of $\rho$, which is $0$.
\end{rmk}

\begin{proof}
The representation $\rho$ is not abelian since $\rho([a,b]) \neq \Id$. Hence it is reducible over $\R$ and has a global fixed point  \(p\in \partial_{\infty} \mathbb H\). If every simple closed curve had parabolic image, then $\rho$ would have abelian image and we would get $\rho(c) = \Id$. Hence our assumption on $\rho(c)$ ensures that there exists a simple closed curve $a$ such that $\rho(a)$ is hyperbolic. Up to reversing the orientation of $a$, we can assume that $p$ is the repelling fixed point of $a$.

As before, we let \(b\in \pi_1 (\Sigma_{1,1}) \) be represented by a simple closed curve intersecting $a$ exactly once, and oriented so that the homological intersection $a\wedge b$ is positive. Up to conjugating $\rho$, we can assume that the axis of $\rho(a)$ is the geodesic $(0,\infty)$, i.e. 
\[\rho(a): z\mapsto \frac{1}{\lambda} z\]
for some $\lambda >1$. Since $\rho(a)$ and $\rho(b)$ do not commute, $\rho(b)$ does not fix $0$, and the axis of $\rho(bab^{-1})$ is the geodesic $(t,\infty)$ for some $t\neq 0$. The fact that $\rho(c)$ is a negative parabolic implies that $t>0$.

Fix parameters $\mu > \nu > 0$ and consider the polygonal domain in $\H$ with consecutive vertices 
\[
\mu i,\, \lambda \mu i,\, \lambda \mu i + t,\,
\mu i + t,\, \mu i + (1-\frac{1}{2\lambda})t,\,
\nu i + (1-\frac{1}{2\lambda})t,\,
\nu i + \frac{1}{2\lambda}t,\,
\mu i + \frac{1}{2\lambda}t{\color{blue} .}
\]
(the purple domain in Figure \ref{f: torus with solvable monodromy}).
\begin{figure}[htbp]
    \centering
    \includegraphics[width=0.6\textwidth]{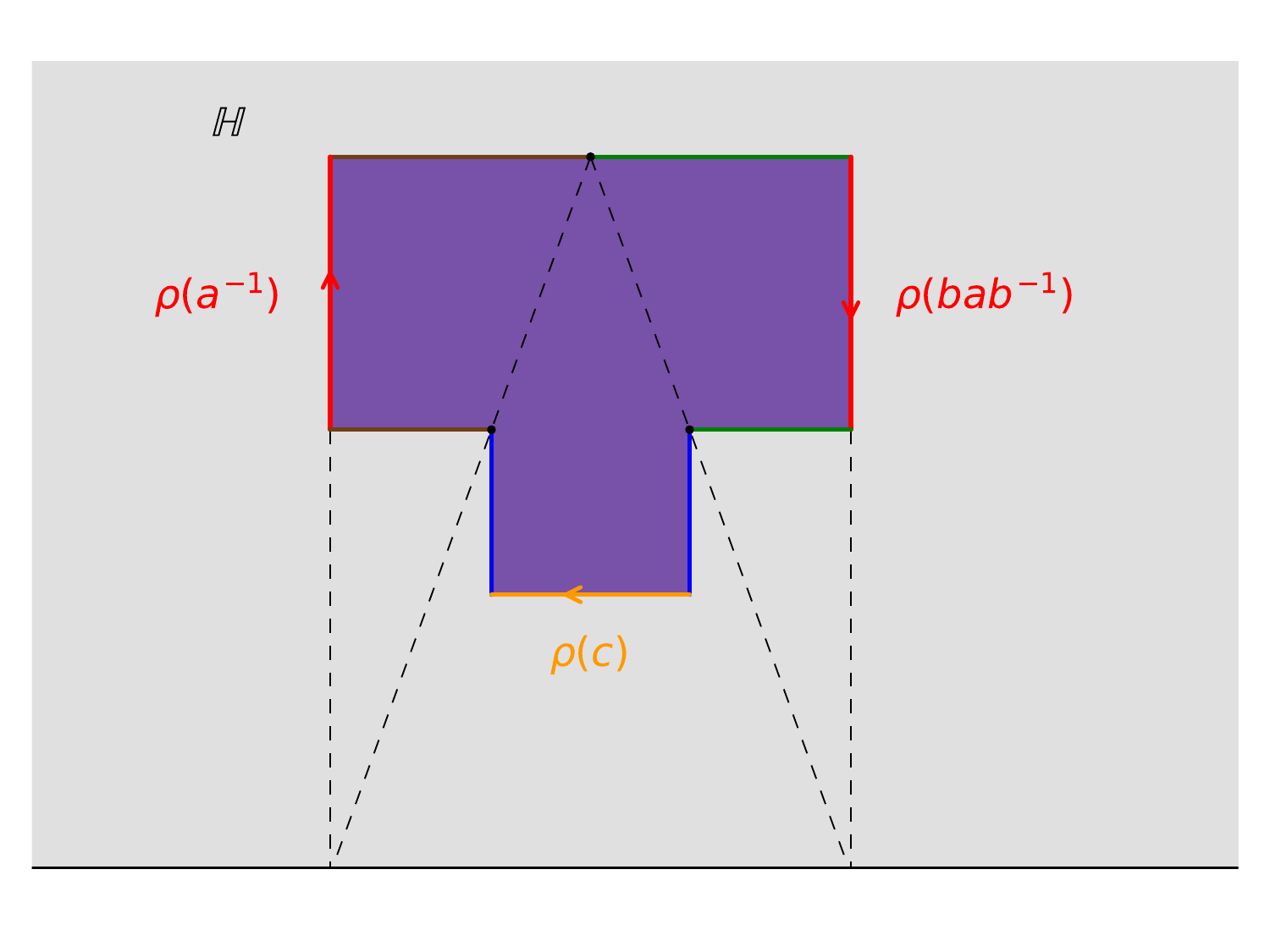}
    \caption{Case where \(\rho\) has reducible image.}
    \label{f: torus with solvable monodromy}
\end{figure}

Note that the yellow horocyclic arc $[\mu i, \mu i + \frac{1}{2\lambda}t ]$ is the image of the yellow horocyclic arc $[\lambda \mu i , \lambda \mu i + \frac{t}{2}]$ by $\rho(a)$. Similarly, the green horocyclic arcs (resp. the blue geodesic segments) are images of one another under $\rho(bab^{-1})$ (resp. $\rho(c): z\mapsto z- (1-\frac{1}{\lambda})t$). Hence, gluing the pairs of arcs of the same color as in Proposition~\ref{p: paste}, one obtains a hyperbolic block $\sigma'$ on $\Sigma_{0,3}$, whose holonomy representation maps $u$ to $\rho(a)^{-1}$, $v$ to $\rho(c)$ and $w$ to $\rho(bab^{-1})$. Hence $\sigma'$ has holonomy $\rho'$. Note that the points $\lambda \mu i + \frac{t}{2}$, $\mu i + \frac{t}{2\lambda}$ and $\mu i + (1-\frac{1}{2\lambda})t$ are identified in $\sigma'$ giving rise to a branch point of angle $4\pi$.

The boundary curves $u$ and $w$ are geodesics of the same length for $\sigma'$. By Lemma~\ref{l: gluing lemma}, they can just be glued together by an orientation-reversing isometry to obtain a hyperbolic block structure $\sigma$ on $\Sigma_{1,1}$ with holonomy $\rho$.
The length of $\partial \Sigma_{1,1}$ for $\sigma$ is the length of the horocyclic arc between $\nu i + (1-\frac{1}{2\lambda})t$ and $
\nu i + \frac{1}{2\lambda}t$, i.e. $(1-\frac{1}{\lambda})\frac{t}{\nu}$. This length can be chosen arbitrarily small if $\nu$ is chosen sufficiently large.
\end{proof}

\section{Representations of co-Euler class 1}

We now piece together the arguments developed in the preceding sections to prove the Main Theorem.

\subsection{Classification of the possible degeneracies}

We now fix a closed surface $S$ of genus $g \geq 2$ and $\rho: \pi_1(S) \to \Isom_+(\H)$ a representation of volume $2\pi(2g-3)$. Choosing a minimizing sequence for $\EE_\rho$ and applying Lemma \ref{lem: Conformal limits harmonic maps minimizing sequence}, we can find a decomposition
\[S= \bigsqcup_{1\leq i \leq I} S_i \sqcup \bigsqcup_{1\leq l \leq L} b_l~,\]
a sequence $(\tau_n)_{n\in \N}$ of elements of $\MCG(S,*)$, sequences $(g^i_n)_{1\leq i \leq I, n\in \N}$ in $\Isom_+(\H)$ and a collection $(j_i)_{1 \leq i \leq I}$ of finite covolume Fuchsian representations of $\pi_1(S_i)$ such that, for each $i$,
\begin{itemize}
\item $\rho^i_n \coloneqq {\Ad_{g^i_n} \circ (\rho \circ \tau_n)}_{|\pi_1(S_i)}$
converges to a representation
\[\rho^i_\infty: \pi_1(S_i) \to \Isom_+(\H)~,\]
\item There exists a $(j_i, \rho_i^\infty)$-equivariant holomorphic or antiholomorphic map $f_i: \H\to \H$.
\end{itemize}

Starting from any representation $\rho$, some of the $\rho_i^\infty$ could have negative volume (in which case the map $f_i$ is antiholomorphic) or vanishing volume (in which case the map $f_i$ is constant).

However, using the fact that $\Vol(\rho)$ is almost maximal, we can prove the following:

\begin{lem} \label{lem: decomposition co-euler class 1}
    We are in one of the following cases: 
    \begin{enumerate}
    \item Each $f_i$ is holomorphic and non-constant,
    
    \item Up to relabelling the decomposition, $\Vol(\rho_\infty^1) =0$, $\rho_\infty^i$ is Fuchsian for all $i\geq 2$, $\chi(S_1) =-1$ and $S_1$ is contained in a $1$-holed torus with boundary $b_1^+$,
    
    \item Up to relabelling the decomposition, $\Vol(\rho_\infty^1) =0$, $\rho_\infty^i$ is Fuchsian for all $i\geq 2$, $\chi(S_1) =-1$ and $S_1$ is a pair of pants with boundaries $b_1^+, b_2^+$, $b_3^+$.
    \end{enumerate}
\end{lem}

\begin{rmk}
    Condition (2) accounts for two similar but distinct situations: either $S_1= T$ or $T$ contains a curve $b_l$, $2\leq l \leq L$ and $S_1$ is the pair of pants $T\setminus b_l$. In any case we find a one-holed subtorus $T$ on which $\Vol({\rho\circ \tau_n}_{\vert T}) \underset{n\to +\infty}{\longrightarrow} 0$.
\end{rmk}

\begin{proof}
    By Proposition \ref{prop: Additivity volume limit}, we know that
    \[\Vol(\rho)
    = \sum_{1\leq i \leq I} \Vol(\rho_\infty^i)~.\]
    Plugging in the additivity of the Euler characteristic, we deduce
    \[
        2\pi = 2\pi \vert \chi(S)\vert - \Vol(\rho)
        = \sum_{1\leq i \leq I}
        \bigl(
        2\pi \vert \chi(S_i)\vert
        - \Vol(\rho_\infty^i)
       \bigr)~.
    \]
    Note that each term of the sum is nonnegative by the
    Milnor--Wood inequality. Moreover, we have
    \begin{itemize}
        \item $2\pi \vert \chi(S_i)\vert
        - \Vol(\rho^i_\infty)
        = 2\pi \vert \chi(S_i)\vert$
        if $\Vol(\rho_\infty^i)=0$,
        \item $2\pi \vert \chi(S_i)\vert
        - \Vol(\rho_\infty^i)>
        2\pi \vert \chi(S_i)\vert$
        if $\Vol(\rho_\infty^i)<0$.
    \end{itemize}
    Since $\vert \chi(S_i)\vert \geq 1$ for each $i$,
    we deduce that:
    \begin{itemize}
        \item Either all the $\rho_\infty^i$ have positive volume,
        in which case the $f_i$ are non-constant and holomorphic
        and we are in case (1),
        \item Or, after relabelling the decomposition,
        we have $\Vol(\rho^1_\infty)=0$
        and $\vert\chi(S_1)\vert=1$, and for all $i\geq2$, $\Vol(\rho_\infty^i)=2\pi\vert\chi(S_i)\vert$.
    \end{itemize}

    In the latter case, for each $i\geq2$, the map $f_i$ is
    holomorphic and the Schwarz--Pick lemma gives
    $f_i^*\sigma_{\H}\leq \sigma_{\H}$.
    The corresponding area forms descend to
    $j_i(\pi_1(S_i))\backslash\H$, and their integrals are equal,
    since
    $\Vol(\rho_\infty^i)=2\pi|\chi(S_i)|$.
    Hence equality holds everywhere, so $f_i$ is an
    orientation-preserving isometry of $\H$.
    By equivariance, $\rho_\infty^i$ is conjugate to $j_i$.
    In particular, it is Fuchsian of finite covolume, with
    non-trivial parabolic peripheral holonomies.

     Finally, if $S_1$ is a one-holed torus or a pair
    of pants of which two boundary components are glued together
    in $S$, then we are in case (2). Otherwise we are in case (3).
\end{proof}

We will conclude the proof of the Main Theorem by treating the three cases of Lemma \ref{lem: decomposition co-euler class 1} separately. The approach, however, is always the same: using the limiting representations $\rho_\infty^i$ together with the Ehresmann--Thurston principle (Lemma \ref{l:perturbation of hyperbolic block}), we prove that $\rho$ restricted to each piece of the decomposition is the holonomy of a hyperbolic block and, controlling the boundary behaviour, we verify that we can apply Lemma \ref{l: gluing lemma} in order to glue the pieces together.

\subsection{The case where nothing collapsed}

Here we conclude the proof of the Main Theorem in case (1) of Lemma \ref{lem: decomposition co-euler class 1}, namely, when for each $1 \leq i \leq I$, there exists a Fuchsian representation $j_i:\pi_1(S_i)\to \Isom_+(\H)$ and a $(j_i, \rho_\infty^i)$-equivariant holomorphic non-constant map $f_i:\H \to \H$.

By Proposition \ref{prop: holomorphic map -> singular metric}, $\rho^i_\infty$ is the holonomy of a singular hyperbolic metric on $S_i^\bullet$ with cusps and cone singularities of angle $\alpha({b_l^\pm})$ at the marked points ${b_l^{\pm}}^\bullet$. By Proposition \ref{prop: Cone metrics -> hyperbolic blocks}, we deduce the existence of some $\epsilon_0>0$ such that $\rho^i_\infty$ is the holonomy of a hyperbolic block with concave boundary lengths $2\epsilon_0$ on $S^\sharp$.

By our Ehresmann--Thurston principle for hyperbolic blocks (Lemma \ref{l:perturbation of hyperbolic block}), we obtain that, for $n$ large enough, $(\rho^i_n)_{1\leq i \leq I}$ (hence ${\rho\circ \tau_n}_{\vert S^\sharp}$) is the holonomy of a hyperbolic block $\sigma_n^\sharp$ on $S^\sharp$ such that, for every boundary component $b = b_l^\pm$ of $S^\sharp$:
\begin{itemize}
    \item either $\rho_\infty(b) = \Id$ and $\ell_{\sigma_n^\sharp}(b) < \epsilon_0$    
    \item or $\rho_\infty(b) \neq \Id$, $b$ is concave for $\sigma_n^\sharp$ and $\ell_{\sigma_n^\sharp}(b) \geq \epsilon_0$.
\end{itemize}

By the Gauss--Bonnet formula, we have
\[\Vol(\rho_\infty) = 2\pi \vert \chi(S^\sharp)\vert - 2\pi\vert B \vert - \left(\sum_{1\leq l \leq L} \alpha(b_l^+) + \alpha(b_l^-)\right)~,\]
where $\vert B\vert $ is the number of branch points counted with multiplicity. Hence 
\[2\pi\vert B \vert + \left(\sum_{1\leq l \leq L} \alpha(b_l^+) + \alpha(b_l^-)\right) = 2\pi~.\]
In particular, since all these terms are nonnegative, there is at most one boundary curve $b$ of $S^\sharp$ such that $\alpha(b) = 2\pi$, hence at most one curve $b$ such that $\rho^\infty(b) = \Id$. In particular, for every $l$ we have $\rho_n(b_l) \neq \Id$ for $n$ large enough (otherwise we would get $\rho-\infty(b_l^+)= \rho_\infty(b_l^-) = \Id$. We thus get that, for all $1\leq l \leq L$,
\begin{itemize}
    \item Either $b_l^+$ and $b_l^-$ are both concave for $\sigma_n^\sharp$
    \item Or (up to reversing the orientation of $b_l$), $\rho_\infty(b_l^+) = \Id$ and $b_l^+$ might be convex for $\sigma_n^\sharp$, but $b_l^-$ is concave and 
    \[\ell_{\sigma_n^\sharp}(b_+) \leq \epsilon_0 \leq \ell_{\sigma_n^\sharp}(b_-)~.\]
\end{itemize}
We can thus apply our gluing lemma for hyperbolic blocks (Lemma \ref{l: gluing lemma}) and deduce that $\rho\circ \tau_n$ is the holonomy of a branched hyperbolic metric on $S$. Hence so is $\rho$.

\subsection{The case where a handle collapsed}

Here we conclude the proof of the Main Theorem in case (2) of Lemma \ref{lem: decomposition co-euler class 1}, namely, when
there exists a partition $S = \bigsqcup_{ 1\leq i\leq I} S_i \sqcup \bigsqcup_{1\leq l \leq L} b_l$ such that:
\begin{itemize}
    \item $S_1$ is a $1$-holed torus with boundary $b_1^+$ and $\Vol(\rho_n^1)\underset{n\to +\infty}{\longrightarrow} 0$,
    \item $\rho_n^i$ converges to a Fuchsian representation $\rho_\infty^i:\pi_1(S_i) \to \Isom_+(\H)$ for all $i\geq 2$.\\
\end{itemize} 

For all $i\geq 2$, there exists $\epsilon>0$ such that $\rho^i_\infty$ is the holonomy of a hyperbolic block with concave boundary of length $2\epsilon$. Moreover, since $\rho^i_\infty(b)\neq \Id$ for every boundary curve $b$, applying the Ehresmann--Thurston principle for hyperbolic blocks, we get that, for $n$ large enough, $\rho_n^i$ is the holonomy of a hyperbolic block $\sigma_n^\sharp$ on $\bigcup_{2\leq i \leq I} S_i^\sharp$ with concave boundaries of length $\geq \epsilon$.

Let us now turn to $S_1$. 
Since $\rho^1_n(b_1^+)$ is conjugate to the inverse of $\rho_n(b_1^-)$ which converges to a positive parabolic, we get that
\[\vert \Tr( \rho^1_n(b_1^+))\vert \underset{n\to +\infty}{\longrightarrow} 2~.\]
Since $\Vol(\rho_n^1) \underset{n\to +\infty}{\longrightarrow} 0$, we deduce that
\[\uptau(\rho_n^1) \underset{n\to +\infty}{\longrightarrow} 2\]
by the classification in Section~\ref{ss: Goldman invariant}. Finally, since
$2\pi - \Vol(\rho^1_n) \leq \sum_{1\leq i \leq I} 2\pi \vert \chi(S_i)\vert - \Vol(\rho^i_n) = 2\pi \vert \chi(S)\vert - \Vol(\rho) = 2\pi~,$
we get that
\[\Vol(\rho^1_n) \geq 0~.\]

Fix a small enough $\eta>0$. For $n$ large enough, we are in one of the following situations:
\begin{itemize}
\item[(a)] Either $2-\eta < \uptau(\rho_n^1) < 2$, in which case $\rho_n^1$ is the holonomy of a hyperbolic structure with a cone singularity on $S_1^\bullet$ by the classification in Section~\ref{ss: Goldman invariant}.

\item[(b)] Or $\uptau(\rho_n^1) = 2$, hence $\rho_n^1$ has solvable image, and $\rho_n^1(b_1^+)$ is a negative parabolic (since $\rho_n(b_1^-)$ is conjugated to $(\rho_n^1(b_1^+))^{-1}$ and converges to a positive parabolic).

\item[(c)] Or $2< \uptau(\rho_n^1) < 2+\eta$, in which case there exists a non-separating simple closed curve $a$ such that $\rho_n^1(a)$ is elliptic.
\end{itemize}

In Case (a), we get that $\rho_n^1$ is the holonomy of a hyperbolic block with arbitrarily small boundary length by Proposition~\ref{prop: Cone metrics -> hyperbolic blocks}. In Case (b), we get the same conclusion by Lemma~\ref{l: geometrization pairs of pants in the parabolic case}. Finally, in Case (c), we get that $\rho_n^1$ is the holonomy of a hyperbolic block with boundary length $\leq \epsilon$ for $\eta$ small enough by Lemma~\ref{l:Uniformization punctured torus elliptic}.

Since, on the other side, the curve $b_1^-$ is concave of length $\geq \epsilon$ for $\sigma_n$ (as well as all the $b_l^\pm$, $l\geq2$), we can apply our gluing lemma \ref{l: gluing lemma} and conclude that $\rho \circ \tau_n$ -- hence $\rho$ -- is the holonomy of a branched hyperbolic structure on $S$.

\subsection{The case where a pair of pants collapsed}

Here we conclude the proof of the Main Theorem in Case (3) of Lemma \ref{lem: decomposition co-euler class 1}, namely, when there exists a partition $S = \bigsqcup_{1\leq i\leq I} S_i \sqcup \bigsqcup_{1\leq l \leq L} b_l$ such that:
\begin{itemize}
    \item $S_1$ is a pair of pants with boundaries $b_1^+, b_2^+, b_3^+$ and $f_n^1$ converges to a constant map,
    \item $\rho_n^i$ converges to a Fuchsian representation $\rho_\infty^i:\pi_1(S_i) \to \Isom_+(\H)$ for all $i\geq 2$.\\
\end{itemize} 

The proof is very similar to the proof of Case (2):
since $\rho^i_\infty$ is Fuchsian for all $i\geq 2$, applying the Ehresmann--Thurston principle for hyperbolic blocks, we obtain, for $n$ large enough, a hyperbolic block $\sigma_n^\sharp$ on $\bigsqcup_{2\leq i \leq I} S_i^\sharp$ of holonomy $(\rho_n^i)_{2\leq i \leq I}$, with concave boundaries of length $\geq \epsilon$.

Turning to $S_1$, we note that $\rho^1_n(b_1^+)$ is conjugate to $\rho_n(b_1^-)^{-1}$. Here $\rho_n(b_1^-)$ converges to a positive parabolic, so its inverse converges to a negative parabolic. Hence, for $n$ large enough, $\rho_n^1$ is not the trivial representation. 
Since $\rho_\infty^1$ fixes a point in $\H$, its minimal displacement vanishes. The sequence $(\rho^1_n)$ thus satisfies the hypotheses of Lemma~\ref{l: geometrization of representations with small translation length}, and we deduce that, for $n$ large enough, $\rho_n^1$ is the holonomy of a hyperbolic block $\sigma_n^1$ on $S_1^\sharp$ such that $\ell_{\sigma_n^1}(b_i^+) \leq \epsilon$ for $i=1,2,3$. Since, on the other side, the curves $b_i^-$ are concave of length $\geq \epsilon$ for $\sigma_n^\sharp$ (as well as all the $b_l^\pm$, $l\geq4$), we can apply our gluing lemma~\ref{l: gluing lemma} and conclude that $\rho \circ \tau_n$ -- hence $\rho$ -- is the holonomy of a branched hyperbolic structure on~$S$.

\bibliographystyle{alpha} 
\bibliography{bibli}

\end{document}